\documentclass[reqno, 11pt]{amsart}

\newtheorem{theorem}{Theorem}[section]
\newtheorem{lemma}[theorem]{Lemma}
\newtheorem{proposition}[theorem]{Proposition}
\newtheorem{corollary}[theorem]{Corollary}

\theoremstyle{definition}
\newtheorem{definition}[theorem]{Definition}
\newtheorem{ex}[theorem]{Example}

\theoremstyle{remark}
\newtheorem{remark}[theorem]{Remark}

\numberwithin{equation}{section}

\usepackage{bbm}
\usepackage{url}
\usepackage{euscript}
\usepackage{pb-diagram}
\usepackage{lamsarrow}
\usepackage{pb-lams}
\usepackage{amsmath}
\usepackage{amsthm}

\usepackage{epsfig}
\usepackage{amssymb}
\usepackage{pifont}
\usepackage{amsbsy}

\newskip\aline \newskip\halfaline
\def\skipaline{\vskip\aline}

\def\qedbox{$\rlap{$\sqcap$}\sqcup$}
\def\qed{\nobreak\hfill\penalty250 \hbox{}\nobreak\hfill\qedbox\skipaline}

\def\proofend{\eqno{\mbox{\qedbox}}}

\newcommand{\one}{\mathbbm{1}}

\newcommand\bC{{\mathbb C}}

\newcommand\bI{{\mathbb I}}

\newcommand{\bP}{{{\mathbb P}}}

\newcommand{\bQ}{{{\mathbb Q}}}
\newcommand\bR{{\mathbb R}}

\newcommand\bZ{{\mathbb Z}}

\DeclareMathOperator{\re}{\mathbf{Re}}
\DeclareMathOperator{\im}{\mathbf{Im}}

\DeclareMathOperator{\uu}{\underline{\mathit{u}}}

\DeclareMathOperator{\tr}{{\rm tr}}
\DeclareMathOperator{\supp}{{\rm supp}}

\DeclareMathOperator{\GL}{{\rm GL}}

\DeclareMathOperator{\codim}{codim}

\DeclareMathOperator{\Gr}{\mathbf{Gr}}
\DeclareMathOperator{\diag}{Diag} \DeclareMathOperator{\Hom}{Hom}
 \DeclareMathOperator{\End}{End}

\newcommand{\bc}{\boldsymbol{c}}
\newcommand{\be}{\boldsymbol{e}}
\newcommand{\bsf}{\boldsymbol{f}}
\newcommand{\bsF}{\boldsymbol{F}}

\newcommand{\ii}{\boldsymbol{i}}

\newcommand{\bp}{\boldsymbol{p}}

\newcommand{\bt}{\boldsymbol{t}}
\newcommand{\bu}{\boldsymbol{u}}
\newcommand{\bv}{\boldsymbol{v}}
\newcommand{\bw}{\boldsymbol{w}}
\newcommand{\bx}{\boldsymbol{x}}

\newcommand{\bsJ}{\boldsymbol{J}}

\newcommand{\balph}{\boldsymbol{\alpha}}
\newcommand{\bgamma}{\boldsymbol{\gamma}}

\newcommand{\bpi}{\boldsymbol{\pi}}

\newcommand{\bom}{\boldsymbol{\omega}}
\newcommand{\bOm}{\boldsymbol{\Omega}}

\newcommand{\si}{{\sigma}}

\newcommand{\ve}{{\varepsilon}}
\newcommand{\eps}{{\epsilon}}
\newcommand{\vfi}{{\varphi}}

\newcommand{\eA}{\EuScript{A}}
\newcommand{\eB}{\EuScript{B}}
\newcommand{\eC}{\EuScript{C}}

\newcommand{\eE}{\EuScript{E}}
\newcommand{\eF}{\EuScript{F}}

\newcommand{\eH}{\EuScript H}
\newcommand{\eI}{\EuScript{I}}

\newcommand{\eL}{\EuScript{L}}

\newcommand{\eN}{\EuScript{N}}
\newcommand{\eO}{\EuScript{O}}
\newcommand{\eP}{\EuScript{P}}
\newcommand{\eQ}{\EuScript{Q}}
\newcommand{\eR}{\EuScript{R}}
\newcommand{\eS}{\EuScript{S}}
\newcommand{\eT}{\EuScript{T}}
\newcommand{\eU}{\EuScript{U}}

\newcommand{\eW}{\EuScript{W}}
\newcommand{\eX}{\EuScript{X}}
\newcommand{\eY}{\EuScript{Y}}
\newcommand{\eZ}{\EuScript{Z}}

\newcommand{\ra}{\rightarrow}
\newcommand{\hra}{\hookrightarrow}
\newcommand{\Lra}{{\longrightarrow}}

\newcommand{\lan}{\langle}
\newcommand{\ran}{\rangle}

\def\inpr{\mathbin{\hbox to 6pt{\vrule height0.4pt width5pt depth0pt \kern-.4pt \vrule height6pt width0.4pt depth0pt\hss}}}

\newcommand{\pa}{\partial}

\newcommand{\ori}{\boldsymbol{or}}

\newcommand{\whA}{\widehat{A}}
\newcommand{\whE}{\widehat{E}}

\DeclareMathOperator{\St}{\mathbf{St}}

\DeclareMathOperator{\cl}{\boldsymbol{cl}}

\DeclareMathOperator{\Cr}{\mathbf{Cr}}
 
\DeclareMathOperator{\pfc}{\widehat{\EuScript{S}}_{an}}

\DeclareMathOperator{\Graff}{\mathbf{Graff}}

\DeclareMathOperator{\Lag}{Lag_h}
\DeclareMathOperator{\Symp}{Sp_h}
\DeclareMathOperator{\usymp}{\underline{sp}_h}

\DeclareMathOperator{\Fl}{\boldsymbol{Fl}}
\DeclareMathOperator{\Flag}{\mathbf{FLAG}}
\DeclareMathOperator{\Flagi}{\mathbf{FLAG}_{iso}}
\DeclareMathOperator{\Sp}{span}

\begin{document}

\title{Schubert calculus on the grassmannian of  hermitian lagrangian spaces}
\date{Started July 3, 2007. Completed August 15, 2007. 2nd Revision September 20, 2007.}

\author{Liviu I. Nicolaescu}

\address{Department of Mathematics, University of Notre Dame, Notre Dame, IN 46556-4618.}
\email{nicolaescu.1@nd.edu}

\begin{abstract}  We  describe a Schubert like stratification  on the Grassmannian of   hermitian lagrangian spaces in $\bC^n\oplus \bC^n$ which is a natural compactification of the  space of  hermitian  $n\times n$ matrices. The closures of the strata  define   integral  cycles and we  investigate  their  intersection theoretic properties. The  methods employed are Morse theoretic.

\end{abstract}

\maketitle

\tableofcontents

\section*{Introduction}
\addtocontents{section}{Introduction}

A hermitian lagrangian subspace  is a subspace $L$  of    the complex Hermitian vector space $\bC^{2n}=\bC^n\oplus \bC^n$ satisfying
\[
L^\perp =\bsJ L,
\]
where $\bsJ:\bC^n\oplus \bC^n\ra \bC^n\oplus \bC^n$ is the unitary operator with the block decomposition
\[
\bsJ=\left[
\begin{array}{cc}
0 &-\one_{\bC^n}\\
\one_{\bC^n} & 0
\end{array}
\right].
\]
We denote by $\Lag(n)$ the Grassmannian of such subspaces.     This space can be identified   with a more familiar space.

Denote by  $F^\pm\subset\bC^{2n}$ the $\pm \ii$ eigenspace  of $\bsJ$,
\[
F^\pm= \bigl\{ (\be, \mp\ii \be);\;\be\in \bC^n,\bigr\}.
\]
Arnold has shown  in \cite{Ar1} that $L\subset \bC^{2n}$ is a hermitian lagrangian subspace if and only  if, when viewed as a subspace of $F^+\oplus F^-$, it is the graph of a unitary operator $F^+ \ra F^-$. Thus we   have a natural diffeomorphism $U(n)\ra \Lag(n)$.

The unitary groups  are    some of the most investigated topological spaces and much is known  about their cohomology rings (see \cite[Chap.IV]{Ep}, \cite[VII.4, VIII.9]{Whit}), and one could fairly ask what else is there to say about these spaces.   To answer this, we need to   briefly explain the  question which gave the impetus for the investigations in this paper.

As is well known $U(\infty)$ is  a classifying space for the functor $K^1$, and its integral cohomology is an exterior algebra $\Lambda(x_1, x_2,\dotsc)$, where $\deg x_i=2i-1$.   If $X$ is a compact, oriented smooth manifold, $\dim X=n$, then the results of Atiyah and Singer \cite{AS} imply that any smooth family $(A_x)_{x\in X}$  Fredholm  selfadjoint  operators defines a  smooth map\footnote{We will not elaborate here on the precise  meaning of smoothness of $U(\infty)$.} $A: X\ra U(\infty)$. We thus obtain  cohomology classes $A^*x_i\in H^{2i-1}(X,\bZ)$.

We are interested in  geometric localization formul{\ae}, i.e.,  in describing \emph{concrete  geometric realizations} of cycles representing the    Poincar\'e duals of these classes. Some of the most interesting situations  arise when $X$ is an odd dimensional sphere $X=S^{2m-1}$.  In this case, the Poincar\'e dual of $A^*x_m$ is a    $0$-dimensional homology class, and we would like to produce an explicit  $0$-cycle  representing it.

For example, in the lowest dimensional case,  $X=S^1$,  we have such  a geometric  realization because the integer $\int_{S^1}A^*x_1$ is the spectral flow of the loop of selfadjoint operators, and as is well known,  in generic  cases,  this can be computed by counting   with appropriate multiplicities the points $\theta\in S^1$ where $\ker A_\theta=0$.  Thus, the  Poincar\'e  dual  of $A^*x_1$  is represented by a  certain $0$-dimensional  degeneracy locus.

The graph of  a  selfadjoint Fredholm operator $A: H\ra H$, $H$ complex Hilbert  space, defines a  hermitian lagrangian  $\Gamma_A$ in the hermitian symplectic space $H\oplus H$, and  we could  view a loop of such operators as a loop in $\Lag(\infty)$.    Adopting this point of view, we can interpret  the  integer $\int_{S^1} A^*x_1$ as a Maslov index, and  using the techniques  developed by   Arnold   in \cite{Ar0} one can explicitly  describe a $0$-cycle  dual to the class $A^*x_1$, \cite{N}.

To the best of our knowledge  there are no such   degeneracy loci descriptions  of the Poincar\'e dual of $A^*x_m$ in the higher dimensional cases $A: S^{2m-1}\ra U(\infty)$, $m>1$, and the existing  descriptions  of the cohomology ring of $U(n)$ do not seem to help in this respect.

With an eye towards  such  applications,  we   describe in this paper    a natural, Schubert like, Whitney regular, stratification of  $\Lag(n)$ and its intersection theoretic properties.

As in the case of   usual Grassmannians, this stratification has a Morse theoretic    description.    We  denote by $(e_i)$ the canonical unitary basis  of $\bC^n$,  and we define the Hermitian operator $A:\bC^{n}\ra \bC^{n}$ by setting
\[
Ae_i=\Bigl(i -\frac{1}{2}\,\Bigr)e_i,\;\;\forall i=1,\dotsc, n.
\]
The operator $A$ defines a function
\[
f=f_A: U(n)\ra \bR,\;\;f(S)=-\re\tr(AS)+\frac{n^2}{2}.
\]
This is a Morse function with one critical point $S_I\in U(n)$ for every subset $I\subset \{1,\dotsc, n\}$.  More precisely
\[
S_Ie_i=\begin{cases}
e_i & i\in I\\
-e_i &i\not\in I.
\end{cases}
\]
Its Morse index is ${\rm ind}\, (S_I)=f(S_I)=\sum_{i\in I^c}(2i-1)$, where $I^c$ denotes  the complement of $I$ in $\{1,2,\dotsc, n\}$. In particular,  this function is self-indexing.

We denote by $W_I^\pm$ the  stable/unstable manifold of $S_I$. These unstable manifolds  are loci of certain Schubert-like incidence relations and they can be  identified with the orbits of a real algebraic group acting  on $\Lag(n)$  so that, according to \cite{Land}, the stratification  given by  these unstable manifolds  satisfies the Whitney regularity condition. In particular, this implies that our  gradient flow satisfies the Morse-Smale  transversality condition. We can thus define the Morse-Floer complex, and it turns out  that the boundary operator of this complex  is trivial. The ideas outlined so far are  classical, going back to the pioneering work of Pontryagin \cite{Pont}, and we recommend \cite{DV} for a nice presentation.

Given that  the Morse-Floer complex is perfect  it is natural  to ask if  the unstable manifolds $W_I^-$  define  geometric cycles in any reasonable way, and if so,     investigate their intersection theory. M. Goresky \cite{Go} has explained how to associate cycles to Whitney stratified objects but this approach seems  difficult to use in concrete computations.

Another approach,  essentially used by  Vasiliev \cite{Vas}  is to   produce  resolutions of $W_I^-$, i.e., smooth maps $f:X_I\ra \Lag(n)$, where $X_I$ is a compact oriented manifold, $f(X_I)=\cl(W_I^-)$, and  $f$ is a diffeomorphism over the smooth part of $\cl(W_I^-)$.    As explained in \cite{Vas},  this approach    reduces the computation    of the intersection cycles $f_*[X_I]\bullet f_*[X_J]$ to classical Schubert calculus on Grassmannians, but  the combinatorial complexity  seemed very discouraging to this author.

 Instead, we chose the most obvious approach, and  we   looked at   the  integration currents    defined by the  semialgebraic sets $W_I^-$ as defining a cycle.   This is where   the theory of intersection of subanalytic cycles developed by R. Hardt \cite{Hardt, Hardt1, Hardt2} comes in very handy.

 The   manifolds $W_I^-$ are semi-algebraic,  have  finite volume, and carry natural orientations $\ori_I$, and thus  define integration currents $[W_I,\ori_I]$.    In Proposition \ref{prop: cycle} we show that the closure of $W_I^-$ is  a naturally oriented pseudo-manifold, i.e., it admits a stratification  by smooth manifolds, with top stratum oriented, while the other strata have (relative) codimension at least $2$. Using the fact that the current $[W_I^-,\ori_I]$ is a subanalytic current  as defined in \cite{Hardt2},  it follows that $\pa[W_I^-,\ori_I]=0$ in the sense of currents. We thus get cycles  $\balph_I\in H_\bullet(U(n),\bZ)$.

The   currents  $[W_I^-,\ori_I]$ define a perfect subcomplex  of the complex of integrally flat currents. This subcomplex is  isomorphic to the Morse-Floer complex, and via the finite-volume-flow technique  of Harvey-Lawson \cite{HL} we  conclude  that the cycles   $\balph_I$ form an integral  basis of $H_\bullet(Un),\bZ)$.  This basis coincides with the basis described in \cite[IV \S 3]{Ep}, and by Vasiliev in  \cite{Vas}.

The cycle $\balph_I$ has codimension $\codim\balph_I=\sum_{i\in I}(2i-1)$. We denote by $\balph_I^\dag\in H^\bullet(U(n),\bZ)$ its  Poincar\'e dual. When $I$ is a singleton, $I=\{i\}$, we use the simpler notation $\balph_i$ and $\balph_i^\dag$ instead of $\balph_{\{i\}}$ and respectively $\balph_{\{i\}}^\dag$.    We call the cycles $\balph_i$ the \emph{basic Arnold-Schubert cycles}.

  It is well known that  the cohomology of $U(n)$ is related via transgression  to the cohomology of its classifying space $BU(n)$. We prove that the basic class  $\balph_i$ is  obtainable by transgression  from the  Chern class $c_i$.

 More precisely, denote by $E$ the rank $n$ complex vector bundle over $S^1\times U(n)$ obtained from the trivial vector bundle
 \[
 \bC^n\times \,\Bigl(\, [0,1]\times U(n)\ra [0,1]\times U(n)\,\Bigr)
 \]
 by identifying the point $\vec{z}\in \bC^n$ in the fiber over $(1,g)\in [0,1]\times U(n)$ with the point $g\vec{z}$ in the fiber over $(0,g)\in [0,1]\times U(n)$. We denote by $p:S^1\times U(n)\ra U(n)$ the natural projection, and by $p_!: H^\bullet(S^1\times U(n),\bZ)\ra H^{\bullet-1}(U(n),\bZ)$ the induced Gysin map.

 The  first  main result of this paper  is  a transgression formula (Theorem \ref{th: TP}) asserting that
 \begin{equation*}
 \balph_i=p_!\, \bigl(\, c_i(E)\,\bigr).
 \tag{$\dag$}
 \label{tag: dag}
 \end{equation*}
In particular, we deduce that  the integral   cohomology ring  is an exterior  algebra with generators $\balph_i^\dag$, $i=1,\dotsc, n$, so that an integral basis of $H^\bullet(U(n),\bZ)$ is given by the exterior monomials
\[
\balph_{i_1}^\dag\cup\cdots \cup\balph_{i_k}^\dag,\;\; 1\leq i_1<\cdots <i_k\leq n,\;\; 0\leq k\leq n.
\]
The second  main result of this paper, Theorem \ref{th: schubert} gives a description of the Poincar\'e  dual of  $\balph_{i_1}^\dag\cup\cdots \cup\balph_{i_k}^\dag$ as a   degeneracy cycle.  More precisely, if $I=\{i_1<\cdots<i_k\}$, then
\begin{equation*}
\balph_I^\dag =\balph_{i_1}^\dag\cup\cdots \cup\balph_{i_k}^\dag.
\tag{$\ddag$}
\label{tag: ast}
\end{equation*}
The last  equality completely characterizes the  intersection ring of $\Lag(n)$ in terms of the integral basis $\balph_I^\dag$.

The sought for localization formul{\ae} are built in  our Morse theoretic approach. More precisely, if $\Phi_t$ denotes the (downward) gradient flow of  the Morse function $f$, then the results of \cite{HL} imply that the forms $\Phi_t^*\balph_I^\dag$ converge as currents when $t\ra -\infty$ to the currents  $\balph_I$.

We want to comment a bit about the flavor of the  proofs.        The   lagrangian grassmannian $\Lag(n)$ has  double incarnation: as the unitary group, and as a   collection of vector subspaces, and each of these points of view has  its uses. The unitary group interpretation is   very well suited for  global    problems, while the grassmannian incarnation is  ideal for local computations.

In this paper we solve  by local means a global    problem, the  computation of  the intersection  of two cycles,   and    not surprisingly, both incarnations of $\Lag(n)$ will play a role in the final solution. Switching between the two points of view requires some lengthy but elementary computations.

The intersection theory investigated in this paper is closely related to the traditional Schubert calculus on complex grassmannians, but uses surprisingly little of  the  traditional technology.       The intersection theory on $\Lag(n)$ has one added layer of difficulty because the  cycles involved could be odd dimensional, and when computing   intersection numbers   one  has to  count a \emph{signed} number of  points, not just a number of geometric points.  Not surprisingly, the computations of these signs  turned out to be a rather tedious job.  Moreover,  given that the cycles involved are represented by    singular real semi-algebraic objects,   the  general position  arguments are a bit more delicate.

Finally,  a few words about the organization of the paper.  The first two sections survey    known material. In Section 1 we describe carefully Arnold's isomorphism $U(n)\ra \Lag(n)$,  while in Section 2 we   describe  the most salient facts concerning the  Morse function $f_A$.

In Section 3 we give an explicit description  of the unstable manifold $W_I^-$ using   Arnold's graph coordinates. In these coordinates,  the unstable manifolds    become identified with certain  vector subspaces of the  vector space of $n\times n$ hermitian matrices. This allows us to identify the unstable manifolds with orbits of a (real) Borel group,  and we use this fact to  conclude that  the  gradient flow satisfies the Morse-Smale condition.

In Section 4 we  investigate the tunnellings of the Morse flow, i.e.,      gradient flow lines connecting two critical points $S_I, S_J$. We introduce a binary relation ``$\prec$'' on the set of critical points by declaring  $S_J\prec S_I$ if and only if  there exists a gradient flow line tunnelling from $S_I$ to $S_J$. In  Proposition \ref{prop: order} we give purely combinatorial description of this  relation  which implies  that ``$\prec$'' is in fact a \emph{partial order}. (The transitivity of ``$\prec$''  is  a reflection of the  Morse-Smale condition.) In Proposition \ref{prop: strat}  we  give a more geometric explanation of this transitivity by  showing that $S_J\prec S_I$ if and only if $W_J^-\subset \cl(W_I^-)$.   This shows that $\prec$  is   very similar to the  Bruhat order  in the classical case of grassmannians,  and  leads to a  natural stratification of the closure of an unstable manifold, where each stratum is itself an unstable  manifold.

In Section \ref{s: 4}  we  introduce the cycles $\balph_I$, and  we show that they form an integral basis of the homology.  This section is rather long since we  had to carefully  describe our orientation conventions. In Section \ref{s: 5a}  we prove the  transgression equality (\ref{tag: dag}) using the    Thom-Porteous formula  describing the duals of  Chern classes as certain degeneracy loci.    Section \ref{s: 5}  contains the proof  of the main intersection   result, the identity  (\ref{tag: ast}).

For the reader's convenience we included  two technical appendices. The first one is about  tame (or $o$-minimal)  geometry and  tame flows.   We need these facts to prove that the   gradient  flow of $f$ is Morse-Stokes in the sense of \cite{HL}.  The criteria for recognizing Morse-Stokes   flows    proved by Harvey-Lawson in \cite{HL} do not apply to our     gradient flow,   and since    our flow  can be described quite explicitly, the  $o$-minimal technology is ideally suited for this task.

The second appendix   is a fast paced overview of R. Hardt's intersection theory of subanalytic currents.    We   describe  a weaker form of transversality, and  explain in Proposition \ref{prop: inter}  how  to compute  intersections  under this milder  conditions

\bigskip

\noindent {\bf Acknowledgements.} I want to thank my colleagues Sam Evens,  Sergei Starchenko and Bruce Williams for useful conversations on   topics related to this paper.

\section*{Notations and conventions}
\addtocontents{section}{Notations and conventions}

\begin{itemize}

\item  For any finite set $I$, we denote by $\# I$ or $|I|$ its cardinality.

\item $\ii:=\sqrt{-1}$.

\item $\bI_n:=\bigl\{\, -n,\dotsc, -1, 1,\dotsc, n\}$, $\bI_n^+=\{1,\dotsc, n\}$.

\item   For an oriented manifold $M$ with boundary  $\pa M$, the induced orientation on the boundary is obtained using the outer-normal first convention.

\item For any subset $S$ of a topological space $X$ we denote by $\cl(S)$ its closure in $X$.

\item For any complex hermitian  vector space we denote by $\End^+(E)$ the space of hermitian linear operators $E\ra E$.

\item  For every complex vector space $E$ and every  nonnegative integer $m\leq \dim_\bC E$ we  denote by $\Gr_m(E)$ (respectively $\Gr^m(E)$) the Grassmannian of  complex subspaces of $E$ of dimension $m$ (respectively codimension $m$).

\item  Suppose $E$ is  a complex  Euclidean vector space  of dimension $n$ and
\[
\Fl:= \bigl\{ \bsF_0\subset \bsF_1\subset \cdots \subset \bsF_n\,\bigr\}
\]
is a complete flag of subspaces of $E$, i.e., $\dim \bsF_i=0$, $\forall i=0,\dotsc, n$.

For  every integer $0\leq m \leq n$, and every partition $\mu =\mu_1\geq \mu_2\cdots$ such that $\mu_1\leq m$ and $\mu_i=0$, for all $i>n-m$, we define the Schubert cell $\Sigma_\mu(\Fl)$ to be the subset of $\Gr^m(E)$ consisting of  subspaces $V$ satisfying the incidence relations
\[
\dim(V\cap \bsF_j)=i,
\]
$\forall i=1,\dotsc, m$, $\forall j,\;\;m+i-\mu_i\leq j\leq m+i-\mu_{i+1}$.

\end{itemize}

\section{Hermitian lagrangians}
\label{s: 1}
\setcounter{equation}{0}
We  would like to collect  in this section a few basic facts
concerning hermitian lagrangian spaces which we will  need  in our
study.  All of the results   are due to V.I. Arnold, \cite{Ar1}. In
this section all vector spaces will be assumed finite dimensional.

\begin{definition} A  \emph{hermitian symplectic  space}  is a pair $(\widehat{E}, J)$, where $\widehat{E}$ is a  complex  hermitian space, and    $J:\widehat{E}\ra\widehat{E}$  is a unitary operator such that
\[
J^2=-\one_{\whE},\;\;\dim_\bC(\ker(J-\ii)=\dim_\bC\ker(J+\ii).
\]
An isomorphism of hermitian symplectic spaces  $(\widehat{E}_i, J_i)$, $i=0,1$,
is  a  map unitary  map $T:\widehat{E}_0\ra \widehat{E}_1$ such that  $TJ_0=J_1T$.\qed
\end{definition}

If $(\widehat{E}, R, J)$ is a hermitian symplectic space, and
$h(\bullet,\bullet)$ is the hermitian metric on $\whE$, then the
\emph{symplectic hermitian form} associated to   this space is  the
form
\[
\omega: \whE\times \whE \ra \bC,\;\;\omega(\bu,\bv)=h(J\bu,\bv).
\]
Observe that $\omega$ is linear in  the first variable and conjugate
linear in the second variable. Moreover
\[
\omega(\bu,\bv)=-\overline{\omega(\bv,\bu)},\;\;\forall\bu,\bv\in
\whE.
\]
The $\bR$-bilinear  map
\[
q(\bu,\bv): =\re h(\ii J \bu,\bv)
\]
is symmetric, nondegenerate and  has signature  $0$. We denote  by
$\Symp(\whE, J)$ the subgroup of $\GL_\bC(\whE)$ consisting   of
complex linear      automorphisms of $\whE$ which preserve $\omega$,
i.e.,
\[
\omega(T\bu,\bv)=\omega(\bu,\bv),\;\;\forall\bu,\bv\in \whE.
\]
Equivalently,
\[
\Symp(\whE, J)=\bigl\{ T \in \GL_\bC(\whE);\;\; T^* J T =J\,\bigr\}.
\]
Observe that $\Symp(\whE, J)$ is isomorphic to the noncompact Lie group $U(n,n)$, $n=\frac{1}{2}\dim_\bC\whE$.  We denote by $\usymp(\whE,J)$ its Lie algebra.

We set $F^\pm :=\ker(\pm \ii -J)$.   We fix an \emph{isometry} $T: F^+ \ra F^-$ and we set
\[
\whE^+= \Bigl\{ \,  \frac{1}{2}(f+Tf);\;\;f\in F^+\,\Bigr\},\;\;  \whE^-= \Bigl\{ \,  \frac{1}{2\ii}(f-Tf);\;\;f\in F^+\,\Bigr\}
\]
Observe that $\whE^-$ is the orthogonal complement of $\whE^+$, and  the operator $J$  induces a  unitary
isomorphism $\whE^+\ra \whE^-$. Thus,  we can think of $\whE^\pm$ as two different copies of the same hermitian space $E$.

 Conversely, given a hermitian space $E$, we can form $\whE=E\oplus E$,  and define  $\bsJ:\whE\ra \whE$  by
with reflection
\[
\bsJ= \left[
\begin{array}{cc}
0 &-\one_{E}\\
\one_{E} & 0
\end{array}
\right].
\]
Note that
\[
F^\pm = \bigl\{  (\bx,\mp \bx)\in  E\oplus E;\;\;\bx \in E\,\bigr\},
\]
and we have a canonical  isometry
\[
F^+\ni (\bx,-\ii \bx)\stackrel{T}{\longmapsto} (\bx, \ii\bx)\in F^-.
\]
For this reason, in the sequel we will assume that our hermitian
symplectic spaces have the  standard form
\[
\widehat{E}= E\oplus E,\;\; \bsJ=\left[
\begin{array}{cc}
0 &-\one_{E}\\
\one_{E} & 0
\end{array}
\right],
\]
We set $\whE^+=E\oplus 0$,  $\whE^-=0\oplus E$. We say that $\whE^+$ (respectively $\whE^-$) is the \emph{horizontal} (respectively \emph{vertical})  component of $\whE$.

\begin{definition} Suppose $(\widehat{E}, \bsJ)$ is a hermitian symplectic space.  A \emph{hermitian lagrangian subspace} of $\widehat{E}$ is a \emph{complex} subspace  $L\subset \widehat{E}$ such that $L^\perp=\bsJ L$. We will denote  by $\Lag(\widehat{E})$ the set of hermitian lagrangian subspaces of  $\widehat{E}$.\qed
\end{definition}

\begin{remark}  If $\omega$ is the symplectic form associated to $(\widehat{E}, \bsJ)$ then a subspace  $L$ is  hermitian lagrangian if and only if
\[
L=\bigl\{ \bu\in \whE;\;\;\omega(\bu,\bx)=0,\;\;\forall \bx\in
L\,\bigr\}.
\]
This shows that the group $\Symp(\whE, J)$ acts  on $\Lag(\whE)$,
and it is not hard to prove that the action is transitive.\qed
\end{remark}

Observe that if $L\in \Lag(\widehat{E})$ then  we have a natural
isomorphism $L\oplus JL\ra \widehat{E}$. It follows that $\dim_\bC
L=\frac{1}{2}\dim_\bC \widehat{E}$, and if we set $2n:=\dim_\bC
\widehat{E}$ we deduce that $\Lag(\widehat{E})$ is  a subset of the
Grassmannian  $\Gr_n(\widehat{E})$ of complex  $n$-dimensional
subspaces of $\widehat{E}$. As such, it  is  equipped with an
induced topology.

\begin{ex}  Suppose $E$ is a complex hermitian  space. To any  linear operator $A:E\ra E$ we associate its graph
\[
\Gamma_A= \bigl\{ (x, Ax)\in E\oplus E;\;\; x\in E\,\bigr\}.
\]
Then $\Gamma_A$ is  a hermitian lagrangian subspace of $E\oplus E$ if
and only if the operator $A$ is self-adjoint.

More generally, if $L$ is a lagrangian subspace in a hermitian
symplectic vector space $\whE$, and $A:L\ra L$ is a linear operator,
then the  graph of $JA: L\ra JL$ viewed as a subspace in $L\oplus
JL=\whE$ is a lagrangian subspace if and only if $A$ is a Hermitian
operator.\qed \label{ex: sym}
\end{ex}

\begin{lemma} Suppose $E$ is a complex hermitian  space, and  $S: E\ra E$ is a linear  operator.    Define
\begin{equation}
\eL_S: =\left\{ \,\bigl(\,(\one+S)x,
 -\ii(\one-S)x\, \bigr);\;\;x\in E\,\right\}\subset E\oplus E.
\label{eq: LS}
\end{equation}
Then $\eL_S\in \Lag(E\oplus E)$ if and only if $S$ is a unitary
operator. \label{lemma: unit-lag}
\end{lemma}

\proof  Observe that $\eL_S$ is the image of the linear map
\[
E\mapsto  E\oplus E,\;\;x\mapsto \eI_S(x) =\left(\, \bigl(\one+S\bigr)x, -\ii\bigl(\one-S\bigr)x\,\right).
\]
This map is   injective because
\[
(\one+S)x=(\one-S)x=0\Longrightarrow  x=Sx=-Sx\Longrightarrow x=0.
\]
Hence $\dim_\bC \eL_S=\dim_\bC E=\frac{1}{2}\dim_\bC E\oplus
E=\dim_\bC\eL_S^\perp$, and we deduce  that $\eL_S$ is a lagrangian
subspace if and only  if $J\eL_S\subset \eL_S^\perp$.

We denote by $(\bullet, \bullet)$ the Hermitian scalar product in
$E$ and by $\lan\bullet,\bullet\ran$ the Hermitian inner product in
$E\oplus  E$. If
\[
u=\eI_S(x)= \left(\, \bigl(\one+S\bigr)x, -\ii\bigl(\one-S\bigr)x\,\right)
\Longrightarrow \bsJ u= \left(\, \ii \bigl(\one-S\bigr)x, \bigl(\one+S\bigr)x\,\right)
\]
For every $v=(\one+S)y\oplus  -\ii(\one-S)y=\eI_S(y)\in \eL_S$ we
have
\[
\lan \bsJ u, v\ran =\ii (x-Sx, y+Sy) +\ii (x+Sx, y-Sy)
\]
Now observe that
\[
(x-Sx, y+Sy)= (x,y)+ (x,Sy)-(Sx,y)-(Sx,Sy)=  (x,Sy)-(Sx,y),
\]
and
\[
(x+Sx, y-Sy)= (x,y) -(x, Sy)+ (Sx,y)-(Sx,y) = -(x, Sy)+ (Sx,y)
\]
so that
\[
\lan J\eI_S(x), \eI_S(y)\ran=2(x,y)-2(Sx,Sy),\;\;\forall x,y\in E.
\]
We deduce that $J\eL_S\subset \eL_S^\perp$ if and only if
\[
(x,y)=(Sx,Sy),\;\;\forall x,y\in E \Longleftrightarrow \mbox{the
operator $S$ is unitary} .\proofend
\]
\begin{lemma}  If $L\in \Lag(\widehat{E})$ then $L\cap F^\pm=\{0\}$.
\label{lemma: trans}
\end{lemma}

\proof  Suppose  $f\in F^\pm\cap L$. Then   $Jf\in L^\perp$ so that
$\lan Jf, f\ran =0$. On the other hand, $Jf=\pm \ii f$ so that
\[
0 =\lan Jf,f\ran =\pm \ii|f|^2\Longrightarrow f=0.\proofend.
\]
Using  the isomorphism $\whE=F^+\oplus F^-$ we deduce from the above
lemma  that  $L$ can be represented as the graph of a linear
isomorphism $T=T_L: F^+\ra F^-$, i.e.,
\[
L=  \bigl\{ f\oplus Tf;\;\;f\in F^+\,\bigr\}.
\]
Denote by $\eI_\pm: E\ra F^\pm$ the unitary map
\[
E\ni x\mapsto \frac{1}{\sqrt{2}}( x,\mp \ii x )\in F^\pm.
\]
We denote by $\eS_L: E\ra E$ the  linear map  given by the
commutative diagram
\[
\begin{diagram}
\node{E} \arrow{s,l}{\eI_+} \arrow{e,t}{\eS_L}\node{E}\arrow{s,r}{\eI_-}\\
\node{F_+}\arrow{e,b}{T}\node{F^-}
\end{diagram},
\]
i.e., $\eS_L= \eI_-^{-1}T\eI_+$.  A simple computation  shows that $L=\eL_{\eS_L}$. From Lemma \ref{lemma: unit-lag} we deduce that the operator $\eS_L$
is unitary, and   that the map
\[
\Lag(\whE)\ni L\mapsto \eS_L\in U(E)
\]
is the inverse of the map $S\mapsto \eL_S$.    This proves the following result.

\begin{proposition}[Arnold] Suppose $E$ is a complex hermitian  space, and denote by $U(E)$ the group of unitary operators $S:E\ra E$. Then the map
\[
U(E)\ni S\mapsto \eL_S\in \Lag(E\oplus E)
\]
is a homeomorphism. In particular, we deduce that $\Lag(\whE)$ is a
smooth,  compact, connected, orientable \emph{real} manifold of
dimension
\[
\dim_\bR \Lag(\whE)=\bigl(\dim_\bC E\,\bigr)^2. \proofend
\]
\label{prop: arnold}
\end{proposition}

Suppose $A:E\ra E$ is a selfadjoint operator. Then, as we know its
graph $\Gamma_A$ is a lagrangian in $\whE=E\oplus E$, and thus there exists a unitary  operator $S\in U(E)$ such that
\[
\Gamma_A=\eL_S=\left\{ \,\bigr(\, (\one+S) x, -\ii(\one-S)x\,\bigr);\;\;x\in E\,\right\}.
\]
Note that the graph $\Gamma_A$ intersects the  ``vertical axis''  $\whE^-=0\oplus E$ only at the origin so that the operator $\one +S$ must be invertible.

Next observe that for every $x\in E$ we have $-\ii(\one-S)x= A(\one+ S)x$ so that
\begin{equation}
A= -\ii (\one-S)(\one+S)^{-1}= 2\ii (\one+S)^{-1}-\ii.
\label{eq: sym-unit0}
\end{equation}
We conclude
\begin{equation}
S=\eS_{\Gamma_A}= \eC(\ii A):=(\one-\ii A)(\one+\ii A)^{-1}=2(\one+\ii A)^{-1}-\one.
\label{eq: sym-unit}
\end{equation}
The   expression $\eC(\ii A)$ is the so called  \emph{Cayley transform} of
$\ii A$.

Observe that we have a left  action ``$\ast$'' of $U(E)\times U(E)$
on $U(E)$ given by
\[
(T_+, T_-)\ast S= T_-ST_+^*,\;\;\forall T_+,T_-, S\in U(E).
\]

To obtain  a  lagrangian  description of this action   we need to
consider  the \emph{symplectic unitary group}
\[
U(\whE,\bsJ):=U(\whE)\cap \Symp(\whE,\bsJ)=\bigl\{ \, T\in
U(\whE);\;\;T\bsJ=\bsJ T\,\bigr\}.
\]
The subspaces $F^\pm$ are invariant subspaces of   any   operator
$T\in U(\whE,\bsJ)$   so that we have an isomorphism
$U(\whE,\bsJ)\cong U(F^+)\times U(F^-)$.   Now identify $F^\pm$ with
$E$ using the isometries
\[
\frac{1}{\sqrt{2}}\eI_\pm :E\ra F^\pm,\;\;\eI_J: E\oplus E\ra
F^+\oplus F^-.
\]
We obtain an isomorphism
\[
U(\whE,\bsJ)\ni T\mapsto (T_+, T_-) \in U(E)\times U(E)..
\]
Moreover, for any lagrangian $L\in \Lag(\whE)$,  and $S\in U(E)$,
and any $T\in U(\whE,\bsJ)$ we have
\begin{equation}
\eS_{T L}= (T_+, T_-)\ast \eS_L ,\;\;\eL_{(T_+,T_-)\ast S}= T\eL_S.
\label{eq: equivar}
\end{equation}

\section{Morse flows on the  Grassmannian of hermitian lagrangians}
\label{s: 2}
\setcounter{equation}{0}
In this  section we would like to describe a few  properties of some
nice Morse functions on the Grassmannian of complex lagrangian
subspaces.    The main source for  all these facts is the very nice
paper by  I.A. Dynnikov and A.P. Vesselov, \cite{DV}.

Suppose $E$ is  complex Hermitian space of complex dimension $n$. We
equip the space $\whE=E\oplus E$ with the canonical complex
symplectic structure. Recal that
\[
\whE^+:= E\oplus 0,\;\; \whE^-:=0\oplus E.
\]
For every symmetric operator $A: \whE^+\ra \whE^+ $ we denote by $\widehat{A}:\whE\ra \whE$ the symmetric operator
\[
\widehat{A} =\left[
\begin{array}{cc}
A  &0\\
0 &-A
\end{array}
\right]:\whE\ra \whE.
\]
Let us point out that $\whA\in \usymp(\whE,\bsJ)$. Define
\[
f_A:U(\whE^+)\ra \bR,\;\;f_A(S):=\re \tr(AS),
\]
and
\[
\vfi_A: \Lag(\whE)\ra  \bR,\;\; \vfi_A(L)= \re\tr(\widehat{A}P_L),
\]
where $P_L$ denotes the orthogonal projection onto $L$. An elementary computation shows that
\begin{equation}
P_{\eL_S}= \frac{1}{2} \left[
\begin{array}{cc}
\one +\frac{1}{2}(S+S^*) & \frac{\ii}{2}(S-S^*)\\
\frac{\ii}{2}(S-S^*) & \one-\frac{1}{2} (S+S^*)
\end{array}
\right],
\label{eq: unit-proj}
\end{equation}
 and we deduce
\[
\vfi_A(\eL_S)=f_A(S),\;\;\forall S\in U(\whE^+).
\]
The following result is classical, and it goes back to Pontryagin \cite{Pont}. For a proof we refer to \cite{DV}.

\begin{proposition} If $\ker A=\{0\}$  then a unitary operator  $S\in U(\whE^+)$ is a critical point  of $f_A$ if and only if   there exists a unitary basis  $e_1,\dotsc, e_n$ of $E$ consisting of eigenvectors of $A$  such that
\[
Se_k=\pm e_k,\;\;\forall k=1,\dotsc, n. \proofend
\]
\end{proposition}

We can reformulate the  above result by saying that when $\ker A\neq
0$, then  a unitary operator  $S$ is a critical point of $f_A$ if
and only if  $S$ is an involution  and both $\ker(\one -S)$ and
$\ker(\one+S)$ are   invariant subspaces of $A$.   Equivalently this means
\[
S=S^*,\;\;S^2=\one_E,\;\;SA=AS.
\]
To obtain more  detailed  results, we  fix an orthonormal basis $e_1,\dotsc, e_n$ of $E$. For any    $\vec{\alpha}\in \bR^n$ such that
\begin{equation}
0<\alpha_1<\cdots <\alpha_n,
\label{eq: ord}
\end{equation}
we denote by $A=A_{\vec{\alpha}}$ the symmetric operator $E\ra E$ defined by $Ae_k= \alpha_k e_k$, $\forall k$, and we set $f_{\vec{\alpha}}:=f_{A_{\vec{\alpha}}}$, and by $\Cr_{\vec{\alpha}}\subset U(\whE^+)$, the set of critical points of $f_{\vec{\alpha}}$.

For every $\vec{\eps}\in \{\pm 1\}^n$ we define $S_{\vec{\eps}}\in U(\whE^+)$ by
\[
S_{\vec{\eps}}e_k=\eps_k e_k,\;\; k=1,2,\dotsc, n.
\]
Then
\[
\Cr_{\vec{\alpha}}=\bigl\{\, S_{\vec{\eps}};\;\vec{\eps}\in\{\pm 1\}^n\,\bigr\}.
\]
Note that  this critical set is independent of the vector $\vec{\alpha}$ satisfying (\ref{eq: ord}). For this reason we will   use the simpler notation $\Cr_n$ when referring to this critical set.

 To compute the index of $f_{\vec{\alpha}}$ at the critical point $S_{\vec{\eps}}$ we need to compute the Hessian
 \[
 \eQ_{\vec{\eps}}(H):=\frac{d^2}{dt^2}|_{t=0}  \re\tr (\, A S_{\vec{\eps}} e^{tH}\,),\;\;H\in \uu(E) \mbox{= the Lie algebra of $U(E)$}.
 \]
We have
\[
\eQ_{\vec{\eps}}(H) = \re\tr A_{\vec{\alpha}}S_{\vec{\eps}} H^2)= -\re\tr A_{\vec{\alpha}}S_{\vec{\eps}} H H^*
\]
Using the  basis $(e_i)$ we can represent $H\in\uu(E)$  as $H= \ii Z$, where $Z$  is  a hermitian matrix $\bigl(\, z_{jk}\,\bigr)_{1\leq ij,\leq n}$, $z_{jk}=\bar{z}_{kj}$. Note that $z_{jj}$ is a \emph{real} number, while $z_{ij}$  can be any complex number if $i\neq j$. Then a simple computation shows
\begin{equation}
\eQ_{\vec{\eps}}(\ii Z)=-\sum_{i,j} (\eps_i\alpha_i+ \eps_j\alpha_j)|z_{ij}|^2=-\sum_i \eps_i \alpha_i |z_{ii}|^2-2\sum_{i<j} (\eps_i\alpha_i+ \eps_j\alpha_j)|z_{ij}|^2.
\label{eq:  hess}
\end{equation}
Hence, the index of $f_{\vec{\alpha}}$ at $S_{\vec{\eps}}$ is
\[
\mu_{\vec{\alpha}}(\vec{\eps}):= \# \{ i;\;\;\eps_i=1\} + 2 \#\bigl\{ (i,j);\;\;i<j,\;\; \eps_i\alpha_i +\eps_j\alpha_j >0\,\bigr\}
\]
Observe that if $i<j$ then   $\eps_i\alpha_i +\eps_j\alpha_j >0$ if and only if $ \eps_j =1$. Hence
\[
\mu_{\vec{\alpha}}(\vec{\eps})=\sum_{\eps_j=1} (2j-1).
\]
In particular, we see that the index is independent of the vector $\vec{\alpha}$ satisfying the conditions (\ref{eq: ord}).

It is convenient to introduce another parametrization of   the critical set. Recall that
\[
\bI_n^+:=\bigl\{\,1,\dotsc, n\,\bigr\}.
\]
For every subset $I\subset \bI_n^+$ we  denote by  $S_I\in U(\whE^+)$ the unitary operator defined   by
\[
S_Ie_j=\begin{cases} e_j & j\in I\\
-e_j &j\not\in I.
\end{cases}
\]
Then $\Cr_n:=\bigl\{ S_I;\;\;I\subset \bI_n^+\,\bigr\}$, and the index of $S_I$ is
\begin{equation}
{\rm ind}\,(S_I)= \sum_{i\in I} (2i-1).
\label{eq: index}
\end{equation}
The co-index is
\begin{equation}
{\rm coind}\,(S_I)={\rm ind}\,(S_{I^c})=n^2-\mu_I,
\label{eq: coind}
\end{equation}
where $I^c$ denotes the complement of $I$, $I^c:=\bI_n^+\setminus I$.

\begin{definition}  We define the \emph{weight} of a finite subset $I\subset \bZ_{>0}$ to be the integer
\[
\bw(I)=\begin{cases}
0 & I=\emptyset\\
\sum_{i\in I}(2i-1) & I\neq \emptyset.
\end{cases} \proofend
\]
\end{definition}

Hence ${\rm ind}\,(S_I)= \bw(I)$. Let us observe a remarkable  fact.

\begin{proposition} Let
\[
\vec{\xi}=\Bigl(\,\frac{1}{2},\frac{3}{2},\cdots, \frac{2n-1}{2}\,\Bigr)\in \bQ^n,
\]
and set
\[
f_0:=f_{\vec{\xi}},\;\;\vfi_0:=\vfi_{\vec{\xi}}.
\]
Then for every  $I\subset \bI^+_n$ we have
\[
\bw(I)=  f_0(S_I)+\frac{n^2}{2}  =\vfi_0(\Lambda_I) +\frac{n^2}{2} .
\]
In other words  the gradient flow of $f_{\xi}$ is selfindexing, i.e.,
\[
f_0(S_I)-f_0(S_J)=\bw(J)-\bw(I).
\]
\label{prop: self-ind}
\end{proposition}

\proof We have
\[
f_0(S_I)= \frac{1}{2}\bigl(\,\bw(I)-\bw(I^c)\,\bigr).
\]
On the other hand, we have
\[
\frac{1}{2}\bigl(\, \bw(I)+\bw(I^c)\,\bigr)=\frac{1}{2}\bw(\bI^+_n)=\frac{n^2}{2}.
\]
Adding up the above equalities  we obtain the desired conclusion.\qed

The (positive) gradient flow of the function  $f_A$ has an  explicit description. More precisely, we have the following result, \cite[Proposition 2.1]{DV}.

\begin{proposition} Suppose $A=A_{\vec{\alpha}}$ where $\vec{\alpha}\in \bR^n$ satisfies (\ref{eq: ord}).  We equip $U(\whE^+)$ with the left invariant metric induced  from the inclusion in the  Euclidean space $\End_{\bC}(E)$ equipped with the inner product
\[
\lan X,Y\ran =\re\tr(XY^*).
\]
We denote by $\nabla f_A$ the gradient of $f_A:U(\whE^+)\ra \bR$  with respect to this metric, and  we denote by
\[
\Phi_A:\bR\times U(\whE^+)\ra U(\whE^+),\;\;S\mapsto \Phi_A^t(S)
\]
the flow defined by $\nabla f_A$, i.e., the flow associated to the o.d.e. $\dot{S}=\nabla f_A(S)$.Then
\begin{equation}
\Phi_A^t(S)=\bigl(\, \sinh(tA) +\cosh(tA) S\,\bigr)(\,\cosh(tA)+\sinh(tA)S\,\bigr)^{-1},\;\;\forall S\in U(\whE^+),\;\; t\in \bR.
\label{eq: flow}
\end{equation}\qed
\end{proposition}

It is convenient to have a lagrangian description of the above results via the diffeomorphism $\eL: U(\whE^+)\ra \Lag(\whE)$.  First, we use this isomorphism to   transport isometrically the metric on $U(\whE^+)$. Next,  for every $I\subset \bI_n^+$ we set $\Lambda_I:=\eL_{S_I}$.     For every $i\in\bI_n^+$  we define
\[
\be_i:= e_i\oplus 0\in E\oplus E,\;\;\bsf_i=0\oplus e_i\in E\oplus E.
\]
Then
\[
\Lambda_I=\ker(\one-S_I)\oplus \ker(\one+S_I)= {\rm span}\, \bigl\{ \be_i;\;\;i\in I\,\bigr\}+ {\rm span}\, \bigl\{ \bsf_j;\;\;j\in I^c\,\bigr\}.
\]
The lagrangians $\Lambda_I$ are the critical points  of the function $\vfi_A:\Lag (\whE)\ra \bR$.

Using (\ref{eq: LS}) and (\ref{eq: flow}) we deduce that for every $S\in U(\whE^+)$ we have
\begin{equation}
\eL_{\Phi^t(S)}= e^{t\widehat{A}} \eL_S.
\label{eq: flow1}
\end{equation}
The above equality  describes the  (positive) gradient flow of $\vfi_A$. We denote this flow  by $\Psi^t_A$.

We can use the lagrangians   $\Lambda_I$ to produce the  \emph{Arnold atlas}, \cite{Ar0}. Define
\[
\Lag(\whE)_I:=\bigl\{ L\in \Lag(\whE);\;\;L\cap \Lambda_I^\perp=0,\bigr\}.
\]
Then $\Lag(\whE)_I$ is an open subset of $\Lag(\whE)$ and
\[
\Lag(\whE)=\bigcup_I \Lag(\whE)_I.
\]
Denote by $\End^+_\bC(\Lambda_I)$ the space of self-adjoint endomorphisms of $\Lambda_I$. We have a diffeomorphism
\[
\End^+_\bC(\Lambda_I)\ra \Lag(\whE)_I,
\]
which associates to each symmetric  operator $T: \Lambda_I\ra \Lambda_I$, the graph $\Gamma_{{\bsJ}T}$  of the operator
\[
{\bsJ}T:\Lambda_I\ra \Lambda_I^\perp
\]
 regarded as a subspace in $\Lambda_I\oplus \Lambda_I^\perp\cong \whE$.  More precisely, if  the operator $T$ is   described in the orthonormal basis $\{ \be_i,\;\bsf_j;\;\;i\in I,\;\;j\in I^c\}$  by the  Hermitian matrix $(t_{ij})_{1\leq i,j\leq n}$, then  the graph of ${\bsJ}T$ is spanned by the vectors
\begin{subequations}
\begin{equation}
\be_i(T):=\be_i +\sum_{i'\in I} t_{i'i}\bsf_{i'} -\sum_{j\in I^c} t_{ji} \be_j,\;\;i\in I,
\label{eq: basee}
\end{equation}
\begin{equation}
\bsf_j(T):=\bsf_j+\sum_{i\in I} t_{ij}\bsf_i -\sum_{j'\in I^c} t_{ j'j} \be_{j'|} ,\;\;j\in I^c.
\label{eq: basef}
\end{equation}
\end{subequations}
The inverse map
\[
\eA_I:\Lag(\whE)_I\ra \End^+_\bC(\Lambda_I)
\]
 is known as the \emph{Arnold coordinates} on $\Lag(\whE)_I$.

Let $I\subset \bI_n^+$. If $L\in \Lag(\whE)_I$ has Arnold coordinates $\eA_I(L)=T$, i.e., $T$ is a symmetric operator
\[
T:\Lambda_I\ra \Lambda_I,
\]
and $L=\Gamma_{{\bsJ}T}$, then $ \Phi_A^{t}L=e^{t\whA}\Gamma_{{\bsJ}T}$ is spanned by the vectors
\[
e^{t\alpha_i}\be_i +\sum_{i'\in I} t_{i'i}e^{-t\alpha_{i'}}\bsf_{i'} -\sum_{j\in I^c} t_{ji} e^{t\alpha_j} \be_j,\;\;i\in I,
\]
\[
e^{-t\alpha_j}\bsf_j+\sum_{i\in I} t_{ij}e^{-t\alpha_i}\bsf_i -\sum_{j'\in I^c} t_{j' j} e^{t\alpha_{j'}}\be_{j'} ,\;\;j\in I^c
\]
or, equivalently, by the vectors
\begin{subequations}
\begin{equation}
\be_i +\sum_{i'\in I} t_{i'i}e^{-t(\alpha_{i'}+\alpha_i)}\bsf_{i'} -\sum_{j\in I^c} t_{ji} e^{t(\alpha_j-\alpha_i)} \be_j,\;\;i\in I,
\label{eq: span1}
\end{equation}
\begin{equation}
\bsf_j+\sum_{i\in I} t_{ij}e^{t(\alpha_i-\alpha_j)}\bsf_i -\sum_{j'\in I^c} t_{j'j} e^{t(\alpha_{j'}+\alpha_j)}\be_{j'} ,\;\;j\in I^c.
\label{eq: span2}
\end{equation}
\end{subequations}
This shows that $e^{t\whA}\Gamma_{{\bsJ}T}\in \Lag(\whE)_I$, so that $\Lag(\whE)_I$ is invariant under the flow $\Psi_A$.

If we denote by $A_I$ the restriction of $\whA$ to $\Lambda_I$,  and we regard $A_I$ as a symmetric operator $\Lambda_I\ra \Lambda_I$, then we deduce  from the above equalities that
\[
e^{t\whA}\Gamma_{\bsJ T} = \Gamma_{\bsJ e^{tA_I}Te^{tA_I}}.
\]
We can rewrite the above equality in terms of Arnold coordinates as
\begin{equation}
\eA_I(\Psi^t L) = e^{tA_I} \eA_I(L) e^{tA_I},\;\;\forall L\in \Lag(\whE)_I.
\label{eq: flow3}
\end{equation}

\section{Unstable manifolds}
\label{s: 3}
\setcounter{equation}{0}
The unstable manifolds  of the positive  gradient flow of $\vfi_A$
have many similarities with the Schubert cells of complex
Grassmannians,  and we would like to investigate  these similarities
in great detail.

The stable/unstable variety of $\Lambda_I$ \emph{with respect to the
positive gradient flow }$\Psi_A^{t}$ is defined by
\[
W^\pm_I:=\bigl\{ L\in \Lag(\whE);\;\;  \lim_{t\ra \infty} e^{\pm
t\widehat{A}}L =\Lambda_I\,\bigr\}\stackrel{(\ref{eq: flow3})}{=}\bigl\{ L\in \Lag(\whE)_I;\;\;\lim_{t\ra\pm \infty} e^{tA_I}\eA_I(L)e^{tA_I}=0\,\bigr\}.
\]
If $\eA_I(L)= (t_{ij})_{1\leq i,j\leq n}$ then  the equalities (\ref{eq: span1}) and (\ref{eq: span2}) imply  that
\[
\lim_{t\ra- \infty} e^{tA_I}\eA_I(L)e^{tA_I}=0 \Longleftrightarrow t_{ij}=0,\;\;\mbox{if  $i,j\in I$, or $j\in I^c$, $i\in I$ and m$j<i$}.
\]
We can rewrite the  last system of equalities  in the more compact
form
\begin{equation}
W^-_I=\bigl\{ T\in \End^+(\Lambda_I);\;\;t_{ji}=0,\;\;\forall 1\leq
j\leq i,\;\;i\in I\,\bigr\}. \label{eq: w-}
\end{equation}
This shows that $W^-_I$ has  real  codimension $\sum_{i\in
I}(2i-1)$. This agrees with our previous computation (\ref{eq:
coind}) of the index of $\Lambda_I$. Thus
\[
\codim_\bR W_I^-=\bw(I),\;\;\dim_\bR W_I^-= n^2-\bw(I)=
\frac{1}{2}\dim_\bR\Lag(\whE)-\vfi_0(\Lambda_I).
\]

For any $L\in \Lag(\whE)$ we set
\[
L^\pm:= L\cap \whE^\pm.
\]
The dimension of $L^+$ is called the \emph{depth} of  $L$, and will
be denoted by $\delta(L)$.

From the description (\ref{eq: w-}) of the unstable variety $W^-_I$,
$\# I=k$ we deduce the following result.\footnote{The
characterization in Proposition \ref{prop: ww}    depends
essentially on the fact that the eigenvalues of $A$ satisfy the
inequalities $0<\alpha_1<\cdots <\alpha_n$.}

\begin{proposition}  Let $L\in \Lag(\whE)$, $I\subset \{1,\cdots, n\}$, $k=\# I$.  We denote by $S\in U(\whE^+)$ the unitary operator corresponding to $L$. The following statements are equivalent.
\begin{enumerate}
\item[(a)] $L\in W^-_I$

\item[(b)] $L \in \Lag(\whE)_I$ and $\lim_{t\ra\infty} e^{-tA}L^+=\Lambda_I^+$.

\item[(c)] $\dim L^+=k$ and $\lim_{t\ra\infty} e^{-tA}L^+=\Lambda_I^+$.

\item[(d)] $\dim\ker(\one-S)=k$ and $\lim_{t\ra\infty} e^{-tA}\ker(\one-S)=\Lambda_I^+$.
\end{enumerate}
\label{prop: ww}
\end{proposition}

\proof  The description (\ref{eq: w-}) shows that  (a) $\Rightarrow$
(b) $\Rightarrow$ (c).    Suppose that $L$ satisfies  (c) and let
$\Lambda_J=\lim_{t\ra \infty} e^{-t\whA} L$, i.e., $L\in W_J^-$.
Then using the implication (a) $\Rightarrow$ (b) for the unstable
manifold $W_J^-$ we deduce
\[
\lim_{t\ra \infty} e^{tA}L^+=\Lambda_J^+.
\]
On the other hand, since $L$ satisfies $(c)$ we have
\[
\lim_{t\ra \infty} e^{tA}L^+=\Lambda_I^+.
\]
This implies $I=J$ which proves the implication (c) $\Rightarrow$
(a). Finally, observe that (d)  is a reformulation of (c) via the
diffeomorphism  $\eS:\Lag(\whE)\ra U(\whE^+)$. \qed

The condition  $\lim_{t\ra\infty} e^{-tA}L^+=\Lambda_I^+$  can be
rephrased  as an incidence condition. We write
\[
I=\{\nu_1 >\cdots >\nu_k\}.
\]
 We have $\lim_{t\ra \infty}e^{-tA}L^+\ra \Lambda_I^+$ if and only  if  $L^+$ is the graph  of a  linear map
\[
X: \Lambda_I^+\ra \Lambda^+_{I^c},\;\;A\be_i= \sum_{j\in I^c}
x^j_i\be_j,\;\;\forall i\in I,
\]
such that,
\[
x^j_i=0,\;\;\forall i\in I,\;\;j\in I^c,\;\; j<i.
\]
We consider the complete flag  $\Fl_\bullet=\{\,\bsF_0\subset \cdots
\subset\bsF_n\,\}$ of subspaces  of $\whE^+$,
\[
\bsF_j={\rm span}_\bC\bigl\{ \be_i;\;\;i\leq j\,\bigr\},
\]
and we form de dual  flag
$\Fl^\bullet=\{\,\bsF^0\subset\cdots\subset\bsF^n\,\}$,
\[
\bsF^j:=\bsF_{n-j}^\perp= {\rm span}_\bC\bigl\{
\be_i;\;\;i>n-j\,\bigr\}.
\]
Then $\lim_{t\ra \infty}e^{-tA}L^+\ra \Lambda_I^+$ if and only if
\[
\forall i=0,1,\dotsc,k:\;\;\dim_\bC(L^+\cap
\bsF_j^\perp)=i,\;\;\forall j,\; \nu_{i+1}\leq j<\nu_{i},
\]
or, equivalently
\[
\forall i=0,1,\dotsc,k:\;\;\dim_\bC(L^+\cap
\bsF^\nu)=i,\;\;\forall\nu,\; n+1-\nu_i\leq\nu \leq
n-\nu_{i+1},\;\;\nu_0=n+1.
\]
We define $\mu_i$ so that
\[
n-k+i-\mu_i= n+1-\nu_i \Longleftrightarrow \mu_i=\nu_i-(k+1-i),
\]
 and we obtain a partition $\mu_I=(\mu_1\geq \mu_2\geq \cdots\geq \mu_k\geq 0)$.  We deduce that $L^+\in \Sigma_I(\Fl^\bullet)$, where $\Sigma_I(\Fl^\bullet)$ denotes the Schubert cell associated to the partition $\mu$, and the flag $\Fl^\bullet$.

\begin{remark} The partition $(\mu_1,\dotsc, \mu_k)$ can be given a very simple intuitive interpretation.  We describe the set $I$ by placing $\bullet$'s  on the positions $i\in I$, and $\circ$'s on the positions $j\in I^c$. If $I=\{\nu_1>\cdots>\nu_k\}$, then $\mu_i$ is equal to the number of $\circ$'s situated to the left of the $\bullet$ located on the position $\nu_i$. Thus
\[
\mu_{\{k\}}= (k-1,0,\dotsc,0),\;\; \mu_{\{1,\dotsc, k-1,k+1,\dotsc,
n\}}= 1^{n-k}=\underbrace{(1,\dotsc,1)}_{n-k}. \proofend
\]
\end{remark}

A critical lagrangian $\Lambda_I$  is completely characterized  by
its depth $k=\delta(\Lambda_I)=\# I$, and  the associated partition
$\mu$. More precisely,
\begin{equation}
I=\{ \mu_1+ k >\mu_2+ k-1>\cdots >\mu_k+1\}. \label{eq: mu-nu}
\end{equation}
The  Ferres diagram  of the partition $\mu_I$ fits inside a $k\times
(n-k)$ rectangle.

We denote by $\eC_n$ the set
\[
\eC_n=\bigl\{ (m,\mu);\;\;k\in \{0,\dotsc, n\},\;\;\mu\in
\eP_{m,n-m}\,\bigr\}
\]
where $\eP_ {m,n-m}$  is the set of partitions whose Ferres diagrams
fit inside  a $k\times (n-k)$ rectangle.     We have a bijection
\[
\bI_n^+\supset I\mapsto \pi_I =(\# I, \mu_I)\in \eC_n.
\]
For every $\pi=(m,\mu)\in \eC_n$  there exists a unique  $I\subset
\bI_n^+$ such that $\pi_I=(m,\mu)$ and we set
\[
\Lambda_{(m,\mu)}:=\Lambda_I,\;\;W^\pm_{(m,\mu)}:=W^\pm_I.
\]
Observe that
\begin{equation}
\codim_\bR
W^-_{(m,\mu)}=m^2+2|\mu|,\;\;\mbox{where}\;\;|\mu|:=\sum_i\mu_i,
\label{eq: codim-part}
\end{equation}
and
\begin{equation}
\dim_\bR W^-_{(m,\mu)}= n^2-m^2-2|\mu|= (n-m)^2+\dim_\bR\Sigma_\mu.
\label{eq: dim-part}
\end{equation}

The involution $I\mapsto I^c$ on the collection of subsets  of
$\bI_n^+$  is mapped to  the involution
\[
\bI_n^+\supset I\longmapsto\pi_I \in\eC_n\ni,\;\; \pi=(m,\mu)
\mapsto\pi^*:=(n-m, \mu^*)\in \eC_n,
\]
where $\mu^*$ is the transpose of the complement of $\mu$ in the
$k\times (n-k)$ rectangle.   In other words
\[
\pi_{I^c}=\pi_I^*.
\]

\begin{remark} There is a remarkable  involution in this story.  More precisely,
the operator $\bsJ:\whE\ra \whE$  defines a diffeomorphism
\[
\bsJ:\Lag(\whE)\ra\Lag(\whE),\;\;L\mapsto {\bsJ}L.
\]
If we use the depth-partition labelling, then to every pair $\pi=(k,
\mu)\in \eC_n$ we can associate a Lagrangian $\Lambda_{k,\mu}$ and
we have
\[
\bsJ \Lambda_{\pi}=\Lambda_{\pi^*}.
\]
We list some of the  properties of this involution.

\begin{itemize}

\item $f_A({\bsJ}L)= -f_A(L)$, $\forall L\in \Lag(\whE)$, because
$P_{\bsJ L}=\one_{\whE}-P_L$ and $\whA$ and $\tr\whA=0$.

\item $e^{t\whA} \bsJ= \bsJ e^{-t\whA}$ because $\bsJ\whA=-\whA
\bsJ$.

\item $\bsJ L^\pm= ({\bsJ}L)^\mp$, $\forall L\in \Lag(\whE)$.

\item $\bsJ\Lambda_I=\Lambda_{I^c}$,  $\forall I\subset
\bI_n^+$.

\item $\bsJ W^\pm_I= W^\mp_{I^c}$, $\forall I\subset \{1,\dotsc,
n\}$.
\end{itemize}
The  involution is transported by  the diffeomorphism
$\eS:\Lag(\whE)\ra U(E)$ to the  involution $S\mapsto -S$ on
$U(E)$.\qed \label{rem: dual}
\end{remark}

Proposition  \ref{prop: ww} can be rephrased as follows.

\begin{corollary} Let $L\in \Lag(\whE)$ and set $S:=\eS(L)\in U(E)$.  Then the following hold.
\begin{enumerate}
\item[(a)] $L\in W^-_{(m,\mu)}$ if and only if
\[
\dim\ker(\one-S)=m\;\; \mbox{and} \;\;\ker(\one-S)\in
\Sigma_\mu(\Fl^*)\subset\Gr_m(E).
\]
\item[(b)]  $L\in W^+_{(k,\lambda)}$ if and only if
\[
\dim\ker(\one+S)=n-k\;\; \mbox{and}\;\;\ker(\one+S)\in
\Sigma_{{\lambda}^*}(\Fl^*)\subset \Gr^k(E).
\]

\end{enumerate}\qed
\label{cor: ww}
\end{corollary}

Finally, we can give an invariant theoretic description of the
unstable manifolds $W_I^-$.

\begin{definition}      (a)  We define the \emph{symplectic annihilator} of a subspace $U\subset \whE$
to be the subspace $U^\dag:= \bsJ U^\perp$, where $U^\perp$ denotes the orthogonal complement.

(b) A subspace $U\subset  \whE$ is called \emph{isotropic}
(respectively \emph{coisotropic}) if $U\subset U^\dag$ (respectively
$U^\dag \subset U$). (Observe that  a  lagrangian subspace is a
maximal isotropic space.)

(c) An \emph{isotropic flag} of $\whE$ is a collection of isotropic
subspaces
\[
0=\eI_0 \subset \eI_1\subset \cdots \subset \eI_n,\;\;
\dim\eI_k=k,\;\;n=\frac{1}{2}\dim_\bC\whE.
\]
The top space $\eI_n$ is called the \emph{lagrangian subspace
associated to $\eI_\bullet$}. \qed
\end{definition}

Consider the   isotropic  flag $\eI_\bullet$ given by
\[
\eI_\ell =   \Sp \bigl\{ \be_i;\;\;i >n-\ell\,\bigr\}.
\]
If $I=\{\nu_k<\cdots <\nu_1\}$, then $L\in W_I^-$ if and only if
\[
\forall i=0,1,\dotsc,k:\;\;\dim_\bC(L\cap
\eI_\nu)=i,\;\;\forall\nu,\; n+1-\nu_i\leq\nu \leq
n-\nu_{i+1},\;\;\nu_0=n+1.
\]
Define the (real) Borel group
\[
\eB=\eB(\eI_\bullet):=\bigr\{  T\in \Symp(\whE,\bsJ);\;\;
T\eI_\ell\subset \eI_\ell,\;\;\forall \ell\in\bI_n^+\,\bigr\}.
\]

\begin{proposition} The unstable manifold  $W_I^-$  coincides with the $\eB$-orbit of $\Lambda_I$.
\label{prop: un-orb}
\end{proposition}

\proof Observe that $W_I^-$ is $\eB$-invariant  so that $W_I^-$ contains the $\eB$-orbit of $\Lambda_I$.     To prove the converse, we need a better understanding of $\eB$.

 Using the unitary basis $\be_1,\dotsc,\be_n,\bsf_1,\dotsc,\bsf_n$ we can identify $\eB$  with the group of $(2n)\times (2n)$ matrices   $\eT$ which, with respect to the direct sum decomposition $\whE^+\oplus \whE^-$, have the block description
\[
\eT=\left[\begin{array}{cc}
A & AS\\
0 & (A^*)^{-1}
\end{array}
\right],
\]
where $A$ is a lower triangular  invertible $n\times n$ matrix, and $S$ is a  hermitian $n\times n$ matrix.

The Lie algebra of $\eB$ is   is the vector space  $\eX$ consisting of matrices  $X$  of the form
\[
X=\left[ \begin{array}{cc}
\dot{A} & \dot{S}\\
0 & -\dot{A}^*
\end{array}
\right],
\]
where $\dot{A}$ is lower triangular, and $\dot{S}$ is hermitian.  In particular, we deduce that
\[
\dim_\bR \eB = n(n+1)+ n^2 = 2n^2+n.
\]
 Observe that the matrix $\whA$ defining the Morse flow $\Psi^t$  on $\Lag(\whE)$ belongs to the Lie algebra of $\eB$, and     for any open neighborhood $\eN$ of $\Lambda_I$ in $W_I^-$ we have
\[
W_I^-=\bigcup_{t\in \bR} \Psi^t(\eN).
\]
Thus, to prove that $\eB\Lambda_I=W_I^-$, it suffices to show that  the orbit $\eB\Lambda_I$    contains a tiny open neighborhood of $\Lambda_I$ in $W_I^-$.  To achieve this we look at the smooth map
\[
\eB\ra W_I^-,\;\;   g\mapsto g\cdot \Lambda_I,
\]
and  it suffices to show that its  differential at $1\in \eB$ is surjective.

The kernel of this differential is the Lie algebra of the stabilizer  of $\Lambda_I$ with respect to the action of $\eB$. Thus, if we denote by $\St_I$ this stabilizer, it suffices to show that
\[
\dim \eB-\dim \St_I= \dim W_I^-=\bw(I^c).
\]
Observe that $X$ belongs to the Lie algebra of $\St_I$ if and only if the subspace $X\Lambda_I$ is contained in $\Lambda_I$, or equivalently, any vector in $\Lambda_I^\perp=\Lambda_{I^c}$ is orthogonal to $X\Lambda_I$.   If we denote by $\lan \bullet,\bullet\ran$ the hermitian inner product on $\whE$ we deduce that $X$ belongs to the Lie algebra of $\St_I$ if and only if
\[
\lan \be_j, X\be_{i}\ran = \lan \bsf_{i'}, X\be_{i}\ran=\lan \bsf_i, X\bsf_{j'}\ran=\lan \be_{j'} ,X\bsf_j \ran =0,\;\;\forall i,i'\in I, j,j'\in I^c.
\]
If we write $X$ in bloc form
\[
\left[ \begin{array}{cc}
\dot{A} & \dot{S}\\
0 & -\dot{A}^*
\end{array}
\right],\;\; \dot{A}=(\dot{a}^j_i)_{1\leq i,j\leq n},\;\; \dot{S}=(\dot{s}^j_i)_{1\leq i, j\leq n},
\]
then we deduce that $X$ is in the Lie algebra of $\St_I$ if and only if
\[
\dot{a}^j_i= \dot{s}^j_{j'}=0,\;\;\forall i\in I,\;\;j,j'\in I^c.
\]
Suppose $I=\{i_k<\cdots <i_1\}$.  The  equalities
\[
\dot{s}^j_{j'}=0,\;\;j,j'\in I^c
\]
impose $(n-k)^2$      \emph{real} constraints on the matrix $\dot{S}$. For an $i_\ell\in I$, the equalities
\[
\dot{a}^j_{i_\ell}=0, \;\;j\in I^c,\;\;i_\ell<j
\]
impose $(n-i_\ell-\ell+1)$ complex constraints  on $\dot{A}$.  The vector space of lower triangular complex $n\times n$ matrices  has real dimension $n(n+1)$ so that the  Lie algebra of $\St_I$ has real dimension
\[
n(n+1)  -2\sum_{\ell=1}^k (n-i_\ell-\ell+1)+n^2-(n-k)^2=n^2+n-k^2 +2\sum_{\ell=1}^k (i_\ell+\ell-1)
\]
\[
=n^2+n-k +2\sum_{\ell_1}^ki_\ell =n^2+n+\bw(I).
\]
We deduce that
\[
\dim_\bR\eB-\dim_\bR\St_I=  n^2-\bw(I)=\dim W_I^-.\proofend
\]

\begin{corollary}   The collection of unstable manifolds $(W_I^-)_{I\subset \bI_n^+}$ defines a Whitney regular stratification of $\Lag(\whE)$. In particular,  the flow $\Psi^t$ satisfies the Morse-Smale transversality condition.
\end{corollary}

\proof    The  statement  about the Whitney regularity follows immediately   from Proposition \ref{prop: un-orb} and the results of Lander \cite{Land}.      For the reader's convenience, we include an alternate argument.

The unstable  varieties   $W_I^-$  are the   the orbits of a  smooth, semialgebraic action of the semialgebraic group $\eB$ on $\Lag(\whE)$.  If  $W_J^-\subset \cl(W_I^-)$  then, according to  the results of C.T.C. Wall \cite{Wall},  the set $\eR$  of points in $W_J^-$ where the pair $(W_I^-, W_J^-)$  is Whitney regular  is nonempty.  Since $\eB$ acts  by diffeomorphisms  of $\Lag(\whE)$   the set $\eR$ is a $\eB$-invariant subset of $W_J^-$, so it must  coincide with $W_J^-$.

Since the  stratification  by  the unstable manifolds  of the flow $\Psi^t$ satisfies the Whitney regularity condition we deduce from \cite[Thm. 8.1]{Ni2}  that  the flow  satisfies the  Morse-Smale transversality condition.\qed

\section{Tunnellings}
\label{s: 3a}
\setcounter{equation}{0}

 The main  problem we want to investigate in this section is  the structure of \emph{tunnellings}of the flow $\Psi^t=e^{t\whA}$ on $\Lag(\whE)$.   Given $M, K\subset \bI_n^+$, then a \emph{tunnelling}  from $\Lambda_M$ to $\Lambda_K$ is a gradient trajectory
 \[
 t\mapsto \Psi_A^t L=e^{t\whA}L, \;\;L\in \Lag(\whE)
 \]
  such that
 \[
 \lim_{t\ra \infty} \Psi^{-t}_AL=\Lambda_{M} ,\;\; \lim_{t\ra \infty} \Psi^{t}_AL=\Lambda_{K}.
 \]
 We denote by $\eT(M, K)$ the set of tunnellings from $\Lambda_M$ to $\Lambda_K$, and we say  that \emph{$M$ covers $K$}, and write this  $K\prec M$, if $\eT(M,K)\neq\emptyset$. Equivalently, $K\prec M$ if and only if  $W_{M}^-\cap W_{K}^+\neq\emptyset$.  Observe that
 \[
 L\in W^+_K\Longleftrightarrow {\bsJ}L\in  W^-_{K^c}.
 \]
 Hence
 \[
W_{M}^-\cap W_{K}^+= W_{M}^-\cap {\bsJ}W^-_{K^c},
\]
so that
\[
K\prec M \Longleftrightarrow W_{M}^-\cap {\bsJ}W^-_{K^c}\neq
\emptyset.
\]
Let us observe that, although the  flow  $\Psi_A^t$ depends on the
choice of the hermitian operator  $A:\whE^+\ra \whE^+$,  the
equality (\ref{eq: w-}) shows that the unstable manifolds $W_I^-$
\emph{are independent of  the  choice  of $A$}.  Thus, we can choose
$A$ such that
\[
A\be_i=\frac{2i-1}{2}\be_i,\;\;\forall i=1,\dotsc, n.
\]
Using Proposition \ref{prop: self-ind} on self-indexing we obtain
the following result.

\begin{proposition} If $J\prec I$, then $\bw(J)> \bw(I)$ so that $\dim W_J^-<\dim W_I^-$. \qed
\label{prop: monotone}
\end{proposition}

\begin{definition} (a) For any nonempty  set $K\subset \bI_n^+$ of cardinality $k$ we denote  by   $\nu_K$ the unique  strictly decreasing function $\nu_K:\{1, \dotsc ,k\} \ra \bI_n^+$ whose range is $K$, i.e.,
\[
K=\bigl\{\, \nu_K(k)< \cdots <\nu_K(1)\,\bigr\}.
\]
(b)   We define a partial order $\vartriangleleft$ on the
collection of subsets of $\bI_n^+$ by declaring $J\vartriangleleft
I$ if either $I=\emptyset$, or $\#  J\geq \# I$, and for every
$1\leq \ell\leq \# I$ we have $\nu_I(\ell)\leq \nu_J(\ell)$. \qed
\end{definition}

We have the following elementary  fact whose proof is left to the
reader.

\begin{lemma}  Let $K, M\subset \bI_n^+$.  Then the following statements are equivalent.

(a)  $K\vartriangleleft M $

(b)  For ever $ \ell \in \bI_n^+$ we have $\# \bigl(\,K \cap
[\ell,n] \,\bigr)\geq \# \bigl(\,M \cap [\ell, n]\,\bigr)$.

(c)  $M^c\vartriangleleft K^c$.\qed
\end{lemma}

\begin{proposition}  Suppose $K, M\subset \bI_n^+$. Then $K\prec M$ if and only if $K\vartriangleleft M$.
\label{prop: order}
\end{proposition}

\proof Suppose $L\in W^-_M$. Then  $(\bsJ L)^+=\bsJ L^-$,  and we
deduce  that
 \[
 L\in W^-_M\cap W^+_K \Longleftrightarrow \lim_{t\ra \infty}e^{-tA}L^+=\Lambda_M^+\;\;{\rm and} \;\; \lim_{t\ra \infty}e^{-tA}\bsJ L^-=\Lambda^+_{K^c}.
 \]
 In other words,
 \[
 L\in W^-_M\cap W^+_K \Longleftrightarrow  L^+\in \Sigma_M(\Fl^*),\;\;\bsJ L^-\in \Sigma_{K^c}(\Fl^*).
 \]
 We denote by $U^+=U^+_L$ the orthogonal complement of $L^+$ in  $\whE^+$, and by $T= (t_{ij})_{1\leq i,j}$ the Arnold coordinates  of $L$ in the chart $\Lag(\whE)_M$.

 Observe that $U^+$ contains $\bsJ L^-$, the subspace $L^+$ is spanned by the vectors
 \[
\bv_i= \be_i  -\sum_{j\in M^c,\, j>i} t_{ji} \be_j,\;\;i\in M,
\]
and  $U^+$ is spanned by the vectors
\[
\bu_j=  \be_j+\sum_{i\in M, i<j} t_{ij}\be_i=\be_j+\sum_{i\in M,
i<j} \bar{t}_{ji}\be_i,\;\;j\in M^c.
\]
If we write
\[
M^c= \{ j_{n-m} <\cdots <j_1\},\;\;K^c=\{
\ell_{n-k}<\cdots<\ell_1\},
\]
then    the condition $\bsJ L^-\in \Sigma_{K^c}$ is equivalent with
the  existence  of linearly independent vectors of the form
\begin{equation}
\bw_k =\be_k+\sum_{s>k} a_{sk}\be_s,\;\; k\in K^c, \label{eq: orth}
\end{equation}
which span $\bsJ L^-$. Arguing exactly as in \cite[\S 3.2.2]{M} we
deduce that the inclusion
\[
\bsJ L^-= {\rm span}\,\bigl\{\,\bw_k;\;\;k\in K^c\,\bigr\}\subset
U^+={\rm span}\,\bigl\{ \bu_j;\;\;j\in M^c\}
\]
 can happen  only if
\begin{equation}
n-m=\dim U^+\geq \dim L^-= n-k\;\;\mbox{and}\;\; j_i\geq
\ell_i,\;\;\forall i=1,\cdots,\dotsc, n-k,\label{eq: order}
\end{equation}
i.e., $M^c \vartriangleleft K^c$, so that $K\vartriangleleft M$.

Conversely, if (\ref{eq: order}) holds then, arguing as in \cite[\S
3.2.2]{M} we can find  vectors $\bw_k$ as in  (\ref{eq: orth}) and
complex  numbers $\tau_{ij}$ $i\in M$, $j\in M^c$, $i<j$, such that
\[
\Sp\{ \bw_k;\;\;k\in K^c\}\subset \Sp\{ \be_j+\sum_{i\in M, i<j}
\tau_{ij}\be_i;\;\;j\in M^c\}.
\]
Next complete  the collection $(\tau_{ij})$ to a  collection
$(t_{ij})_{1\leq i, j\leq n}$ such that $t_{ij}=\bar{t}_{ji}$,
$\forall i,j$, and $t_{ij}=0$ if $i\in M$ and $j<i$. The collection
$(t_{ij})$  can be viewed as the Arnold coordinates in the chart
$\Lag(\whE)_M$ of a Lagrangian $L\in W^-_M\cap W_K^+$. \qed

\begin{remark} Proposition \ref{prop: order} implies that if  $K\prec M$ and $M\prec N$ then $K\prec N$, so that $\prec$ is a partial order relation. This fact has an interesting consequence.

If $K_0,K_1,\cdots  K_\nu\subset \bI_n^+$ are such that for every
$i=1,\dotsc, \nu$ there exists  tunnelling from $\Lambda_{K_{i-1}}$
to $\Lambda_{K_i}$ then there must exist   tunnelling  from
$\Lambda_{K_0}$ to $\Lambda_{K_\nu}$.\qed
\end{remark}

\begin{proposition}  Suppose $M, K\subset \bI_n^+$.  The the following statements are equivalent.

\begin{enumerate}

\item[(a)] $K\prec M$.

\item[(b)] $W_K^-\subset \cl(W_M^-)$.
\end{enumerate}
\label{prop: strat}
\end{proposition}

\proof The implication (b) $\Rightarrow$ (a) follows from the
above remark. Conversely  assume $K\prec M$. Then  we deduce
$\Lambda_K\subset\cl(W_M^-)$.   Since $\cl(W_M^-)$ is $\eB$ invariant, where $\eB$ is the  real Borel group defined at the end of  Section \ref{s: 3}, we deduce  that $\eB\Lambda_K  \subset \cl(W_M^-)$. We now conclude by invoking Proposition \ref{prop: un-orb}.\qed

\begin{corollary} For any $M\subset \bI_N^+$ we have
\[
\cl(W_M^-)=\bigsqcup_{K\preceq M} W_K^-.\proofend
\]
\label{cor: strat}
\end{corollary}

\begin{corollary}  Let $K, M\subset  \bI_n^+$, and set $k=\# K$, $m=\# M$. The following statements are
equivalent.

\begin{itemize}
\item   $W_K^-\subset \cl(W_M^-)$ and $\dim W^-_K=\dim W^-_M -1$.

\item $\{1\}\in   K$ and $M=K\setminus \{1\}$.

\end{itemize}           \qed
\label{cor: codim1}
\end{corollary}

\begin{corollary}  Let $K, M\subset  \bI_n^+$, and set $k=\# K$, $m=\# M$. The following statements are
equivalent.

\begin{itemize}
\item   $W_K^-\subset \cl(W_M^-)$ and $\dim W^-_K=\dim W^-_M -2$.

\item $k=m$, $\#(K\cap M)=k-1$ and there exists $i\in \{1,\dotsc, n-1\}$ such that $K=(K\cap M)\cup\{i+1\}$ and $M=(K\cap M)\cup\{i\}$.

\end{itemize}           \qed
\label{cor: codim2}
\end{corollary}

\section{Arnold-Schubert cells, varieties and cycles}
\label{s: 4}
\setcounter{equation}{0}

We want to use the results we have proved so far to  describe a very
useful collection  of subsets
 of  $\Lag(\whE)$.    We begin by describing this collection  using the identification $\Lag(\whE)\cong U(\whE^+)$.

 For every complete flag $\Fl_\bullet=\bigl\{\,0=\bsF_0\subset\bsF_1\subset \cdots \subset \bsF_n=E\,\bigr\}$ of $\whE^+$,   and for every subset $I=\{\nu_k<\cdots <\nu_1\}\subset \bI_n^+$, we  denote $\Fl_\bullet^\dag$ the dual flag $\bsF_j^\dag:=\bsF_{n-j}^\perp$, we set $\nu_0=n+1$, $\nu_{k+1}=0$,  and we denote by $\eW_I^-(\Fl_\bullet)$ the set
 \[
\bigl\{ S\in U(\whE^+);\;\;\dim_\bC\bsF_\nu^\dag\cap \ker(\one -S)=
j,\;\;\nu_{j+1}\leq n-\nu <\nu_j,\;\;j=0,\dotsc, k\,\bigr\}
\]
\[
=\bigl\{ S\in U(\whE^+);\;\;\dim_\bC\bsF_\ell^\perp\cap \ker(\one
-S)= j,\;\;\nu_{j+1}\leq \ell <\nu_j,\;\;j=0,\dotsc, k\,\bigr\}.
 \]
We say that $\eW_I(\Fl_\bullet)$ is the \emph{Arnold-Schubert (AS)
cell of type $I$ associated to the flag $\Fl_\bullet$.}  Its
closure, denoted by $\eX_I(\Fl_\bullet)$  is called the \emph{AS
variety of type $I$, associated to the flag $\Fl_\bullet$}. We want
to point out that
\[
S\in \eW^-_I(\Fl_\bullet)\Longrightarrow \dim_\bC\ker(\one-S)=\# I.
\]
If we fix a unitary basis $\underline{\be}=\bigl\{\be_1,\dotsc,
\be_n\,\bigr\}$ of $\whE^+$ we obtain a flag
\[
\Fl_\bullet(\underline{\be}),\;\;\Fl_\nu(\underline{\be}):=\Sp_\bC\bigl\{\be_j;\;\;j\leq
\nu\,\bigr\}.
\]
We set
\[
\eW_I^-(\underline{\be}):=\eW_I^-\bigl(\,\Fl_\bullet(\underline{\be})\,\bigr).
\]
As we know, the unitary symplectic group $U(\whE,\bsJ)\cong
U(\whE^+)\times U(\whE^+)$ acts on $U(\whE^+)$, by
\[
(U_+,U_-)\ast S=U_-SU_+^*,
\]
 and we set
\[
\eW_I^-(\Fl_\bullet, U_+, U_-):= (U_+,U_-)\ast \eW_I^-(\Fl_\bullet).
\]
We denote by $\eX_I(\Fl_\bullet,U_+,U_-)$ the closure of
$\eW_I^-(\Fl_\bullet, U_+, U_-)$. When $I$ is a singleton,
$I=\{\nu\}$, we will use the simpler notation $\eW_\nu^-$ and
$\eX_\nu$ instead of $\eW^-_{\{\nu\}}$ and $\eX_{\{\nu\}}$.  For
every unit complex number $\rho$ we set
\[
\eW_I^-(\Fl_\bullet,\rho)= W^-_I(\Fl_\bullet,\bar{\rho}\one,\one)=
W^-_I(\Fl_\bullet,\one,\rho\one)
\]
\[
=\bigl\{\, S\in U(\whE^+);\;\;\dim_\bC\bsF_\ell^\perp\cap \ker(\rho
-S)= j,\;\;\nu_{j+1}\leq \ell <\nu_j,\;\;j=0,\dotsc, n\,\bigr\}.
\]
When $\Fl_\bullet= \Fl_\bullet(\underline{\be})$    we will use the
alternative notation
\begin{equation}
\eW_I^-(\underline{\be},\rho):=
\eW_I^-\bigl(\,Fl_\bullet(\underline{\be},\rho)\,\bigr),\;\;\eX_I(\underline{\be},\rho)=\eX_I\bigl(\,Fl_\bullet(\underline{\be},\rho)\,\bigr).
\label{eq: ex}
\end{equation}

\begin{ex} Suppose  $\be_1,\dotsc, \be_n$ is an orthonormal basis of $\whE^+$. We form the flag $\Fl_\bullet$ given by
\[
\Fl_j:= \Sp_\bC\{\be_i;\;\;i\leq j\,\}.
\]
For every $\nu\in \bI_n^+$  and every unit complex number  $\rho$ we
have
\[
\eW_\nu^-(\Fl_\bullet,\rho)=\Bigl\{ S\in U(\whE^+);\;\;\exists
z_{\nu+1},\dotsc, z_n\in \bC:\;\;\ker(\rho-S)=\Sp\bigl\{\,
\be_\nu+\sum_{j>\nu}z_j\be_j\,\bigr\}\,\Bigr\}.
\]
Moreover
\[
\eX_\nu(\Fl_\bullet,\rho) =\Bigl\{ S\in U(\whE^+);\;\;
\ker(\rho-S)\cap \Sp\{\be_\nu,\dotsc,\be_n\}\neq
0\,\Bigr\}.\proofend
\]
\label{ex: wnu}
\end{ex}

\begin{definition}  Let $I\subset \bI_n^+$.  We say that a subset $\Sigma\subset \Lag(\whE)=U(E)$ is an \emph{Arnold-Schubert ($AS$)  cell, respectively variety, of type $I$} if there exists a flag $\Fl_\bullet$  of $E$, and $U_\pm\in U(E)$ such that such that  $\Sigma= \eW^-_I(\Fl_\bullet, U_+,U_-)$, respectively $\Sigma= \eX_I(\Fl_\bullet, U_+, U_-)$.  We  will refer to  $\eX_{\nu}$ as the \emph{basic} $AS$ varieties. \qed
\end{definition}

Note that an $AS$  cell of type $I$ is a non-closed, smooth,
semi-algebraic  submanifold  of $\Lag(\whE)$, semialgebraically
diffeomorphic to $\bR^{n^2-\bw(I)}$. The  $AS$ cells can be given a
description    as incidence loci of  lagrangian subspaces of $\whE$.

We denote by $\Flagi(\whE) $ the collection of isotropic  flags of
$\whE$. The  unitary symplectic group
\[
U(\whE,\bsJ)=\bigl\{ T\in U(\whE);\;\;T\bsJ=\bsJ T\,\bigr\},
\]
maps isotropic subspaces to isotropic subspaces and thus acts on
$\Flag_{iso}$. It is easily seen that this action  is transitive.

For any flag $\eI_\bullet\in \Flagi$, and any  subset
\[
I=\{\nu_1>\cdots >\nu_k\}\subset \bI_n^+
\]
we set $\nu_0=n+1$, $\nu_{k+1}=0$, and we define
\[
\eW_I^-(\eI_\bullet):=\bigl\{ \,L\in \Lag(\whE);\;\;\dim (\eI_n\cap
L)=k,\;\; \dim L\cap \eI_\nu=i,\;\;\forall  n-\nu_{i+1}\leq \nu \leq
n+1-\nu_i\,\bigr\}.
\]
If we choose  a  complete flag $\Fl_\bullet$ of $\whE^+$, then the
dual flag $\Fl_\bullet^\dag$ is  an isotropic flag, and we observe
that the  diffeomorphism $\eL: U(\whE^+)\ra \Lag(\whE)$ sends
$\eW_I^-(\Fl_\bullet)$ to $W_I^-(\Fl^\dag_\bullet)$. If
$\underline{\be}$ is a unitary basis, then we ill write
\[
W_I^-(\underline{\be}):=
W_I^-\bigl(\Fl_\bullet(\underline{\be})^\dag\,\bigr)
\]
As we explained  earlier, the unitary symplectic group
$U(\whE,\bsJ)$ is isomorphic to $U(\whE^+)\times U(\whE^+)$, so that
every $T\in U(\whE,\bsJ)$ can be identified with a pair
$(T_+,T_-)\in  U(\whE^+)\times U(\whE^+)$, such that for every $S\in
U(\whE^+)$ we have
\[
T\eL_S=\eL_{T_-ST_+^*}.
\]
We deduce that
\[
\eW_I^-(\Fl_\bullet, T_+,T_-)\stackrel{\eL}{\longleftrightarrow}
T\eW_I^-(\Fl^\bullet).
\]
Since $U(\whE,\bsJ)$ acts transitively on $\Flagi$ we conclude that
any  $AS$ cell is of the form $\eW_I(\eI_\bullet)$ for some  flag
$\eI_\bullet\in \Flagi$.

\smallskip

\ding{45} \emph{In the sequel  we will use the notation $\eW_I^-$
when referring to  $AS$ cells viewed as subsets of the unitary group
$U(\whE^+)$, and the notation $W_I^-$ when referring to $AS$ cells
viewed as subsets of the Grassmannian $\Lag_(\whE)$.}

\smallskip

We would like to associate cycles to the $AS$  cells, and to do this
we must first fix some orientation conventions.  First we need to
fix an orientation on $\Lag(\whE)$ which is orientable because it is
diffeomorphic to the  connected Lie group $U(\whE^+)$.

To fix an orientation on $U(\whE)$ it suffices to pick an
orientation on the Lie algebra $\uu(\whE^+)=T_\one U(\whE^+)$.
This induces an orientation on each tangent space $T_SU(\whE)$ via
the \emph{left translation isomorphism}
\[
T_\one U(\whE^+)\ra T_SU(\whE^+),\;\; T_\one U(\whE^+)\ni X\mapsto
SX\in T_SU(\whE^+).
\]
 To produce such an orientation we first  choose a unitary basis  of $\whE^+$,
\[
\underline{\be}=\bigl\{\,\be_1,\dotsc,\be_n\,\bigr\}.
\]
 We can  then describe any $X\in \uu(\whE^+)$ as a skew-hermitian matrix
\[
X=(x_{ij})_{1\leq i,j\leq n}
\]
We identify $\uu(\whE^+)$ with the space  of  Hermitian operators
$\whE^+\ra \whE^+$, by associating to the  skew-hermitian operator
$X$ the hermitian operator $Z=-\ii X$. Hence $X=\ii Z$, and we write
\[
z_{ij}(X):=-\ii x_{ij}.
\]
 Note that $z_{ii}\in \bR$, $\forall i$, but $z_{ij}$ may not be real if $i\neq j$.  The functions $(z_{ij})_{1\leq i\leq j}$  define linear coordinates on $\uu(\whE^+)$ which via the exponential map define coordinate  in an open neighborhood of $\one$ in $U(\whE^+)$. More precisely, to any sufficiently small  Hermitian matrix $Z=(z_{ij})_{1\leq i,j\leq n}$  one associates the  unitary operator $e^{-\ii Z}$.

 Using the above linear coordinates  we  obtain a decomposition of $\uu(\whE^+)$, as a direct sum between the real  vector space  with coordinates $z_{ii}$, and  a complex vector space  with complex coordinates $z_{ij}$, $i<j$.   The complex summand has a canonical orientation, and we orient the real summand  using the ordered basis
\[
\pa_{z_{11}},\dotsc, \pa_{z_{nn}}.
\]
Equivalently, if we set $\theta^{ij}=\re z_{ij}$, $\vfi^{ij}=\im
z_{ij}$, $\theta^i=z_{ii}$, then the linear functions
$\theta^i,\theta^{ij},\vfi^{ij}:\uu(\whE^+)\ra \bR$ form a basis
of the real dual  of $\uu(\whE^+)$.  The function
$z_{ij}:\uu(\whE^+)\ra \bC$ are $\bR$-linear and we have
\[
\theta^{ij}\wedge \vfi^{ij}= \frac{1}{2\ii} z_{ij}\wedge
\bar{z}_{ij}= \frac{1}{2\ii} z_{ij}\wedge  z_{ji}.
\]
The above orientation  of $\uu(\whE^+)$ is described by the volume
form
\[
\bOm_n= \Bigl(\,\bigwedge_{i=1}^n\theta^i\,\Bigr)\wedge
\Bigl(\,\bigwedge_{1\leq i<j\leq n}\theta^{ij}\wedge
\vfi^{ij}\,\Bigr).
\]
The volume form $\bOm_n$  on $\uu(\whE^+)$ is uniquely determined by
the unitary basis  $\underline{\be}$, and depends continuously on
$\underline{\be}$.  Since the set of unitary bases  is connected, we
deduce that the orientation  determined by $\bom$ is independent of
the choice of the unitary basis $\underline{\be}$. We will refer to
this as the \emph{canonical orientation on the group $U(\whE^+)$}.
Note that when $\dim_\bC\whE^+=1$, the canonical orientation of
$U(1)\cong S^1$ coincides with the  counterclockwise orientation on
the unit circle in the plane.

 We will need to have a  description of this orientation    in terms of Arnold coordinates.       For  a lagrangian $\Lambda\in \Lag(\whE)$ we denote by $\Lag(\whE)_\Lambda$ the Arnold chart
 \[
 \Lag(\whE)_\Lambda=\bigl\{ \, L\in \Lag(\whE);\;\;L\cap \Lambda^\perp\,\bigr\}.
 \]
 The Arnold coordinates identify  this open set with the space $\End^+_\bC(\Lambda)$ of hermitian operators $\Lambda\ra \Lambda$. By choosing a unitary basis of $\Lambda$ we can identify such an operator $A$ with a Hermitian matrix  $(a_{ij})_{1\leq i,j\leq n}$, and  we can  coordinatize  $\End^+_\bC(\Lambda)$ using the   functions $(a_{ij})_{1\leq i \leq j\leq n}$. We the orient  $\Lag(\whE)_\Lambda$ using the form
 \[
(-1)^{n^2}\Bigl(\,\bigwedge_{i=1}^n da_{ii} \,\Bigr)\wedge
\Bigl(\,\bigwedge_{1\leq i<j\leq n} \frac{1}{2\ii} da_{ij}\wedge d
a_{ji}\,\Bigr)
\]
We will refer to this as the \emph{canonical orientation on the
chart $\Lag(\whE)_\Lambda$}. We want to  show that this orientation
convention agrees  with the canonical orientation on the group
$U(\whE)$.

The relationship between the Arnold coordinates    on the chart
$\Lag(\whE)_{\whE^+}=\Lag(\whE)_{\bI_n^+}$ and the  above
coordinates  on $U(\whE^+)$ is given  by the Cayley transform. More
precisely, if $S=e^{\ii Z}$, $Z$ Hermitian matrix,   and $A$ are
the Arnold coordinates of the associated lagrangian $\eL_S$, the
according to (\ref{eq: sym-unit0}) we have
\[
\one +S =2 (1+\ii A)^{-1} \Longleftrightarrow \ii A
=2(\one+S)^{-1}-\one.
\]
To see whether this correspondence is orientation preserving we
compute its differential at $S=\one$, i.e., $A=0$.  We set
\[
S_t:= e^{t\ii Z},  \;\;\;\ii A_t= 2(\one+S_t)^{-1}-\one
\]
we deduce upon  differentiation at $t=0$ that $\dot{A}_0= -2Z$.

Thus, the differential at $\one $ of the  Cayley transform is
represented by a \emph{negative} multiple of  identity matrix in our
choice of coordinates.  This shows that the canonical orientation
on the  chart $\Lag(\whE)_{\whE^+}$  agrees  with the canonical
orientation on the group $U(\whE^+)$.

To show that this happens for any chart $\Lag(\whE)_\Lambda$ we
choose $T=(T_+,T_-)\in U(\whE,\bsJ)$ such that $T\whE=\Lambda$. Then
\[
\Lag(\whE)_\Lambda=   T \Lag(\whE)_{\whE^+},
\]
We fix a unitary basis $\{\be_1, \dotsc, \be_n\}$ of $\whE^+$ and we
obtain unitary  basis  $\be_i'=T\be_i$ of $\Lambda$. Using these
bases we obtain Arnold coordinates
\[
\eA: \Lag(\whE)_{\whE^+}\ra
\End^+_\bC(\bC^n),\;\;\eA':\Lag(\whE)_\Lambda\ra \End^+_\bC(\bC^n).
\]
 Let $L\in \Lag(\whE)_{\whE^+}\cap \Lag(\whE)_\Lambda$.  The Arnold coordinates of $L$  in the chart  $\Lag(\whE)_\Lambda$ are equal to the  Arnold coordinates of $L'=T^{-1}L$ in the chart $\Lag(\whE)_{\whE^+}$, i.e.,
 \[
 \eA'(L)= \eA(T^{-1} L).
 \]
Using (\ref{eq: equivar}) we deduce
\[
\eS_{T^{-1}L}=T_-^* \eS_L T_+=T^{-1}\ast \eS_L.
\]
Form (\ref{eq: sym-unit0}) and (\ref{eq: sym-unit}) we deduce
\[
\eS_L =\eC_{\ii} \bigl(\,\eA(L)\,\bigr):=(\one-\ii \eA(L))(\one+\ii
\eA(L))^{-1},
\]
\[
\ii\eA'(L)=\ii \eA(T^{-1}L)= 2(\one
+T_-^*\eS_LT_+)^{-1}-\one=\eC_{\ii}^{-1}(T_-^*\eS_LT_+).
\]
We seen that the transition map
\[
\End^+_\bC(\bC^n)\ni \eA(L)\mapsto \eA'(L)\in \End^+_\bC(\bC^n)
\]
is  the composition  of the maps
\[
\End^+_\bC(\bC^n)\stackrel{\eC_{\ii}}{\Lra} U(n)\stackrel{
T\ast}{\Lra} U(n)\stackrel{\eC_{\ii}^{-1}}{\Lra} \End^+_\bC(\bC^n).
\]
This composition is orientation preserving if and only if the map
$S\mapsto T\ast S$ is such.  Now we  remark that the map $S\mapsto
T\ast S$ is indeed orientation preserving because it is homotopic to
the identity map since $U(\whE^+)$ is connected.

Fix $I=\{\nu_k<\cdots <\nu_1\}\subset \bI_n^+$, and a unitary basis
of $\whE^+$,
\[
\underline{\be}=\bigl\{\be_1,\dotsc, \be_n\,\bigr\}.
\]
We  want to describe a canonical orientation on the  $AS$ cell
$W_I^-=W_I^-(\underline{\be})$.  We will achieve this by  describing
a canonical co-orientation.

The cell $W_I^-$ is  contained in the   Arnold chart $\Lag(\whE)_I$,
and it is described in the Arnold coordinates $(t_{pq})_{1\leq p\leq
q\leq n}$ on this chart by the system of linearly independent
equations
\[
t_{ji}=0,\;\; i\in I,\;\;j\leq i.
\]
We set
\[
u_{pq}=\re t_{pq},\;\; v_{pq}=\im t_{pq},\;\;\forall 1\leq p <q.
\]
The \emph{conormal bundle}  $T^*_{W_I^-} \Lag(\whE)$ of
$W_I^-\subset \Lag(\whE)$ is the kernel of the natural  restriction
map $T^*\Lag(\whE)|_{W_I^-}\ra T^*W_I^-$. This bundle morphism is
surjective and thus we have a short  exact sequence of bundles over
$W_I^-$,
\begin{equation}
0\ra T^*_{W_I^-} \Lag(\whE)\Lra T^*\Lag(\whE)|_{W_I^-}\Lra
T^*W_I^-\ra 0. \label{eq: con-seq}
\end{equation}
The  $1$-forms $d u_{ji}$, $dv_{ji}$, $dt_{ii}$,  $i\in I$, $j<i$,
trivialize  the conormal bundle. We can orient  the conormal  bundle
$ T^*_{W_I^-}\Lag(\whE)$ using the form
\begin{equation}
\bom_I=(-1)^{\bw(I)} d\bt_I \wedge \Bigl(\bigwedge_{j<i,\;i\in
I}du_{ji}\wedge dv_{ji}\,\Bigr), \label{eq: bomi}
\end{equation}
where $d\bt_I$ denotes the wedge  product of the $1$-forms
$dt_{ii}$, $i\in I$, written in increasing order,
\[
d\bt_I= dt_{\nu_k\nu_k}\wedge \cdots\wedge dt_{\nu_1\nu_1}.
\]
We denote by $\ori_I^\perp$ this co-orientation, and we will refer
to it  as the \emph{canonical co-orientation}. As explained in
Appendix \ref{s: b} this   co-orientation induces a canonical
orientation $\ori_I$ on $W_I^-$. We denote by $[W_I^-,\ori_I^\perp]$
the current of integration thus defined.

To understand how to detect this co-orientation in the unitary
picture we need to   give a unitary description of the  Arnold
coordinates on $\eW_I^-(\underline{\be})\subset U(\whE^+)$.

\begin{definition} Fix a unitary basis $\underline{\be}$ of $\whE^+$. For every subset  $I\subset \bI_n^+$  we define $\eU_I\in U(\whE, \bsJ)$ to be the symplectic unitary operator defined by
\[
\eU_I(\be_k)=\begin{cases} \be_k & k\in I\\
\bsJ \be_k & k\not\in I,
\end{cases},\;\; \eU_I(\bsf_k)=\begin{cases}
\bsf_k & k\in I\\
\bsJ\bsf_k & k\not \in I.
\end{cases}\proofend
\]
\label{def: ui}
\end{definition}
Via the diffeomorphism
\[
U(\whE,\bsJ)\ni T\mapsto (T_+, T_-)\in  U(\whE^+)\times U(\whE^),
\]
the operator $\eU_I$ corresponds to the pair of unitary operators
\[
\eU_I^+=\eT_I,\;\; \eU_I^-=\eT_I^*,
\]
where
\begin{equation}
\eT_I(\be_k)= \begin{cases}
 \be_k & k\in I\\
 \ii\be_k  & k\not\in I.
 \end{cases}
 \label{eq: transit}
 \end{equation}
 Observe that $\eU_I\whE^+= \Lambda_I$, and that the Arnold coordinates $\eA_I$ on $\Lag(\whE)_I$ are related
to the Arnold coordinates on $\Lag(\whE)_{\whE^+}$ via the
equality
 \[
 \eA_I=  \eA\circ\eU_I^{-1}.
 \]
We deduce that if $S\in U(\whE^+)$ is such that
$\eL_S\in\Lag(\whE^+)_I$ then
\[
\eA_I(\eL_S)=\eC_{\ii}(\eU_I^{-1} \ast S) =\eC_{\ii}( \eT_I S\eT_I).
\]

\begin{ex} Let us describe the orientation of $\eW_I^-\subset U(\whE^+)$ at certain special points.  To any  map $\vec{\rho}: I^c\ra S^1\setminus\{1\}$, $j\mapsto \rho_j$,  we associate the  diagonal unitary operator $D=D_{\vec{\rho}}\in\eW_I^-$  defined by
\[
D\be_j=\begin{cases}
1 & j\in I\\
\rho_j & j\in I^c.
\end{cases}
\]
Every tangent vector $\dot{S}\in T_{D}U(\whE^+)$ can be written as
$\dot{S}=\ii D Z$, $Z$ hermitian matrix, so that
\[
Z= -\ii D^{-1} \dot{S}.
\]
The cotangent space  $T^*_{D}U(\whE^+)$   has a natural  basis given
by the $\bR$-linear forms
\[
\theta^p, \theta^{pq}, \vfi^{pq}: T_{S} U(\whE)\ra \bR,\;\;
\theta^p(Z)= ( Z\be_p,\be_p),
\]
\[
\theta^{pq}(Z)=\re (Z\be_q,\be_p),\;\;\vfi^{p,q}=\im (Z\be_q,\be_p).
\]
To describe the orientation of the conormal bundle
$T_{S_I}U(\whE^+)$ we use the above prescription. The Arnold
coordinates on $\eW_I^-$ are given by
\[
\eW_I\ni S\mapsto \eA_I(S):=\eC_{\ii}
(\eT_IS\eT_I)=-\ii(\one-\eT_IS\eT_I)(\one+\eT_IS\eT_I)^{-1}\in
\End^+(\whE^+).
\]
Using the equality
\[
\eC_{\ii}(\eT_IS\eT_I)
=-2\ii\bigl(\one+\eT_IS\eT_I\,\bigr)^{-1}+\ii\one
\]
 we deduce
\[
\frac{d}{dt}|_{t=0} \eA_I(D e^{\ii t Z})= -2\ii\frac{d}{dt}|_{t=0}
\bigr(\,\one+\eT_I De^{\ii t Z}\eT_I\,\bigl)^{-1}
=-2\ii\frac{d}{dt}|_{t=0} \bigr(\,\one+\eT_I D(\one+\ii
tZ)\eT_I\,\bigl)^{-1}
\]
\[
=-2\ii(\one+\eT_ID\eT_I)^{-1}\frac{d}{dt}|_{t=0}
\bigr(\,\one+t\ii\eT_IDZ\eT_I (\one+\eT_ID\eT_I)^{-1}\,\bigr)^{-1}
\]
\[
=-2(\one+\eT_ID\eT_I)^{-1}\eT_IDZ\eT_I
(\one+\eT_ID\eT_I)^{-1}=\dot{A}.
\]
Hence
\[
Z=-\frac{1}{2} D^*\eT_I^* (\one+\eT_ID\eT_I)\dot{A}
(\one+\eT_ID\eT_I)= -\frac{1}{2} D^*\eT_I^* (\one+\eT_I^2D)\dot{A}
(\one+\eT_I^2D).
\]
Note that
\[
\eT_I^2\be_j=\begin{cases}
\be_j & j\in I\\
-\be_j &j\in I^c
\end{cases}
\]
If $i\in I$, then for every $j\leq i$ we have
\[
\lan \dot{Z}\be_j,\be_i\ran =-\frac{1}{2}\bigl(\dot{A}
(\one+\eT_I^2D)\be_j,\eT_I D (\one+\eT_I^2D)\be_i,\bigr)
\]
\[
-\bigl(\dot{A} (\one+\eT_I^2D)\be_j,\be_i)=\begin{cases}
-(\dot{A}\be_j,\be_i) & j\in I\\
\frac{1}{2}(\rho_j-1)(\dot{A}\be_j,\be_i) & j\in I^c.
\end{cases}
\]
If  for $i\in I$ we set
\[
u^i(\dot{A})= (\dot{A}\be_i,\be_i),
\]
\[
u^{ij}(\dot{A})=\re(\dot{A}\be_j,\be_i),\;\;v^{ij}(\dot{A}=\im
(\dot{A}\be_j,\be_i)
\]
then  we deduce
\[
u^i=-\vfi^i
\]
\[
u^{ij}\wedge v^{ij}= k_j \theta^{ij}\wedge \vfi^{ij}
\]
where $k_j$ is the \emph{positive} constant
\[
k_j=\begin{cases}
1 & j\in I\\
\frac{1}{4}|\rho_j-1|^2 & j\in I^c.
\end{cases}
\]
Using (\ref{eq: bomi}) we conclude  that the conormal  bundle to
$\eW_I^-$ is oriented at $D$ by the     exterior monomial
\[
\theta^I \wedge \Bigl(\bigwedge_{j<i,\, i\in I} \theta^{ji}\wedge
\vfi^{ji} \Bigr),
\]
where $\theta^I$ denotes the wedge product of $\{\theta^i\}_{i\in
I}$ written in increasing order.

In particular, if $I=\{\nu\}$ and $D=S_\nu$, i.e., $\rho_j=-1$,
$\forall j\in I^c$, then the conormal orientation of $\eW_\nu^-$ is
given at $S_{\nu}=S_{\{\nu\}}$ by the exterior monomial
\[
\bom^\perp=\theta^\nu \wedge \bigl(\,\theta^{1\nu}\wedge
\vfi^{1\nu}\,\bigr)\wedge \cdots \wedge
\bigl(\,\theta^{\nu-1,\nu}\wedge \vfi^{\nu-1,\nu}\,\bigr).
\]
The tangent space $T_{S_{\{\nu\}}} \eW_{\{\nu\}}^-$ is oriented by
the exterior  monomial
\begin{equation}
\bom_T= (-1)^{\nu-1} \theta^1\wedge\cdots \theta^{\nu-1}\wedge
\theta^{\nu+1}\wedge \cdots \wedge \theta^n \wedge
\Bigl(\bigwedge_{j<k,\,k\neq \nu} \theta^{jk}\wedge \vfi^{jk}\Bigr).
\label{eq: orient1}
\end{equation}
because
\[
\bom^\perp\wedge
\bom_T=\bOm_n=\Bigl(\,\bigwedge_{i=1}^n\theta^i\,\Bigr)\wedge
\Bigl(\,\bigwedge_{1\leq i<j\leq n}\theta^{ij}\wedge
\vfi^{ij}\,\Bigr). \proofend
\]

\label{ex: orient}
\end{ex}

\begin{proposition} We have an equality of currents
\[
\pa[W_I^-,\ori_I]=0.
\]
In other words, using the terminology of Definition \ref{def:
elcyc}, the pair $(W_I^-,\ori_I^\perp)$ is an \emph{elementary
cycle}. \label{prop: cycle}
\end{proposition}

\proof   The proof  relies on the   theory of subanalytic currents
developed by R. Hardt \cite{Hardt2}. For the  reader's convenience
we have gathered in  Appendix \ref{s: b} the basic properties of
such currents.

 Here is  our strategy. We will prove that there exists  an oriented, smooth, subanalytic submanifold $\eY_I$ of $\Lag(\whE)$    with the following properties.

\begin{enumerate}

\item [(a)] $W_I^- \subset \eY_I\subset \cl(W_I^-)=\eX_I$.

\item [(b)]  $\dim(\eX_I\setminus \eY_I) <\dim W_I^- -1$.

\item [(c)] The orientation on $\eY_I$ restricts to the orientation $\ori_I$ on $W_I^-$.

\end{enumerate}

Assuming the existence of such a $\eY_I$ we observe first that,
$\dim\eY_I=\dim W_I^-$, and  that we have an equality of currents
\[
[W_I^-,\ori_I]=[\eY_I,\ori_I].
\]
Moreover
\[
\supp \pa [\eY_I,\ori_I]\subset \cl(\eY_I)\setminus
\eY_I=\eX_I\setminus  \eY_I
\]
so that
\[
\dim \supp \pa [\eY_I,\ori_I] < \dim \eY_I -1.
\]
This proves  that
\[
\pa [W_I^-,\ori_I]=\pa [ \eY_I,\ori_I]=0.
\]
To prove the existence of an $\eY_I$ with the above properties we
recall  that we have a stratification  of $\eX_I$, (see Corollary
\ref{cor: strat})
\begin{equation}
\eX_I=\bigsqcup_{J\prec I} W_J^-, \label{eq: strat}
\end{equation}
where
\[
\dim W_J^-=  \dim W_I^- + \bw(I)-\bw(J).
\]
We distinguish two cases.

\smallskip

\noindent {\bf A.} $1\in I$. In this case,  using  Corollary
\ref{cor: codim1} we deduce  that all the lower strata in the above
stratification have codimension  at least $2$.    Thus, we can
choose $\eY_I=W_I^-$, and the properties (a)-(c) above are trivially
satisfied.

\noindent {\bf B.} $1\not\in I$. In this case, Corollary \ref{cor:
codim1} implies that the stratification  (\ref{eq: strat}) had a
unique codimension  $1$-stratum, $W_{I_*}^-$, where $I_*:=\{1\}\cup
I$. We set
\[
\eY_I:= W_I^-\cup W_{I_*}^-.
\]
We have to prove that this $\eY_I$ has all the  desired properties.
Clearly (a) and (b) are trivially satisfied.  The rest of the
properties follow from  our next result.

\begin{lemma}  The set $\eY_I$ is a smooth, subanalytic, orientable manifold.
\label{lemma: smooth}
\end{lemma}

\proof Consider the Arnold chart $\Lag(\whE)_{I^*}$.   For any $L\in
\Lag(\whE)_{I^*}$ we denote by $t_{ij}(L)$ its Arnold coordinates.
This means that $t_{ij}=\bar{t}_{ji}$ and that $L$ is spanned by the
vectors
\[
\be_i(L)= \be_i+ \sum_{i'\in I_*} t_{i'i}\bsf_{i'}-\sum_{j\in I_*^c}
t_{ji}\be_{ji},\;\;i\in I_*,
\]
\[
\bsf_j(L) =\bsf_j +\sum_{i\in I_*} t_{ji}\bsf_i -\sum_{j'\in I_*^c}
t_{j'j}\be_{j'},\;\; j\in I_*^c.
\]
The $AS$ cell $W_{I^*}^-$ is described by  the equations
\[
t_{ji}=0,\;\;\forall i\in I_*,\;\;j\leq i.
\]
We will prove that
\begin{equation}
W_I^-\cap \Lag(\whE)_{I_*}=\Omega:=\bigl\{ L\in \Lag(\whE)_{I^*}
;\;\;t_{ji}(L) =0,\;\;\forall i\in I,\;\;j\leq i,\;\;t_{11}(L)\neq
0\,\bigr\}. \label{eq: clos}
\end{equation}
Denote by $A\in \End^+_\bC(\whE^+)$ the hermitian operator defined
by $A\be_i =\alpha_i$, $\forall i\in\bI_n^+$, where the real numbers $\alpha_i$ satisfy
\[
0<\alpha_1<\cdots <\alpha_n.
\]
Extend $A$ to  $\whA:\whE\ra \whE$ by setting
$\whA\bsf_i=-\alpha_i\bsf_i$. Note that $L\in W_I^-$ if and only if
\[
\lim_{t\ra \infty} e^{-t\whA}L=\Lambda_I.
\]
Clearly, if $L\in \Omega$, then $ L_t=e^{-t\whA}L$ is spanned by the
vectors
\[
\be_1(L)=\be_1 + t_{11} e^{2\alpha_1 t}\bsf_1 -\sum_{j\in I^c, j\neq
1} t_{j1} e^{(t(\alpha_1-\alpha_j)}\be_j,
\]
\[
\be_i(L_t)=\be_i -\sum_{j\in I_*^c, j>i} t_{ji}
e^{t(\alpha_i-\alpha_j)}\be_i,\;\;i\in I
\]
\[
\bsf_j(L_t)= \bsf_j+\sum_{i\in
I_*,\,i<j}t_{ij}e^{t(\alpha_i-\alpha_j)}\bsf_i-\sum_{j'\in I_*^c}
t_{j'j} e^{-t(\alpha_j+\alpha_{j'})} \be_{j'},\;\;j\in I_c^*.
\]
We note that as $t\ra \infty$ we have
\[
\Sp\{\be_1(L_t)\}\ra \Sp\{\bsf_1\},\;\; \Sp\{\be_i(L_t)\}\ra
\Sp\{\be_i\},\;\;\forall i\in I,
\]
\[
\Sp\{\bsf_j(L_t)\}\ra \Sp\{\bsf_j\},\;\;\forall j\in  I_*^c.
\]
This proves $L_t\ra \Lambda_I$ so that $\Omega\subset W_I^-\cap
\Lag(\whE)_{I_*}$.

Conversely, let $L\in W_I^-\cap \Lag(\whE)_{I_*}$. Then
\[
\lim_{t\ra \infty} L_t =\Lambda_I,\;\;\mbox{where}\;\;L_t=
e^{-t\whA} L.
\]
The space $L_t$ is spanned by the vectors
\[
\be_1(L_t)=\be_1+ t_{11} e^{2\alpha_1 t}\bsf_1+ \sum_{i\in I} t_{i1}
e^{t(\alpha_1+\alpha_i)} \bsf_i -\sum_{j\in I_*^c} t_{j1}
e^{t(\alpha_1-\alpha_j)}\be_j,
\]
\[
\be_i(L_t) =\be_i + t_{1i} e^{t(\alpha_i+\alpha_1)}\bsf_i
+\sum_{i'\in I} t_{i'i} e^{t(\alpha_i+\alpha_{i'})} \bsf_{i'}
-\sum_{j\in I_*^c} t_{ji}e^{t(\alpha_i-\alpha_j)} \be_j,\;\;i\in I,
\]
\[
\bsf_j(L_t)=\bsf_j + t_{1j} e^{t(\alpha_1-\alpha_j)}\bsf_1+
\sum_{i\in I} t_{ij} e^{t(\alpha_i-\alpha_j)}\bsf_i -\sum_{j'\in
I_*^c} t_{j'j} e^{-t(\alpha_j+\alpha_{j'})}\be_{j'},\;\;j\in I_*^c.
\]
Observe that
\[
\be_1(L_t),\bsf_j(L_t)\perp\Sp\bigl\{\,\be_i;\;\;i\in I\,\bigr\}
\subset \Lambda_I,\;\;\forall j\in I_*^c,
\]
and using the condition $L_t\ra \Lambda_I$ we deduce
\[
\Sp\{ \be_1(L_t),\bsf_j(L_t);\;\;j\in I_*^c\}\ra
\Sp\bigl\{\,\bsf_j;\;\;j\in I^c\,\bigr\}\subset \Lambda_I.
\]
On the other hand, the line spanned by $\be_1(L_t)$ converges as
$t\ra \infty$ to either the line spanned by $\be_1$,  or to the line
spanned by $\bsf_i$, $i\in I_*$. Since the line spanned by $\be_1$,
and the line spanned $\bsf_i$, $i\in I$ are  orthogonal to
$\Lambda_I$ we deduce
\[
\Sp\{\be_1(L_t)\}\ra \Sp\{\bsf_1\},
\]
which implies
\[
t_{11}\neq 0,\;\;t_{i1}=0,\;\;\forall i\in I.
\]
Hence
\[
\be_i(L_t) =\be_i +\sum_{i'\in I} t_{i'i}
e^{t(\alpha_i+\alpha_{i'})} \bsf_{i'} -\sum_{j\in I_*^c}
t_{ji}e^{t(\alpha_i-\alpha_j)} \be_j,\;\;\forall i\in I.
\]
Now observe that
\[
\be_i(L_t)\perp \bsf_j,\;\;\forall j\in I^c
\]
and we conclude that
\[
\Sp\bigl\{\be_i(L_t);\;\;i\in I\,\bigr\}\ra
\Sp\bigl\{\be_i;\;\;i\in I\,\bigr\}.
\]
Since $(\be_i(L_t)-\be_i)\perp \be_{i'}$, $\forall i,i'\in I$,
$i\neq i'$, we deduce
\[
\Sp\{\be_i(L_t)\}\ra \Sp\{\be_i\}
\]
which implies
\[
t_{ii'}=0,\;\;t_{ji}=0,\;\;\forall i,i'\in I,\;\;j\in I^c,\;\;j<i.
\]
This proves that $W_I^-\cap \Lag(\whE)_{I_*}\subset \Omega$ and
thus, also the equality (\ref{eq: clos}).   In particular, this
implies  that $\eY_I=W_I^-\cup W_{I_*}^-$  is smooth, because  in
the Arnold chart  $\Lag(\whE)_{I_*}$ which contains the stratum
$W_{I_*}^-$ is  described by the linear equations
\begin{equation}
t_{ji}=0,\;\; i\in I,\;\;j\leq i. \label{eq: cono1}
\end{equation}
To prove that $\eY_I$ is orientable,  we will  construct an
orientation $\ori_{I_*}$ on $\eY_I\cap \Lag(\whE)_{I^*}$ with the
property that its restriction to
\[
\eY_I\cap \Lag(\whE)_I\cap\Lag(\whE)_{I_*}\subset W_I^-
\]
coincides with the canonical orientation $\ori_I$ on $W_I^-$.

We define  an orientation $\ori_{I_*}$  on $\eY_I\cap
\Lag(\whE)_{I_*}$ by orienting the conormal bundle of this
submanifold using the  conormal  volume form
\[
\bom_{I_*} =(-1)^{\bw(I)}d\bt_I\wedge \Bigl(\bigwedge_{i\in
I,\;\;j\leq i}(\frac{1}{2\ii} dt_{ji}\wedge dt_{ij}\,\Bigr).
\]
Let
\[
L\in  \Lag(\whE)_I\cap\Lag(\whE)_{I_*}\subset W_I^-.
\]
We denote by $t_{ij}(L)$ its coordinates in the chart
$\Lag(\whE)_{I_*}$. Then $L$ is spanned by the vectors
\[
\be_1(L) =\be_1+ t_{11}\bsf_1+\sum_{i\in I} t_{i1}\bsf_i -\sum_{j\in
I_*^c} t_{j1}\be_j,\;\;
\]
\[
\be_i(L)= \be_i+ t_{1i}\bsf_1+\sum_{i'\in I}
t_{i'i}\bsf_{i'}-\sum_{j\in I_*^c} t_{ji}\be_{j},\;\;i\in I_*,
\]
\[
\bsf_j(L) =\bsf_j +t_{1j}\bsf_1+\sum_{i\in I} t_{ij}\bsf_i
-\sum_{j'\in I_*^c} t_{j'j}\be_{j'},\;\; j\in I_*^c.
\]
The space $L$ belongs to $\Lag(\whE)_I$ if and only if
\[
L\cap \Lambda_I^\perp= L\cap \Sp\bigl\{\, \be_j,\;\;\bsf_i;\;\;i\in
I,\;\;j\in I^c\}=0.
\]
This is possible if and only if $t_{11}\neq 0$.

We set
\[
\bsf_1'(L)= \frac{1}{t_{11}}\be_1(L)= \bsf_1+
\frac{1}{t_{11}}\be_1+\sum_{i'\in
I}\underbrace{\frac{t_{i'1}}{t_{11}}}_{=:x_{i'1}}\bsf_{i'}
-\sum_{j\in I_*^c}
\underbrace{\frac{t_{j1}}{t_{11}}}_{=:x_{j1}}\be_j
\]
\[
\be_i'(L)=\be_i(L)-t_{1i}\bsf_1(L)
\]
\[
= \be_i+ \sum_{i'\in I}
\underbrace{\Bigl(t_{i'i}-\frac{t_{1i}t_{i'1}}{t_{11}}\Bigr)}_{=:x_{i'i}}\bsf_{i'}-\underbrace{\frac{t_{1i}}{t_{11}}}_{=:x_{1i}}\be_1-\sum_{j\in
I_*^c}\underbrace{\Bigl(t_{ji}-
\frac{t_{1i}t_{j1}}{t_{11}}\Bigr)}_{=:x_{ji}}\be_j,\;\;i\in I
\]
\[
\bsf_j'(L)=\bsf_j(L)-t_{1j}\bsf_1(L)
\]
\[
=\bsf_j+\sum_{i\in I}\underbrace{\Bigl(t_{ij}-\frac{t_{1j}t_{i1}}{t_{11}}\Bigr)}_{=:x_{ij}}\bsf_i-\underbrace{\frac{t_{1j}}{t_{11}}}_{=:x_{1j}}\be_1
-\sum_{j'\in I_*^c}\underbrace{\Bigl(t_{j'j}-\frac{t_{1j}t_{j'1}}{t_{11}}\Bigr)}_{=:x_{j'j}}\be_{j'}, j'\in I_*^c.
\]
The  space $L$ is thus spanned by the vectors $\be_i'(L)$, $i\in
I$ and $\bsf_j'(L)$, $j\in I^c$, where  we  recall that $1\in
I^c$. Also, since $t_{11}=\bar{t}_{11}$ we deduce that
\[
x_{pq}=\bar{x}_{qp},\;\;\forall 1\leq p,q\leq n.
\]
This implies  that  $x_{pq}$  must be  the Arnold coordinates of $L$
in the chart $\Lag(\whE)_I$.

In these coordinates the canonical orientation $\ori_I$  of $W_I^-$
is obtained from the orientation of the conormal bundle given by the
form
\[
\bom_I=(-1)^{\bw(I)}d\boldsymbol{x}_I\wedge
\Bigl(\bigwedge_{j<i,\,i\in I}\frac{1}{2\ii} dx_{ji}\wedge
dx_{ij}\Bigr),
\]
where $d\boldsymbol{x}_I$ denotes the wedge product of the forms $d
x_{ii}$, $i\in I$, in increasing order  with respect to $i$.

Observe that
\[
x_{ii}=t_{ii}-\frac{t_{1i}t_{i1}}{t_{11}},\;\;\forall i\in I,
\]
\[
x_{i1} =\frac{t_{i1}}{t_{11}},\;\; \forall i\in I,
\]
\[
x_{ij}= t_{ij}- \frac{t_{1j}t_{i1}}{t_{11}},\;\;\forall i\in I, j\in
I^c\setminus\{1\}.
\]
Observe that along $W_I^-$ we have
\[
t_{i1}=t_{ij}=0,\;\;\forall i\in I,\;\;j\in I^c,\;\;j<i.
\]
We will  denote  by $O(1)$  any  differential form  on
$\Lag(\whE)_I\cap \Lag(\whE)_{I_*}$ which is a linear combination
of  differential forms of the type
\[
f(t_{p,q}) dt_{p_1q_1}\wedge \cdots \wedge
dt_{p_mq_m},\;\;f|_{W_I^-}=0.
\]
Then
\[
dx_{ii}=dt_{ii}+O(1),\;\;\forall i\in I,
\]
\[
dx_{ij}= dt_{ij} -\frac{t_{1j}}{t_{11}} d t_{i1} +O(1) ,\;\;\forall
i\in I,\;\;j\in I^c,\;\;j\neq 1,
\]
\[
dx_{i1}=\frac{1}{t_{11}} d t_{i1} +O(1).
\]
We deduce that
\[
\bom_I =\frac{1}{t_{11}^{2\#I}} \bom_{I_*} + O(1).
\]
The last  equality shows that the orientations $\ori_I$ and
$\ori_{I^*}$ coincide on the overlap $W_I^-\cap \Lag(\whE)_{I^*}$.
This concludes the proofs of both Lemma \ref{lemma: smooth} and of
the Proposition \ref{prop: cycle}. \qed

\begin{remark}  Arguing as in the first part of Lemma \ref{lemma: smooth} one can prove that for every  $k\in \bI_n^+$ the smooth locus of $\eX_\nu$ contains the strata $\eW_m^-$, $m\geq k$, and $\eW_{\{1,k\}}^-$. In particular, the singular locus of $\eX_k$ has codimension at least $3$ in $\eX_k$.

The codimension $3$ is optimal. For example, the  Maslov variety
$\eX_1\subset \Lag(2)$ is a union of three strata
\[
\eX_1=\eW_1^-\cup \eW_2^-\cup \eW_{\{1,2\}}^-.
\]
The smooth locus is $\eW_1^-\cup \eW_2^-$. The stratum $\eW_2^-$ is one dimensional and  its closure is a smoothly embedded circle. The stratum $\eW_{\{1,2\}}^-$ is zero dimensional. It consists of a point in $\eX_1$ whose link is homeomorphic to a disjoint union of  two  $S^2$-s. One can prove that $\eX_1$ is a $3$-sphere with two  distinct points identified. \qed
\end{remark}

We see that   any  $AS$ cell  $\eW_I(\Fl_\bullet, U_+, U_-)$ defines
a subanalytic cycle in $U(\whE^+)$.  For fixed $I$, any two such
cycle are homologous since    any one of them is the image of
$[\eW_I^-,\ori_I]$ via a   real analytic map,    real analytically
homotopic to the identity. Thus they all determine the same homology
class
\[
\boldsymbol{\alpha}_I\in H_{n^2-\bw(I)}(\Lag(\whE),\bZ),
\]
called the \emph{$AS$ cycle of type $I\subset \bI_n^+$}.  By
Poincar\'e duality we obtain cocycles
\[
\boldsymbol{\alpha}_I^\dag\in H^{\bw(I)}(\Lag(\whE),\bZ).
\]
We will refer to these as \emph{$AS$ cocycles of type $I$}.  When
$I=\{\nu\}$, $\nu\in \bI_n^+$ we will use the simpler notations
$\balph_\nu$ and $\balph_\nu^\dag$ to denote the $AS$ cycles and
cocycles of type $\{\nu\}$. We will refer to these cycles  as the
\emph{basic} $AS$ (co)cycles.

\begin{ex} Observe that the $AS$ cycle $\balph_\emptyset$ is the orientation cycle of $\Lag(\whE)$.

The codimension $1$ basic cycle $\balph_1$ is the  so called
\emph{Maslov cycle}. It defined  by the same incidence  relation
as the Maslov cycle in the case of real  lagrangians, \cite{Ar0}.

The top codimension  basic cycle $\balph_n$ can be identified  with
the integration cycle  defined by the embedding
\[
U(n-1)\hra U(n),\;\; U(n-1)\ni T\mapsto T\oplus  1\in U(n).\proofend
\]
\end{ex}

\section{A transgression formula}
\setcounter{equation}{0}
 \label{s: 5a}

The basic cycles  have a remarkable property. To formulate it we
need to introduce  some fundamental concepts.    We denote by $\eE$
the rank $n=\dim_\bC\whE^+$ complex vector bundle over $S^1\times
\whE$  obtained from the trivial  vector bundle
\[
\whE^+\times [-\pi,\pi]\times  U(\whE^+)\ra  [-\pi,\pi]\times
U(\whE^+),
\]
by identifying the   point $\bu\in \whE^+$ in the fiber over
$(-\pi,g)\in [-\pi,\pi]\times U(\whE^+)$ with the point
$\bv=g\bu\in\whE^+$  in the fiber over $(\pi, g)\in [-\pi,\pi]\times
U(\whE^+)$. Equivalently, consider the  $\bZ$-equivariant  bundle
\[
\widetilde{\eE}=\whE^+\times \bR\times U(\whE)  \ra  \bR\times
U(\whE),
\]
where  the $\bZ$-action is given by
\[
\bZ\times \bigl(\,\whE^+\times \bR\times U(\whE)\,\bigr)\ni (n; \bu,
\theta ,S) \longmapsto  ( S^n \bu, \theta+2n\pi, S)\in \whE^+\times
\bR\times U(\whE),.
\]
Then $\eE$ is the bundle
\[
\bZ\backslash\widetilde{\eE}\ra \bZ\backslash\bigl(\bR\times
U(\whE)\,\bigr).
\]
The sections of this bundle can be identified with  maps $\bu:
\bR\times U(\whE^+)\ra \whE^+$ satisfying the equivariance condition
\[
\bu(\theta+2\pi, S)= S\bu(\theta, S),\;\;\forall (t,S)\in\bR\times
U(\whE^+).
\]
Denote by
\[
\bpi_!: H^\bullet\bigl(\,S^1\times  U(\whE^+),\,\bZ\,\bigr)\ra
H^{\bullet-1}\bigl(\,U(\whE^+),\,\bZ\,\bigr)
\]
 the Gysin map determined  by the natural projection
\[
\bpi:S^1\times  U(\whE^+)\ra U(\whE^+).
\]
For every $\nu=\{1,\dotsc, n\}$ we define $\bgamma_\nu\in
H^{2\nu-1}\bigl( \,U(\whE^+),\bZ\,\bigr)$ by setting
\[
\bgamma_\nu:=\bpi_!  c_\nu(\eE),
\]
where $c_\nu(\whE)\in H^{2\nu}(\whE,\bZ)$ denotes the $\nu$-th Chern
class of $\eE$.

\begin{theorem}[Transgression Formula] For every $\nu=\{1,\dotsc, n\}$ we have the equality
\[
\balph_\nu^\dag=\bgamma_\nu.
\]
\label{th: TP}
\end{theorem}

\proof   Here is briefly the  strategy. Fix a unitary basis
$\underline{\be}=\{\be_1,\dotsc,\be_n\}$ of $\whE^+$,   and consider the
$AS$  variety
\[
\eX_\nu(-1):= \bigl\{ S\in U(\whE^+);\;\;\ker (\one +S)\cap
\Sp_\bC\{ \be_\nu,, \dotsc, \be_n\}\geq 1 \}\,\bigr\}.
\]
It defines a subanalytic cycle
$[\eX_\nu(-1),\ori_\nu]$.  We will prove that there exists a
subanalytic cycle  $\bc$ in $S^1\times U(\whE)$ such that the
following happen.

\begin{itemize}

\item The (integral) homology class determined by $\bc$ is Poincar\'{e} dual to $c_\nu(\eE)$.

\item  We have an equality of  \emph{subanalytic  currents} $\bpi_*\bc= [\eX_\nu(-1),\ori_\nu]$.

\end{itemize}

To construct this  analytic cycle we will use the interpretation of
$c_\nu$ as  the Poincar\'e dual of a degeneracy cycle \cite{HL1,M}.

We set $V:= \Sp_\bC\{ \be_\nu,, \dotsc, \be_n\}$, and we denote by
$\underline{V}$ the  trivial vector bundle with fiber $V$ over
$S^1\times U(\whE)$.   Denote by $\bP(V)$ the    projective space of
lines in $V$, and by $\bp$ the natural projection
\[
\bp: \bP(V)\times S^1 \times (\whE^+)\ra S^1 \times (\whE^+).
\]
We have a tautological line bundle $\eL\ra  \bP(V)\times S^1 \times
(\whE^+)$ defined as the pullback  to $\bP(V)\times S^1 \times
(\whE^+)$ of the tautological line bundle over $\bP(V)$.

To  any  bundle morphism $T: \underline{V} \ra \eE$ we can associate
in a canonical fashion a bundle  morphism $\tilde{T}:\eL\ra
\bp^*\eE$.  We regard $\tilde{T}$ as a section of the bundle
$\eL^*\otimes \bp^*\eE$. If $T$  is a $C^2$, subanalytic section
such that the associated section $\tilde{T}$ vanishes transversally,
then  the zero set $\eZ(\tilde{T})$ is a $C^1$ subanalytic manifold
equipped with a natural orientation   and defines  a subanalytic
current $[\eZ(\tilde{T})]$. Moreover (see \cite[VI.1]{HL1}), the
subanalytic current $\bp_*[\eZ(\tilde{T})]$ is Poincar\'e dual  to
$c_\nu(\eE)$. We will   produce  a  $C^2$, subanalytic  bundle
morphism   $T$ satisfying the above   transversality condition, and
satisfying the additional  equality of currents
\[
\bpi_*\bp_*[\eZ(\tilde{T})]= [\eX_\nu(-1),\ori_\nu].
\]
To construct    such a morphism $T$ we first choose a  polynomial
$\eta\in \bR[\theta]$ satisfying the following conditions
\[
\eta'(\theta)\geq 0, \;\;\forall \theta\in [-\pi,\pi],
\]
\[
\eta(-\pi)=0, \;\;\eta(\pi)=1,\;\;  \eta(0)=\frac{1}{2},
\]
\[
\eta'(0)=\frac{1}{4}, \;\;\eta'(\pm \pi)=\eta''(\pm \pi)=0.
\]
Note that a   bundle  morphism $T :\underline{V} \ra \eE$   is
uniquely  determined by the sections  $T\be_j$, $\nu\leq j\leq n$,
of $\eE$.  Now define a  vector  bundle  morphism
\[
\eT :V\times \bigl(\, [-\pi,\pi]\times U(\whE^+)\,\bigr)\Lra
\whE^+\times \bigl(\, [-\pi,\pi]\times U(\whE^+)\,\bigr)
\]
given by
\[
V\times \bigl(\, [-\pi,\pi]\times U(\whE^+)\,\bigr)\ni (\bv; \theta,
S) \mapsto (S(\theta)\bv; \theta, S)\in \whE^+\times \bigl(\,
[-\pi,\pi]\times U(\whE^+)\,\bigr),
\]
where
\[
S(\theta)= \one + \eta(\theta)(S-\one)=
\bigl(\,1-\eta(\theta)\,\bigr)\one + \eta(\theta) S.
\]
Observe that $S(-\pi)$ is the inclusion of $V$ in $\whE$, while
$S(\pi)=S$. Thus, for every $\bv\in V$  the map
\[
\Psi_{\bv}:[-\pi,\pi]\times U(\whE^+)\ra \whE^+,\;\;
\Psi_{\bv}(\theta,S)= S(\theta)\bv
\]
satisfies
\[
\Psi_{\bv}(\pi, S)=S\Psi_{\bv}(-\pi, S),
\]
and defines a $C^2$- semialgebraic section of $\eE$. Hence $\eT$
determines a $C^2$-semialgebraic  bundle morphism
$T:\underline{V}\ra \eE$.

Now let $(\ell, \theta, S)\in \bP(V)\times S^1 \times U(\whE^+)$
which in the zero set of $\tilde{T}$. This means that   restriction
of $S_\theta$ to the line $\ell\subset V$ is trivial, i.e.,
\[
\ell \subset  \ker  \bigl(\, (\, 1-\eta(\theta)\,)\one +
\eta(\theta) S\,\bigr) .
\]
 Clearly when $\eta(t)=0,1$ this is not possible.  Hence $\eta(\theta)\neq 0,1$ and thus $-\frac{1-\eta(\theta)}{\eta(\theta)}$ must be an eigenvalue of the unitary operator $S$.  Since $\eta(\theta)\in (0,1)$,  and the eigenvalues of $S$ are complex numbers of   norm $1$, we deduce that $-\frac{1-\eta(\theta)}{\eta(\theta)}$  can be an eigenvalue of $S$ if and only if  $-\frac{1-\eta(\theta)}{\eta(\theta)}=-1$, so that $\eta(\theta)=\frac{1}{2}$. From the properties of $\eta$ we conclude that this happens if and only if $\theta=0$.  Thus
 \[
 \eZ(\tilde{T})=\bigl\{\, (\ell, \theta, S)\in \bP(V)\times S^1\times  U(\whE^+);\;\;\theta=0,\;\; \ell\subset \ker (\one +S)\,\bigr\}.
 \]
\begin{lemma}  The section $\tilde{T}$ constructed above vanishes transversally.
\end{lemma}

\proof Let $(\ell_0, 0,S_0)\in \eZ(\tilde{T})$.   Fix $\bv_0\in V$
spanning $\ell_0$. Then we can identify an open neighborhood of
$\ell_0$ in $\bP(V)$ with an open neighborhood of $0$ in the
hyperplane $\ell_0^\perp\cap V$: to any  $\bu\in \ell_0^\perp\cap V$
we associate the line $\ell_{\bu}$ spanned by $\bv_0+\bu$.  We
obtain in this fashion a map
\begin{equation}
(\,\ell_0^\perp\cap V \,) \times (-\pi,\pi)\times U(\whE^+)\ni
(\bu,\theta, S)\stackrel{F}{\longmapsto} \Bigl(\one +
\eta(\theta)(S-\one)\,\Bigr)(\bv_0 +\bu)\in \whE^+, \label{eq: map}
\end{equation}
and we have to prove that  the point $(0,0, S_0)\in(\,\ell_0^\perp\cap V\,)  \times (-\pi,\pi)\times U(\whE^+)$ is a regular point of this
map.

Choose a smooth  path $(-\ve,\ve)\ni t\mapsto (\bu_t,\theta_t,
S_t)\in (\,\ell_0^\perp\cap V \,) \times (-\pi,\pi)\times U(\whE^+)$ such
that
\[
\bu_0=0, \;\;\theta_0=0,\;\; S_{t=0}=S_0.
\]
We set
\[
\dot{\bu}:=
\frac{d}{dt}|_{t=0}\bu_t,\;\;\dot{\theta}:=\frac{d}{dt}|_{t=0}\theta_t,\;\;\dot{S}_0:=\frac{d}{dt}|_{t=0}
S_t
\]
and
\[
X:=S_0^{-1}\dot{S}_0= S_0^*\dot{S}_0
,\;\;\mbox{i.e.},\;\;\dot{S}_0=S_0 X.
\]
Observe that $X$ is a skew-hermitian operator $\whE^+\ra \whE^+$,
and we can identify the tangent space to $\ell_0^\perp\cap V  \times
(-\pi,\pi)\times U(\whE^+)$ at $(0,0,S_0)$ with the space of vectors
\[
(\dot{u},\dot{\theta}, X)\in \ell_0^\perp\cap V\times \bR\times
\uu(\whE^+).
\]
Then
\[
\frac{d}{dt}|_{t=0} F(\bu_t, \theta_t, S_t)=\frac{1}{2}(\one
+S)\dot{u}_0+\eta'(0)\dot{\theta}_0(S_0-\one)(\bv_0)+
\eta(0)\dot{S}_0\bv_0
\]
($-\one\bv_0=S_0\bv_0$, $\eta'(0)=\frac{1}{4}$)
\[
=\frac{1}{2}(\one +S_0)\dot{u}_0+\frac{1}{2}\dot{\theta}_0S_0\bv_0
+\frac{1}{2}S_0X\bv_0=\frac{1}{2}(\one
+S_0)\dot{u}_0+\frac{1}{2}S_0\bigl(\,\dot{\theta}_0\one
+X\,\bigr)\bv_0.
\]
The surjectivity of the differential of $F$ at $(0,0,S_0)$ follows
from the fact that the $\bR$-linear map
\[
\bR\times \uu(\whE^+)\ni(\dot{\theta}_0, X)\longmapsto  \bigl(\,
\dot{\theta}_0+ X\,\bigr)\bv_0\in \whE^+
\]
is surjective for any nonzero vector $\bv_0\in \whE^+$. \qed

The above lemma proves that $\eZ(\tilde{T})$ is a $C^1$ submanifold
of $\bP(V)\times S^1\times U(\whE^+)$. It carries a natural
orientation which we will describe a bit later. It thus defines a
subanalytic current $[\eZ(\tilde{T})]$.  Observe that
\[
\eZ(\tilde{T})\subset  \bP(V)\times \{\theta=0\}\times
U(\whE^+)\subset \bP(V)\times S^1\times U(\whE^+).
\]
The current $\bp_*[\eZ(\tilde{T})]$ is the integration current
defined by   $\eZ(\tilde{T})$ regarded as submanifold of
$\bP(V)\times U(\whE^+)$.   As such, it has the description
\[
\eZ(\tilde{T})=\bigl\{ \,(\ell, S)\in \bP(V)\times
U(\whE);\;\;(\one+S)|_{\ell}=0\,\bigr\}.
\]
We set
\[
\eZ(\tilde{T})^*:=\bigl\{\,(\ell,S)\in
\eZ(\tilde{T});\;\;\ell=\ker(\one +S),\;\;\be_\nu\not\in
\ell^\perp\,\bigr\}.
\]
Note that the projection
\[
\bpi:=\bP(V)\times  U(\whE^+)\ra U(\whE),\;\;(\ell, S)\mapsto S,
\]
maps  $\eZ(\tilde{T})$ surjectively onto $\eX_\nu(-1)$. Moreover,
$\eZ(\tilde{T})^*$ is the preimage under $\bpi$ of the top stratum
$\eW_\nu^-(-1)$ of $\eX_\nu(-1)$,
\[
\eZ(\tilde{T})^*=\bpi^{-1}\bigl(\,\eW_\nu^-(-1)\,\bigr),
\]
and the restriction of $\bpi$ to $\eZ(\tilde{T}^*)$ is a bijection
with inverse
\[
W^-_\nu(-1)\ni S\mapsto (\ker(\one+ S), S)\in \eZ(\tilde{T})^*.
\]
\begin{lemma} The map $\bpi: \eZ(\tilde{T})^*\ra \eW_\nu^-$ is a diffeomorphism.
\label{lemma: diffeo}
\end{lemma}

\proof  It suffices to show that the differential of $\bpi$ is
everywhere injective. Let $\zeta_0=(\ell_0,S_0)\in
\eZ(\tilde{T})^*$.   Suppose $\ell_0=\Sp\{\bv_0\}$.  We have to
prove  that if
\[
(-\ve,\ve)\ni (\ell_t, S_t)\in
\]
a is a smooth path $\eZ(\tilde{T})^*$ passing through $\zeta_0$ at
$t=0$ and $\dot{S}_0:=\frac{d}{dt}|_{t=0}=0$, then
$\frac{d}{dt}|_{t=0}\ell_t=0$.

We write $\ell_t=\Sp\{\bv_0+\bu_t\}$, where $t\mapsto\bu_t\in
\ell_0^\perp\cap V$ is a $C^1$-path such that $\bu_0=0$.  Then
\[
S_t(\bv_0+\bu_t) =-\bv_0-\bu_t,\;\;\forall t,
\]
and  differentiating with respect to  $t$ at $t=0$ we get
\[
-\dot{\bu}_0= S_0\dot{\bu}_0+ \dot{S}_0(\bv_0)= S_0\dot{\bu}_0.
\]
Hence  $\dot{\bu}_0\in \ker(\one+S_0)$. We conclude  that
$\dot{\bu}_0=0$ because  $\ker(\one+S)$ is the line spanned by
$\bv_0$, and $\dot{\bu}_0\perp\bv_0$. \qed

Lemma  \ref{lemma: diffeo} implies that we have an equality of
currents
\[
\bpi_*\bp_*[\eZ(\tilde{T})]=\pm [\eX_\nu(-1)].
\]
To  eliminate the sign ambiguity  we need to   understand the
orientation of $\eZ(\tilde{T})$.

 We begin by describing  the conormal orientation  of  $\eZ(\tilde{T})$ at a  special point $\xi_0=(\ell_0,0, S_0)$, where
 \[
 \ell_0=\Sp_\bC\{\be_\nu\},\;\;\mbox{and}\;\; S_0\be_i= \begin{cases}
 -1 &i=\nu\\
 1 &i\neq \nu.
 \end{cases}
 \]
Observe that $S_0$ is \emph{selfadjoint} and belongs to the top
dimensional stratum $\eW_\nu^-(-1)$ of $\eX_\nu(-1)$. Denote by
$\underline{F}$ the differential at $\xi_0$ of  the map $F$
described in (\ref{eq: map}).

The fiber at $\xi_0$ of the conormal bundle to $\eZ(\tilde{T})$  is
the image of the \emph{real} adjoint of $\underline{F}$,
\[
\underline{F}^\dag : T^*_0\whE^+\ra T^*_{\xi_0}\bigl(\, \bP(V)\times
S^1\times U(\whE^+)\,\bigr).
\]
Since $\underline{F}$ is surjective, its \emph{real} dual
$\underline{F}^\dag$ is injective. The fiber at $\xi_0$ of the
conormal bundle is the image of $\underline{F}^\dag$,  and we have
an orientation on this fiber induced via $\underline{F}^\dag$ by the
canonical orientation of $\whE^+$ as a complex vector space.

The  canonical orientation of the real cotangent space $T^*_0\whE^+$
is  described  by the top degree exterior monomial
\[
\alpha^1\wedge \beta^1\wedge \cdots \wedge \alpha^n\wedge \beta^n,
\]
where $\alpha^k,\beta^k\in \Hom_\bR(\whE^+,\bR)$ are defined by
\[
\alpha^k(\bx)= \re
(\bx,\be_k),\;\beta^k(x)=\im(\bx,\be_k),\;\;\forall
\bx\in\whE^+,\;\;k=1,\dotsc, n.
\]
For every $\dot{u}_0\in V$, $\dot{u}_0\perp\be_\nu$,
$\dot{\theta}_0\in \bR$ and $\ii Z\in \uu(\whE^+)$ we have
\[
\underline{F}^\dag\alpha^k(\dot{\bu}_0,\dot{\theta}_0,  \ii Z)
=\re\bigl( \,\underline{F}(\dot{u}_0,\dot{\theta}_0,  \ii Z)\,
,\,\be_k\,\bigr)
\]
\[
=\frac{1}{2}\re \bigl(\,
(\one+S_0)\dot{u}_0,\be_k\,\bigr)+\frac{1}{2}\re
\bigl(\,S_0(\dot{\theta}_0+\ii Z)\be_\nu,\be_k\,\bigr)
\]
\[
=\frac{1}{2}\re \bigl(\,
\dot{u}_0,(\one+S_0)\be_k\,\bigr)+\frac{1}{2}\re
\bigl(\,(\dot{\theta}_0+\ii Z)\be_\nu,S_0\be_k\,\bigr).
\]
\[
\underline{F}^\dag\beta^k(\dot{\bu}_0,\dot{\theta}_0,  \ii Z)
=\frac{1}{2}\im\bigl(\,
\dot{u}_0,(\one+S_0)\be_k\,\bigr)+\frac{1}{2}\im
\bigl(\,(\dot{\theta}_0+\ii Z)\be_\nu,S_0\be_k\,\bigr)
\]
To simplify the final result observe that the restrictions to
$\ell_0^\perp\cap V$ of the  $\bR$-linear functions
$\alpha^k,\beta^k$, $k\geq \nu$ determine a basis of
$\Hom(\ell_0^\perp\cap V,\bR)$ which we will continue by the same
symbols.  We denote by $dt$ the  tautological linear map $T_0\bR\ra
\bR$.

Recall  (see Example \ref{ex: orient}) that the real  dual of
$\uu(\whE^+)$   admits a natural basis given by the  $\bR$-linear
forms
\[
\theta^j(Z)=(Z\be_j,\be_j),\;\;
\theta^{ij}(Z)=\re(Z\be_j,\be_i),\;\;
\vfi^{ij}=\im(Z\be_j,\be_i),\;\;i<j\in \bI_n^+,\;\;\ii Z\in
\uu(\whE^+).
\]
Observe that
\[
\theta^{ij}=\theta^{ji},\;\;\vfi^{ij}=-\vfi^{ji},\;\;\forall i\neq
j.
\]
For every $\dot{u}_0\in \ell_0^\perp\cap V$,  $\dot{\theta}_0\in
T_0\bR$, $\ii Z\in \uu(\whE^+)$  we have
\[
\re \bigl(\, \dot{u}_0,(\one+S_0)\be_k\,\bigr)= \begin{cases}
 2\alpha^k(\dot{u}_0) & k>\nu\\
0 & k\leq \nu,
\end{cases}
\]
\[
\im \bigl(\, \dot{u}_0,(\one+S_0)\be_k\,\bigr)=\begin{cases}
 2\beta^k(\dot{u}_0) & k>\nu\\
0 & k\leq \nu,
\end{cases}
\]
\[
\re \bigl(\,(\dot{\theta}_0+\ii Z)\be_\nu,S_0\be_k\,\bigr)=
\delta_{\nu k}dt(\dot{\theta}_0)-
\begin{cases}
0  & k=\nu\\
\vfi^{k\nu}(Z) & k\neq \nu,
\end{cases}
\]
\[
\im \bigl(\,(\dot{\theta}_0+\ii Z)\be_\nu,S_0\be_k\,\bigr)=
\begin{cases}
\theta^\nu(Z)  & k=\nu\\
\theta^{k\nu}(Z) & k\neq \nu,
\end{cases}
\]
We deduce  the following.

\begin{itemize}

\item If $k <\nu$ then
\[
\underline{F}^\dag
\alpha^k=-\frac{1}{2}\vfi^{k\nu},\;\;\underline{F}^\dag\beta^k=\frac{1}{2}\theta^{k\nu}
\]
\item If $k=\nu$ then
\[
\underline{F}^\dag \alpha^\nu= \frac{1}{2}dt,\;\;
\underline{F}^\dag\beta^\nu= \frac{1}{2}d\theta^\nu
\]
\item If $j>\nu$ then
\[
\underline{F}^\dag \alpha^j=\alpha^j +\frac{1}{2}\vfi^{\nu
j},\;\;\underline{F}^\dag\beta^k=\beta^j+\frac{1}{2}\theta^{\nu j}.
\]

\end{itemize}
Thus, the conormal  space   of $\eZ(\tilde{T})\hra \bP(V)\times S^1
\times U(\whE)$ at $\xi_0$  has an orientation given by the oriented
basis
\[
-\vfi^{1\nu},\theta^{1\nu},\dotsc,
-\vfi^{\nu-1,\nu},\theta^{1,\nu},dt,d\theta^\nu, \alpha^{\nu+1}
+\frac{1}{2}\vfi^{\nu, \nu+1},  \beta^{\nu+1}
+\frac{1}{2}\theta^{\nu, \nu+1},\dotsc, \alpha^{n}
+\frac{1}{2}\vfi^{\nu, n},  \beta^{n} +\frac{1}{2}\theta^{\nu, n}
\]
which is equivalent with the orientation given by the oriented basis
\[
\theta^{1\nu},\vfi^{1\nu},\dotsc, \theta^{1,\nu},
\vfi^{\nu-1,\nu},dt,d\theta^\nu, \alpha^{\nu+1}
+\frac{1}{2}\vfi^{\nu, \nu+1},  \beta^{\nu+1}
+\frac{1}{2}\theta^{\nu, \nu+1},\dotsc, \alpha^{n}
+\frac{1}{2}\vfi^{\nu, n},  \beta^{n} +\frac{1}{2}\theta^{\nu, n}.
\]
We  will  represent this oriented basis by the exterior polynomial
\[
\bom^{\rm \footnotesize{norm}}\in \Lambda^{2n} T^*_{\xi_0}
\bigl(\,\bP(V)\times S^1\times U(\whE^+)\,\bigr),
\]
\[
\bom^{\rm \footnotesize{norm}}:= \Bigl(\bigwedge_{k<\nu}
\theta^{k\nu}\wedge \vfi^{k\nu}\Bigr)\wedge dt\wedge  d\theta^{\nu}
\wedge \Bigl(\bigwedge_{j=\nu+1} (\alpha^j+\frac{1}{2}\vfi^{\nu
j})\wedge (\beta^j+\frac{1}{2}\theta^{\nu j})\Bigr).
\]
The zero set $\eZ(\tilde{T})$ is a smooth manifold of dimension
\[
\dim_\bR\eZ(\tilde{T}) =\dim_\bR\bigl(\,\bP(V)\times S^1 \times
U(\whE^+)\,\bigr)-\dim_{\bR} \whE^+
\]
\[
=  (2n-2\nu)+ 1+n^2- 2n= n^2-(2\nu-1)=n^2-\bw(\nu).
\]
The orientation   of $T_{\xi_0}\bigl(\, \bP(V)\times S^1\times
U(\whE^+)\,\bigr)$ is described by the exterior monomial
\[
\bOm= dt \wedge \Bigl(\bigwedge_{j=\nu+1}^n \alpha^j\wedge
\beta^j\Bigr)\wedge \underbrace{\Bigl(\bigwedge_{i=1}^n \theta^i
\Bigr)\wedge \Bigl(\bigwedge_{j<i} \theta^{ji}\wedge
\vfi^{ji}\Bigr)}_{\bOm_n}
\]
The orientation of $T_{\xi_0}\eZ(\tilde{T})$ is given by any
$\bom\in \Lambda^{n^2-\bw(\nu)} T_{\xi_0}^*\bigl(\bP(V)\times
S^1\times U(\whE^+)\,\bigr)$ such that
\[
\widehat{\bOm} = \bom^{\rm norm}\wedge \bom.
\]
We can take  $\bom$ to be
\begin{equation}
\bom=\bom_{{\rm tan}} :=(-1)^{\nu-1} \theta^1\wedge\cdots
\theta^{\nu-1}\wedge \theta^{\nu+1}\wedge \cdots \wedge \theta^n
\wedge \Bigl(\bigwedge_{j<k,\,k\neq \nu} \theta^{jk}\wedge
\vfi^{jk}\Bigr). \label{eq: orient2}
\end{equation}
If we now think of $\eZ(\tilde{T})$ as an \emph{oriented}
submanifold of $\bP(V)\times U(\whE)\subset \bP(V)\times S^1\times
U(\whE^+)$, we see that its conormal bundle   has a natural
orientation given by the  exterior form
\[
\bom_0^{\rm norm}= \Bigl(\bigwedge_{k<\nu} \theta^{k\nu}\wedge
\vfi^{k\nu}\Bigr)\wedge  d\theta^{\nu} \wedge
\Bigl(\bigwedge_{j=\nu+1} (\alpha^j+\frac{1}{2}\vfi^{\nu j})\wedge
(\beta^j+\frac{1}{2}\theta^{\nu j})\Bigr).
\]
The  discussion in Example \ref{ex: orient} shows that the tangent
space $T_{S_0}\eW_\nu^-(-1)\subset T_{S_0}U(\whE^+)$ is also
oriented by $\bom_{\rm tan}$.  This proves that the differential $D\bpi: T_{(\ell_0, S_0)} \eZ(\tilde{T})\ra T_{S_0}\eW_\nu^-(-1)$ is orientation preserving.
This concludes the proof of the Theorem \ref{th: TP}. \qed

\begin{remark} (a)  The proof  of  Theorem \ref{th: TP}  shows that we have a resolution  $\widetilde{\eX}_\nu\stackrel{\pi}{\ra}\eX_\nu$ of  $\eX_\nu$, where
\[
\widetilde{\eX}_\nu=\bigl\{ (\ell, S)\in \bP(V_\nu)\times
U(\whE^+);\;\;(\one -S)\mid_\ell=0\,\bigr\},
\]
and $\pi$ is induced by the  natural projection $\bP(V_\nu)\times
U(\whE)\ra U(\whE^+)$. Here $V_\nu:=\Sp_\bC\bigl\{\be_\nu,\dotsc,
\be_n\,\bigr\}$. We call the map $\pi$ a resolution, because it is
semi-algebraic,  proper, and it is a diffeomorphism over the top
dimensional stratum  $\eW_\nu^-$ of $\eX_\nu$.  The map $\pi$ is
also a Bott-Samelson  cycle (see \cite{Duan, PT} for a definition)
for the  Morse function
\[
U(\whE^+)\ni S \mapsto -\re\tr AS\in \bR
\]
and its critical point $S_\nu\in U(\whE^+)$ given by
\[
S_\nu\be_k=\begin{cases}
\be_\nu & k=\nu\\
-\be_k & k\neq \nu.
\end{cases}
\]
All the $AS$-varieties $\eX_I$ admit such Bott-Samelson
resolutions (see \cite{PT}),  $\rho_I:\tilde{X}_I\ra \eX_I$, and
it is essentially trough these resolutions  that the cycles
$\balph_I$ were defined by  Vasiliev \cite{Vas}.

(b)  We have a natural group morphism
\[
\pi_{2n-1}(U(n))\ra   \bZ,\;\;\; \bigl(S^{2n-1}\stackrel{f}{\Lra}
\bZ\,\bigr)\longmapsto \int_{S^{2n-1}} f^*\balph_n^\dag
=f_*[S^{2n-1}]\bullet \balph_n.
\]
Using Bott  divisibility theorem,\cite{Bott}, \cite[Thm.
24.5.2]{Hirz}   and the Theorem \ref{th: TP} we deduce that this is
an injective morphism whose range is the subgroup $(n-1)!\bZ\subset
\bZ$.

(c) The  rank $n$ complex  vector  bundle $\eE\ra  S^1\times U(n)$
has an interesting homotopic theoretic    significance. Its
classifying map $F: S^1\times U(n)\ra BU(n)$ can be  viewed by adjunction as a map  $F_1:U(n) \ra {\rm Map}\,(S^1, BU(n))$.
The space ${\rm Map}\,(S^1, BU(n))$ is the total space of a fibration
\[
\Omega BU(n)\hra {\rm Map}\,(S^1, BU(n))\ra BU(n).
\]
The fiber  of this fibration is homotopic to  $U(n)$. The map  $F_1$
will sent $U(n)$ to a fiber of this fibration, and will be  a
homotopy equivalence $U(n)\ra \Omega BU(n)$. To understand this map
we use the adjunction
\[
{\rm Map}\,(U(n), \Omega BU(n))\cong {\rm Map}\,(\Sigma G,  BU(n))
\]
so the map $F_1$ corresponds to a map $F_2:\Sigma G\ra BU(n)$ which
classifies a rank $n$ complex vector bundle $\eE_2$ over the
suspension $\Sigma U(n)$. This bundle is obtained  via the clutching
construction, where the clutching  over the ``Equator'' $U(n)\hra
\Sigma U(n)$ is given by the identity map $U(n)\ra U(n)$.
Equivalently, the  map $\Sigma U(n)\ra B U(n)$   is  a special case
of the inclusion (see \cite{Sta}) $\Sigma G \hra BG$ for any compact
Lie group.

Finally, we want to point out that, when $n=1$, $\eE$ is the degree
$1$ line bundle over the torus $S^1\times U(1)$.\qed
\end{remark}

\section{The Morse-Floer complex and intersection theory}
\label{s: 5}
\setcounter{equation}{0}
It is well known that the integral cohomology ring of $U(n)$ is an
exterior algebra    freely generated by elements $x_i\in
H^{2i-1}(U(n),\bZ)$, $i=1,\dotsc, n$. The transgression formula
implies that as generators $x_i$ of this ring we can take the $AS$
cocycles $\balph_i^\dag$. In this section we would like to prove
this   by direct geometric  considerations,  and then investigate
the cup product   of  two arbitrary $AS$ cocycles.

\begin{proposition} The $AS$-cycles $\balph_I$, $I\subset \bI_n^+$, form a $\bZ$-basis of $H_\bullet\bigl(\,\Lag(\whE)\,,\bZ\,\bigr)$.
\label{prop: basis}
\end{proposition}

\proof    We  will use a Morse theoretic approach.  Consider again
the  Morse flow $\Psi_A^t=e^{t\whA}$ on  $\Lag(\whE)$.

\begin{lemma} The flow $\Psi^t$ is a Morse-Stokes flow, i.e., the following hold.

\noindent (a) The flow $\Phi_t$ is a finite volume flow,  i.e., the
$(n^2+1)$-dimensional  manifold
\[
\Bigl\{ \,\bigl(\,t, \Psi^{1/t}(L), L\,\bigr);\;\;t\in
(0,1],\;\;L\in \Lag(\whE)\,\Bigr\}\subset (0,1]\times
\Lag(\whE)\times \Lag(\whE),
\]
has finite    volume.

\noindent (b)  The  stable and unstable manifold $W_I^\pm$ have
finite volume.

\noindent (c) If there exists a tunnelling  from $\Lambda_I$ to
$\Lambda_J$ then $\dim W_J^- <\dim  W_I^-$. \label{lemma: mor-stok}
\end{lemma}

\proof From  Theorem \ref{th: tame-flow} we deduce that  $\Psi^t$ is
a tame flow.  Proposition \ref{prop: tamefl}  now implies that the
flow  satisfies (a) and (b).    Property (c) follows from
Propositiob \ref{prop: strat}. \qed

As in  Harvey-Lawson  \cite{HL},  we consider   the subcomplex
$C_\bullet(\Psi^t)$ of the complex
$\eC_\bullet\bigl(\,\Lag(\whE)\,\bigr)$ of subanalytic chains
generated by the    analytic chains $[W_I^-,\ori_I]$, and their
boundaries.  According to \cite[Thm. 4.1]{HL}, the inclusion
\[
C_\bullet(\Psi^t)\hra \eC_\bullet
\]
induces an isomorphism in homology.

 Proposition \ref{prop: cycle} implies that the complex $C_\bullet(\Psi^t)$ is perfect. Hence the $AS$ cycles, which form an integral basis of the complex $C_\bullet(\Psi^t)$, also form an integral basis of the integral homology of $\Lag(\whE)$.\qed

\begin{remark} (a) The complex $C_\bullet(\Psi^t)$ is  isomorphic to the  Morse-Floer complex of the flow $\Phi^t$, \cite[\S 2.5]{Ni1}. \qed
\end{remark}

Using the  Poincar\'e duality  on $U(\whE^+)$ we  obtain
intersection products
\[
\bullet: H_{n^2-p}\bigl(\,U(\whE^+),\bZ\bigr)\times
H_{n^2-q}\bigl(\, U(\whE^+),\bZ\,\bigr)\ra
H_{n^2-p-q}\bigl(\,U(\whE^+),\bZ\,\bigr).
\]
For every pair of  nonempty, disjoint subsets $I, J\subset \bI_n^+$
such that
\[
I=\bigl\{ i_1< \cdots <i_p\bigr\},\;\; J= \bigl\{ j_1<\cdots
<j_q\,\bigr\},
\]
 we  define $\eps(I,J)=\pm 1$ to be the signature of the permutation $(i_1,\dotsc, i_p;j_1,\dotsc, j_q)$ of $I \cup J$.

\begin{proposition} Let $I, J\subset \bI_n^+$ such that $\bw(I)+\bw(J)=\bw(\bI_n^+)=n^2$. then
\[
\alpha_I\bullet \alpha_J = \begin{cases}
0 & I\cap J\neq \emptyset\\
\eps(I,J)  &  I=J^c.
\end{cases}
\]
\label{prop: pd}
\end{proposition}

\proof Fix unitary basis  $\{\be_1,\dotsc, \be_n\}$ of $\whE^+$, and
consider the symmetric operator $A_0:\whE^+\ra \whE^+$ given by
\[
A_0\be_i=\frac{2i-1}{2}\be_i.
\]
We   form as usual the associated symmetric operator $\whA_0:\whE\ra
\whE$, and the  positive gradient flow $e^{t\whA_0}$ on $\Lag(\whE)$
associated to the Morse function
\[
\vfi_0: \Lag(\whE)\ra \bR,\;\; L\mapsto \re\tr (\whA_0 P_L).
\]
For every  critical point $\Lambda_K$ of $\vfi_0$ we have
\[
\dim W_K^-=\bw(K)= \vfi_0(\Lambda_K)+\frac{n^2}{2}.
\]
For every $M\subset \bI^+_n$ we denote  by $W_M^+$ the stable
manifold at $\Lambda_M^+$.

Let $\bw(I)+\bw(J)=\bw(\bI_n^+)$.  Using the equality
\[
W_{J^c}^+=\bsJ W_J^-
\]
we deduce  that $W_{J^c}^+$ is also an $AS$ cell  of type $J$, so
that  that we can  represent the  homology class $\balph_J$    by
the subanalytic cycle given as the  integration  over the  stable
manifold $W_{I^c}^+$ equipped with the orientation  induced by the
diffeomorphism $\bsJ: W_J^-\ra W^+_{J^c}$. We denote by
$\eX_{J^c}^+$ its closure.

We have $\vfi_0(\Lambda_{J^c})=-\vfi_0(\Lambda_J)$, and the
equality $\bw(I)+\bw(J)=n^2$ translates into the equality
\[
\vfi_0(\Lambda_{J^c})=\vfi_0(\Lambda_I)=:\kappa.
\]
Observe that,
\[
\eX_{J^c}^+\setminus \{\Lambda_{J^c}\} \subset \bigl\{ \vfi_0
>\kappa\},\;\;\eX_I^-\setminus \{\Lambda_I\}\subset \{\vfi_0 <\kappa\}.
\]
This shows that  if $I^c\neq J$  and $\bw(I^c)=\bw(J)$ the supports
of the subanalytic currents  $[\eX_{J^c}^+]$ and $[\eX_I^-]$ are
disjoint, so that, in this case,
\[
\alpha_I\bullet\alpha_J=0.
\]
When $J=I^c$  we see that the supports of the above subanalytic
cycles intersect only at $\Lambda_I$.  In fact,  only the top
dimensional strata  of their supports intersect, and they do so
transversally.   Hence  the intersection of the analytic cycles
$[\eX_{J^c}^+]$ and $[\eX_I^-]$ is well defined, and from
Proposition \ref{prop: inter}   we deduce
\[
[\eX_{J^c}^+]\bullet [\eX_I^-]=\pm [\Lambda_I],
\]
where  $[\Lambda_I]$ denotes  the Dirac $0$ dimensional current
supported at $\Lambda_I$.   The fact that the correct choice of
signs is $\eps(I,I^c)$ follows from our orientation conventions.\qed

From the  above result we  deduce that for every cycle $c\in
H_k(U(n),\bZ)$ we have a decomposition
\begin{equation}
c=\sum_{\bw(I)=k} \eps(I,I^c)(c\bullet \balph_{I^c})\balph_I.
\label{eq: decomp}
\end{equation}

\begin{theorem}[Odd Schubert calculus]   If $I=\{i_k<\cdots <i_1\}\subset \bI_n^+$ then
\begin{equation}
\balph_I=\balph_{i_k}\bullet\cdots \bullet \balph_{i_1}, \label{eq:
schubert}
\end{equation}
or equivalently,
\begin{equation}
\balph_I^\dag=\balph_{i_k}^\dag\wedge \cdots \wedge
\balph_{i_1}^\dag, \label{eq: schubert'}
\end{equation}
\label{th: schubert}
\end{theorem}

\proof Let us first describe our strategy.  Fix a unitary basis
$\underline{\be}$ of $\whE^+$,  an injection $\rho: I\ra
S^1\setminus \{\pm 1\}$, $i\mapsto \rho_i$  and we consider the $AS$
varieties $\eX_{i}(\rho_i)=\eX_i(\rho_i,\underline{\be})$, $i\in I$,
defined in (\ref{eq: ex}).  We denote by $[\eX_{i}(\rho_i)]$ the
associated   subanalytic cycles.    We will  prove the following
facts.

\begin{enumerate}

\item[{\bf A.}] The varieties  $\eX_{i_\ell}$ intersect quasi-transversally, i.e., for any subset $J\subset I$ we have
\[
\codim \bigcap_{j\in J}\eX_j(\rho_j) \geq \bw(J).
\]
\item [{\bf B.}] There exists a continuous  semialgebraic map  $\Xi:U(\whE^+)\ra U(\whE^+)$, semialgebraically homotopic to the identity such that
\[
\Xi\Bigl(\bigcap_{i\in I}\eX_i(\rho_i) \,\Bigr)\subset \eX_I=
\eX_I(1).
\]
\item[{\bf C.}] The intersection current $[\eX_{i_k}(\rho_{i_k})]\bullet\cdots\bullet[\eX_{i_1}(\rho_{i_1})]\bullet[\eX_{I^c}(1)]$ is  a well defined  zero dimensional  subanalytic current consisting of a single point with multiplicity $\eps(I,I^c)$.
\end{enumerate}
We claim that the above  facts imply (\ref{eq: schubert}).  To see
this, note first that  {\bf A}  implies that,   according to
\cite{Hardt} (see also  Appendix \ref{s: b}), we can form the
intersection current
\[
\eta=[\eX_{i_k}(\rho_{i_k})]\bullet \cdots
\bullet[\eX_{i_1}(\rho_{i_1})].
\]
The current $\eta$ is a subanalytic current  whose  homology class
is $\balph_{i_k}\bullet\cdots \bullet \balph_{i_1}$, and its support
is
\[
\supp(\eta)=  \Bigl(\bigcap_{i\in I}\eX_i(\rho_i)\,\Bigr).
\]
The pushforward  $\Xi_*(\eta)$ is also a subanalytic current and it
represents the same homology class  since $\Xi$ is homotopic to the
identity. Moreover, property {\bf B} shows that
\[
\supp \Xi_*(\eta) \subset \eX_I(1).
\]
Consider again the  dual  $AS$ varieties  $\eX^+_J$,
$\bw(J)=\bw(I)$. In the proof of Proposition \ref{prop: pd} we have
seen that
\[
\eX_I\cap \eX_J^+=\emptyset,\;\;\mbox{if $J\neq I$}.
\]
Hence, the equality (\ref{eq: decomp})  implies that there exists an
integer $k=k(I)$ such that
\[
\balph_{i_k}\bullet\cdots \bullet \balph_{i_1}=k_I\balph_I,
\]
where
\[
k_I=\eps(I,I^c)\bigl(\, \balph_{i_k}\bullet\cdots \bullet
\balph_{i_1}\,\bigr)\bullet\balph_{I^c}.
\]
The equality $k_I=1$ now follows from {\bf C}.

\medskip

\noindent {\bf Proof of  A.} Since the set of unitary operators with
simple eigenvalues is open  and dense, we deduce that  the set
\[
\eY_J:=\bigcap_{j\in J}\eX_j(\rho_j)
\]
 contains a  dense open subset  $\eO_J$ consisting of operators  $S$ such that
 \[
\dim_\bC\ker(\rho_j-S)= 1,\;\;\forall j\in J.
\]
For $\nu\in \bI^+_n$  we set $\bsF_\nu:=\Sp\{\be_i,\;\;i\leq \nu\}$.

Observe that   if $S\in \eO_J$, then for every $j\in J$ we have
$\ker(\rho_j-S)\subset \bsF_{j-1}^\perp$. Suppose that
$J=\{j_m<\cdots<j_1\}$, and define
 \[
 \Phi: \eO_J\ra \bP(\bsF_{j_m}^\perp)\times \cdots \bP(\bsF_{j_1-1}^\perp),\;\; S\mapsto \bigl(\,\ker(\rho_{j_m}-S),\dotsc, \ker(\rho_{j_1}-S)\,\bigr).
 \]
 The image    of $\Phi$  is
 \[
 \Phi(\eO_J)=\bigl\{ (\ell_m,\dotsc,\ell_1)\in \bP(\bsF_{j_m}^\perp)\times \cdots \bP(\bsF_{j_1-1}^\perp);\;\; \ell_i\perp \ell_{i'},\;\;\forall i\neq i'\,\bigr\}.
 \]
 The  resulting map  $\eO_J\ra \Phi(\eO_J)$ is a fibration with  fiber  over $(\ell_m,\dotsc,\ell_1)$ diffeomorphic to  the      manifold $\eF$ consisting of the unitary operators on  the subspace  $(\ell_m\oplus \cdots\oplus\ell_1)^\perp\subset \whE^+$  which  do not have the numbers $\rho_j$, $j\in J$, in their spectra.  The manifold $\eF$ is open in  this group of unitary operators. Now observe that
 \[
 \dim_\bR\Phi(\eO_J)= 2(n-j_m)+2(n-1-j_{m-1})+ \cdots +2(n-(m-1)-j_1)
 \]
 \[
 = 2nm -m(m-1)-2\sum_{j\in J} j.
 \]
 The fiber $\eF$ has dimension $(n-m)^2$ so that
 \[
 \dim_\bR\eO_J= (n-m)^2+2nm -m(m-1)-2\sum_{j\in J} j= n^2 -\sum_{j\in J} (2j-1).
 \]
 Hence
 \[
\codim \eY_J= \codim \eO_J=\bw(J).
\]

\smallskip

\noindent {\bf Proof of  B.} In the proof we will need the following
technical result.

\begin{lemma} There exists a   $C^2$-semialgebraic map $\xi:S^1\ra S^1$   semialgebraically homotopic to the  identity such that
\[
\xi^{-1}(1)=\{\rho_i;\;\;i\in I\}.
\]
\label{lemma: deform}
\end{lemma}

\begin{figure}[h]
\centerline{\epsfig{figure=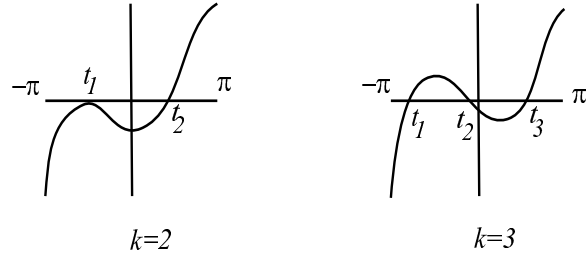,height=1.3in,width=3in}}
\caption{\sl Constructing degree $1$ self maps of the circle.}
\label{fig: local2}
\end{figure}

\proof   We write $ \rho_i=e^{\ii t_i}$, $t_i\in(-\pi,\pi)$ ,and
we consider a  $C^2$-semialgebraic map
\[
f: [-\pi,\pi]\ra [-\pi,\pi]
\]
satisfying the following conditions (see Figure \ref{fig: local2},
where $k=\# I$ )

\begin{itemize}

\item $f(\pm\pi)=\pm\pi$.
\item $f^{-1}(0)=\{t_i;\;\;i\in I\}$.

\end{itemize}

Now define $\xi: S^1\ra S^1$ by setting
\[
\xi(e^{\ii t})=e^{\ii f(t)},\;\; t\in [-\pi,\pi].
\]
We have a $C^2$ semialgebraic  homotopy between $\xi$ and the
identity map given by
\[
\xi_s(e^{\ii t})= e^{\ii( (1-s)f(t)+st)},\;\;s\in [0,1], \;\;t\in
[-\pi,\pi].\proofend
\]

Using the map $\xi$ in Lemma \ref{lemma: deform} we define
\[
\Xi: U(\whE^+)\ra U(\whE^+),\;\;  S\mapsto \Xi(S).
\]
The map $\Xi$ is   semialgebraic because its graph
$\Gamma_\Xi\subset U(\whE^+)\times U(\whE^+)$ can be given the
description
\[
\Gamma_\Xi=\Bigl\{\, (S,S');\;\exists A\in
U(\whE^+),\;ASA^*=\diag(\lambda_1,\dotsc,\lambda_n),\;AS'A^*=\diag(\xi(\lambda_1),\dotsc,\xi(\lambda_n)\,\bigr)\,\Bigr\}.
\]
The continuity of $\Xi$  is classical; see \cite[Theorem X.7.2]{DS}.

If we consider the  set $\eO_I$ defined in the proof of {\bf A} then
we notice that
\[
\Xi(\eO_I)\subset \eW_I^-(\underline{\be},1)
\]
and thus
\[
\Xi(\bigcap_{i\in I}\eX_i(\rho_i)
\,\Bigr)=\Xi\bigl(\,\cl(\eO_I)\,\bigr)\subset\cl(\,\eW_I^-(1)\,)=\eX_I(1).
\]
This proves {\bf B}.

\noindent {\bf Proof of C.}   For $\nu\in \bI_n^+$ and  $\rho\neq
-1$ we set
\[
{}^*\eW_\nu^-(\rho):=\bigl\{\, S\in
\eW_\nu^-;\;\;\dim_\bC\ker(\rho-S)=1,\;\;\ker(1+S)=0\,\bigr\}
\]
Note that ${}^*\eW_\nu^-(\rho)$ is an open and dense subset of
$\eW^-_\nu(\rho)$. We  first want to produce  a natural
trivializing frame of the conormal bundle of ${}^*\eW_\nu^-(\rho)$.
Set $\lambda:=-\ii\frac{1-\rho}{1+\rho}$.

The Cayley transform
\[
S\mapsto -\ii (\one-S)(\one+S)^{-1}
\]
maps ${}^*\eW_\nu^-(\rho)$  onto the  subset $\eR_\nu^*$ of the
space of Hermitian operators $A:\whE^+\ra \whE^+$  such that

\begin{itemize}

\item $\dim \ker (\lambda-A)=1$.

\item $\be_j\perp \ker (\lambda-A)$, $\forall j<\nu$.

\item $(\be_\nu,\bu)\neq 0$, $\forall u\in \ker(\lambda- A)$, $\bu\neq 0$.

\end{itemize}

Note that for any $A\in \eR_\nu^*$ there exists a unique vector
$\bu=\bu_A\in \ker(\lambda- A)$ such that $(\bu,\be_\nu)=1$.  For
$A\in \eR_\nu^*$ we denote by $(\lambda-A)^{[-1]}$ the unique
Hermitian operator $\whE^+\ra \whE^+$ such that
\[
(\lambda-A)^{[-1]}\bu_A=0,\;\;(\lambda- A)^{[-1]}
(\lambda-A)\bv=\bv,\;\;\forall\bv\perp\bu_A.
\]
If $(-\ve, \ve)\ni t\mapsto A_t \in \eR_\nu^*$ is a smooth path, and
$\bu_t:=\bu_{A_t}$,  then differentiating the equality
$A_t\bu_t=\lambda \bu_t$ at $t=0$ we deduce
\[
\dot{A}_0\bu_0=(\lambda-A_0)\dot{\bu}_0.
\]
Taking the inner product with $\bu_0$ we deduce
\[
( \dot{A}_0\bu_0,\bu_0)=0.
\]
We write $\dot{\bu}_0=c \bu_0+\dot{\bv}_0$, where
$(\bv_0,\bu_0)=(\bv_0, \be_j)=0$, $\forall j<\nu$. We deduce
\[
\bv_0= (\lambda-A_0)^{[-1]}\dot{A}_0\bu_0,
\]
so that
\[
\bigl(\, \dot{A}_0\bu_0,
(\lambda-A_0)^{[-1]}\be_j\,\bigr)=(\bv_0,\be_j)=0,\;\;\forall j<\nu.
\]
This shows that the fiber at $A_0$ of the conormal bundle of
$\eR_\nu^*$  contains the  $\bR$-linear forms
\[
\dot{A}_0\mapsto u^\nu(\dot{A}_0) = u^\nu_{A_0}(\dot{A}_0)=\bigl(
\dot{A}_0\bu_0,\bu_0\,\bigr),
\]
\[
\dot{A}_0\mapsto u^{j\nu}(\dot{A}_0)=
u^{j\nu}_{A_0}(\dot{A}_0)=\re\bigl(\, \dot{A}_0\bu_0,(\lambda-
A_0)^{[-1]}\be_j\,\bigr),
\]
\[
\dot{A}_0\mapsto v^{j\nu}(\dot{A}_0)=
v^{j\nu}_{A_0}(\dot{A}_0)=\im\bigl(\, \dot{A}_0\bu_0,
(\lambda-A_0)^{[-1]}\be_j\,\bigr).
\]
Since the vectors $\be_j$, $j<\nu$ lie in the orthogonal complement
of $\ker(\lambda-A_0)$ we deduce that the  vectors
$(\lambda-A_0)^{[-1]}\be_j$, $j<\nu$ are linearly independent over
$\bC$. A dimension count now implies that the above linear forms
form a basis  of the  fiber at $A_0$ of the conormal bundle of
$\eR^*_\nu$. Since  the forms $u^\nu_A$, $u^{j\nu}_A$, $v^{j\nu}_A$
depend smoothly on $A\in \eR^*_\nu$, we deduce that they define a
smooth frame  of the  conormal bundle. Moreover, the canonical
orientation  of $\eR_\nu^*$ is given by
\[
(-1)^{\bw(\nu)} u^\nu\wedge \bigwedge_{j<\nu}u^{j\nu}\wedge
v^{j\nu}.
\]
In particular, we deduce that if $S\in {}^*\eW_\nu^-(\rho)$ is a
unitary operator such that the vectors $\be_i$ are eigenvectors of
$S$ then the canonical orientation of the fiber of $S$ of the
conormal bundle of  $ {}^*\eW_\nu^-(\rho)$ is given by
\[
\theta^\nu\wedge \bigwedge_{j<\nu}\theta^{j\nu}\wedge \vfi^{j\nu}
\]
where, for any $\dot{S}\in T_SU(\whE^+)$ we have
\[
\theta^\nu(\dot{S})=  (-\ii S^{-1}\dot{S}\be_\nu,\be_\nu),
\]
\[
\theta^{j\nu}(\dot{S}=\re  (-\ii S^{-1}\dot{S}\be_\nu,\be_j),
\]
\[
\vfi^{j\nu}(\dot{S})= \im  (-\ii S^{-1}\dot{S}\be_\nu,\be_j).
\]
We deduce that if $\rho: I\ra S^1\setminus\{\pm 1\}$ is an injective
map  then the  manifolds ${}^*\eW_i^-(\rho_i)$, $i\in I$ intersect
transversally.

Now observe that  the manifolds  ${}^*\eW_i^-(\rho_i)$, $i\in I$,
and  $\eW_I^-(1)$ intersect  at a unique point $S_0\in U(\whE^+)$,
where
\[
S_0\be_j=\begin{cases}\rho_j\be_j & j\in I\\
\be_j & j\in I^c.
\end{cases}
\]
The computations in Example  \ref{ex: orient} show  that this
intersection is transversal, and moreover, at $S_0$ we have
\[
\ori_{i_k}^\perp\wedge \cdots\wedge \ori_{i_1}^\perp\wedge
\ori_{I^c}^\perp=\eps(I,I^c)\ori\bigl(\,T^*_{S_0}U(\whE^+)\,\bigr)
\]
The equality  $k_I=1$ now follows by invoking Proposition \ref{prop:
inter}.\qed

\begin{remark} Observe that on the  collection of  subsets   of $\bI_n^+$ we have  two partial order relations.
\[
 K\prec M\Longleftrightarrow \bsJ W_{K^c}^-\cap W_M^-\neq \emptyset\Longleftrightarrow W_K^-\subset \cl(W_M^-).
 \]
 and
 \[
 K\supset  M \Longleftrightarrow  \balph_{K^c}\bullet \balph_M\neq 0.
 \]
Note that $K\supset M \Longrightarrow K\prec M$, but the converse
is not true.    Following the analogy with  the complex
Grassmannian, we could refer to the partial order $\prec$ as  the
\emph{Bruhat order} on the set of parts of $\bI_n^+$.    It would be
very interesting to investigate the combinatorial properties of the
poset $(\bI_n^+,\prec)$. What is its M\"{o}bius function? Is this  a
Cohen-Macaulay poset?

These  combinatorial questions are  special cases of a the more
general   problem concerning the nature of the singularities of
the $AS$ varieties.   We believe that a good geometric
understanding of these singularities  could lead to a more refined
information concerning the homotopy type  of $\Lag_h(E)$. In
particular, we believe that   from the structure of the
singularities we can extract enough information about the Steenrod
squares  to be able to distinguish $\Lag(\whE)$ from the product
of odd dimensional spheres which has the same cohomology ring.\qed
 \end{remark}

\appendix

\section{Tame geometry}
\label{s: a}
\setcounter{equation}{0}
Since the subject of tame  geometry is  not     very familiar to
many  geometers we  devote this section to a   brief  introduction
to this topic.   Unavoidably, we will  have to omit many interesting
details and  contributions, but  we refer to \cite{Co, Dr1, DrMi2}
for  more systematic presentations. For every set $X$ we will denote
by $\eP(X)$ the   collection of all subsets of $X$

An  \emph{$\bR$-structure}\footnote{This is a highly condensed and
special version of  the traditional definition of structure.     The
model theoretic definition  allows for  ordered fields, other than
$\bR$, such as extensions of $\bR$ by ``infinitesimals''. This can
come in handy even if  one is interested  only in the  field $\bR$.}
is a collection $\eS=\bigl\{\, \eS^n\,\bigr\}_{n\geq 1}$,
$\eS^n\subset \eP(\bR^n)$, with the following properties.

\begin{description}

\item[${\bf E}_1$]   $\eS^n$ contains all the real algebraic subvarieties of $\bR^n$, i.e., the zero sets of  finite collections of polynomial in $n$ real variables.

\item[${\bf E}_2$]  For every   linear map $L:\bR^n\ra \bR$, the half-plane $\{\vec{x}\in \bR^n;\;\;L(x)\geq 0\}$  belongs to $\eS^n$.

\item[${\bf P}_1$] For every $n\geq 1$, the family $\eS^n$ is closed under  boolean operations, $\cup$, $\cap$ and complement.

\item[${\bf P}_2$]  If $A\in \eS^m$, and $B\in \eS^n$, then $A\times B\in \eS^{m+n}$.

\item[${\bf P}_3$]  If $A\in \eS^m$, and $T:\bR^m\ra \bR^n$ is an affine map, then $T(A)\in \eS^n$.

\end{description}

\begin{ex}[Semialgebraic sets]    Denote by $\eS_{alg}$ the collection of real semialgebraic sets.  Thus,  $A\in \eS^n_{alg}$ if and only if  $A$  is a finite  union of sets,  each of which is described by finitely many polynomial equalities and inequalities. The celebrated Tarski-Seidenberg theorem states that $\eS_{alg}$ is a structure.\qed
\end{ex}

 Given a structure $\eS$, then an $\eS$-\emph{definable} set is a set that belongs to one of the  $\eS^n$-s. If $A, B$ are $\eS$-definable, then a function $f: A\ra B$ is called $\eS$-\emph{definable} if its graph
 \[
 \Gamma_f :=\bigl\{ (a,b)\in A\times B;\;\;b=f(a)\,\bigr\}
 \]
 is $\eS$-definable. The reason   these sets are called definable has to do with mathematical logic.

Given an $\bR$-structure $\eS$, and a collection $\eA=(\eA_n)_{n\geq
1}$, $\eA_n\subset\eP(\bR^n)$, we can form a new structure
$\eS(\eA)$, which is the smallest structure containing   $\eS$ and
the sets in $\eA_n$. We say that $\eS(\eA)$ is obtained  from  $\eS$
by \emph{adjoining the collection $\eA$}.

\begin{definition} An $\bR$-structure is called \emph{$o$-minimal} (order minimal) or \emph{tame} if  it satisfies the property

\begin{description}
\item[{\bf T}]    Any set $A\in  \eS^1$ is a \emph{finite} union of open intervals $(a,b)$, $-\infty \leq a <b\leq \infty$, and singletons $\{r\}$. \qed
\end{description}

\end{definition}

\begin{ex} (a) The collection  $\eS_{alg}$ of real semialgebraic sets  is a tame structure.

\noindent (b)(Gabrielov-Hironaka-Hardt)   A \emph{restricted} real
analytic function is a  function $f:\bR^n\ra \bR$ with the property
that there exists a real analytic function $\tilde{f}$ defined in an
open  neighborhood $U$ of the cube $C_n:=[-1,1]^n$ such that
\[
f(x)=\begin{cases}
\tilde{f}(x) & x\in C_n\\
0 & x\in \bR^n\setminus C_n.
\end{cases}
\]
we denote by $\eS_{an}$ the structure obtained from $\eS_{alg}$ by
adjoining the  graphs of all the restricted real analytic functions.
Then $\eS_{an}$ is a tame structure, and the $\eS_{an}$-definable
sets are called \emph{globally subanalytic sets}.

\noindent(c)(Wilkie, van den Dries, Macintyre, Marker)   The
structure  obtained  by adjoining to  $\eS_{an}$ the graph of the
exponential function $\bR\ra \bR$, $t\mapsto e^t$,  is a tame
structure.

\noindent(d)(Khovanski-Speissegger)   There   exists  a tame
structure  $\pfc$ with the following properties
\begin{enumerate}
\item[($d_1$)] $\eS_{an} \subset \pfc$

\item[($d_2$)]  If $U\subset \bR^n$ is open, connected  and $\pfc$-definable, $F_1,\dotsc, F_n: U\times \bR\ra \bR$ are $\pfc$-definable and $C^1$, and  $f: U\ra \bR$ is a $C^1$ function satisfying
\begin{equation}
\frac{\pa f}{\pa x_i} = F_i(x, f(x)),\;\;\forall x\in
\bR,,\;\;i=1,\dotsc, n, \label{eq: pfaff}
\end{equation}
then  $f$ is $\pfc$-definable.

\item[($d_3$)] The structure $\pfc$ is the minimal structure satisfying ($d_1$) and ($d_2$).
\end{enumerate}
The  structure  $\pfc$  is called the \emph{pfaffian
closure}\footnote{Our definition  of pfaffian closure   is more
restrictive than the original one in \cite{Kho, Sp1}, but it
suffices for the geometrical applications we have in mind.}
  of $\eS_{an}$.

Observe that if $f: (a,b)\ra \bR$ is $C^1$,  $\pfc$-definable,  and
$x_0\in (a,b)$ then the antiderivative $F:(a, b)\ra \bR$
  \[
  F(x)=\int_{x_0}^x f(t)dt,\;\; x\in (a,b),
  \]
  is also $\pfc$-definable. \qed

\end{ex}

 The  definable sets  and function of a tame structure have  rather remarkable \emph{tame} behavior which prohibits  many pathologies.  It is perhaps   instructive to give an example of function which is not definable in any tame structure. For example, the function $x\mapsto \sin x$ is not definable in a tame structure because the intersection of its graph with the horizontal axis is the  countable set $\pi\bZ$  which  violates  the $o$-minimality condition ${\bf O}$.

 We will list below some of the nice properties of the sets and function definable  in a tame structure  $\eS$. Their proofs can be found in \cite{Co, Dr1}.

\smallskip

\noindent $\bullet$ (\emph{Curve selection.}) If $A$ is an
$\eS$-definable set, and $x\in\cl(A)\setminus A$, then there exists
an $\eS$ definable  continuous map
\[
\gamma:(0,1)\ra A
\]
such that $x=\lim_{t\ra 0} \gamma(t)$.

\noindent $\bullet$ (\emph{Closed graph theorem.})  Suppose $X$ is a
tame set and $f: X\ra \bR^n$ is a  tame bounded function.  Then $f$
is continuous if and only if  its graph is closed in $X\times
\bR^n$.

\noindent $\bullet$   (\emph{Piecewise  smoothness of  tame
functions.})  Suppose $A$ is an $\eS$-definable  set, $p$ is a
positive integer, and $f: A\ra \bR$ is a definable function. Then
$A$ can be partitioned into finitely many  $\eS$ definable sets
$S_1,\dotsc, S_k$,     such that each  $S_i$ is a $C^p$-manifold,
and each of the restrictions $f|_{S_i}$ is a $C^p$-function.

\noindent $\bullet$ (\emph{Triangulability.})  For every   compact
definable set $A$, and any finite collection of definable  subsets
$\{S_1,\dotsc, S_k\}$, there exists  a compact simplicial complex
$K$, and a  definable homeomorphism
\[
\Phi: |K|\ra A
\]
such that  all the sets $\Phi^{-1}(S_i)$ are unions of  relative
interiors of faces of $K$.

\noindent $\bullet$ (\emph{Dimension.})  The  dimension of an
$\eS$-definable  set $A\subset \bR^n$ is the supremum over all the
nonnegative integers $d$ such that there exists a $C^1$  submanifold
of $\bR^n$ of dimension $d$ contained in $A$.  Then $\dim A
<\infty$, and
\[
\dim (\cl(A)\setminus A) <\dim A.
\]

\noindent $\bullet$ (\emph{Crofton formula}, \cite{BK}, \cite[Thm.
2.10.15, 3.2.26]{Feder}.) Suppose $E$ is an Euclidean space, and
denote by $\Graff^k(E)$ the Grassmannian of affine subspaces of
codimension $k$ in $E$.  Fix an invariant measure $\mu$ on
$\Graff^k(E)$.\footnote{ The measure $\mu$ is unique up to a
multiplicative constant.} Denote by $\eH^k$ the $k$-dimensional
Hausdorff measure. Then there exists a constant $C>0$, depending
only on $\mu$, such that for every compact, $k$-dimensional  tame
subset $S\subset E$ we have
\[
\eH^k(S)= C\int_{\Graff^k(E)} \chi(L\cap S) d\mu(L).
\]

\noindent $\bullet$ (\emph{Finite volume.}) Any compact
$k$-dimensional tame set has finite $k$-dimensional Hausdorff
measure.

\qed

In the remainder of this section, by a  tame set we will understand
a $\pfc$-definable set.

\begin{definition}  A tame flow  on a tame set $X$ is a   topological flow $\Phi:\bR\times X\ra X$, $(t,x)\mapsto \Phi_t(x)$, such that the  map $\Phi$ is $\pfc$-definable.\qed
\end{definition}

We list below  a  few properties of tame flows. For proofs we  refer
to \cite{Ni2}.

\begin{proposition} Suppose $\Phi$ is a tame flow on a compact tame set $X$. Then the following hold.

\smallskip

(a)  The  flow $\Phi$ is a finite volume flow  in the sense of
\cite{HL}.

(b) For every  $x\in X$ the limits $\lim_{t\ra \pm \infty}\Phi_t(x)$
exist and are stationary points of $\Phi$. We denote  them by
$\Phi_{\pm \infty}(x)$.

(c) The  maps $x\mapsto \Phi_{\pm \infty}(x)$ are definable.

(d) For any stationary point $y$ of $\Phi$, the unstable variety
$W_y^-=\Phi_{-\infty}^{-1}(y)$ is a definable subset of $X$.   In
particular, if $k=\dim W_y^-$, then  $W_y^-$ has finite $k$-th
dimensional Hausdorff measure. \qed \label{prop: tamefl}
\end{proposition}

\begin{theorem}[Theorem 4.3, \cite{Ni2}]  Suppose $M$ is a compact, connected,  real analytic, $m$-dimensional manifold, $f: M\ra\bR$ is a real analytic Morse function, and $g$ is a real analytic metric  on $M$ such that in the neighborhood of   each critical point $p$ there exists real analytic coordinates $(x^i)_{1\leq i\leq m}$ and  nonzero real numbers $(\lambda_i)_{1\leq i\leq m}$ such that,
\[
\nabla^g f= \sum_{i=1}^m \lambda_i\pa_{x^i},\;\;\mbox{near $p$}.
\]
Then the   flow generated by the gradient $\nabla^g f$ is a tame
flow.\qed \label{th: tame-flow}
\end{theorem}

\section{Subanalytic  currents}
\label{s: b}
\setcounter{equation}{0}
 In this appendix  we gather  without proofs a few facts  about  the subanalytic currents introduced by R. Hardt in \cite{Hardt2}.  Our  terminology    concerning currents closely  follows that of Federer \cite{Feder} (see also the more accessible \cite{Morg}).  However, we changed some  notations to better resemble notations used in  algebraic topology.

Suppose $X$ is a $C^2$, oriented Riemann manifold of dimension $n$.
We denote by $\Omega_k(X)$ the space of $k$-dimensional currents in
$X$, i.e., the topological dual space of the space
$\Omega^k_{cpt}(X)$ of smooth, compactly supported $k$-forms on $M$.
We will denote by
\[
\lan\bullet,\bullet\ran: \Omega^k_{cpt}(X)\times \Omega_k(X)\ra \bR
\]
the natural pairing.  The boundary of a current $T\in \Omega_k(X)$
is  the $(k-1)$-current defined via the Stokes formula
\[
\lan \alpha, \pa T\ran :=\lan d\alpha, T\ran,\;\;\forall \alpha\in
\Omega^{k-1}_{cpt}(X).
\]
For every  $\alpha\in \Omega^k (M)$,  $T\in \Omega_m(X)$,  $k\leq m$
define $\alpha \cap T\in \Omega_{m-k}(X)$ by
\[
\lan\beta , \alpha \cap T  \ran =\lan \alpha\wedge \beta, T
\ran,\;\;\forall \beta\in \Omega^{n-m+k}_{cpt}(X).
\]
We have
\[
\lan \beta, \pa (\alpha \cap T)\ran = \lan \,d\beta, (\alpha\cap
T),\ran = \lan \alpha\wedge d\beta, T\ran
\]
\[
=(-1)^k \lan d(\alpha\wedge \beta) -d\alpha\wedge \beta, T\ran =
(-1)^k \lan \beta, \alpha \cap\pa T\ran +(-1)^{k+1} \lan \beta,
d\alpha \cap T\ran
\]
which yields the \emph{homotopy formula}
\begin{equation}
\pa (\alpha\cap T)= (-1)^{\deg \alpha} \bigl(\, \alpha \cap \pa
T-(d\alpha) \cap T\,\bigr). \label{eq: homotop}
\end{equation}

We say that a set  $S\subset \bR^n$ is \emph{locally subanalytic}
if for any $p\in \bR^n$ we can find an open ball $B$ centered at
$p$ such that $B\cap S$ is globally subanalytic.

\begin{remark} There is a rather subtle distinction between globally subanalytic and locally subanalytic sets. For example, the graph of the function $y=\sin(x)$ is a locally subanalytic subset of $\bR^2$, but it is not a globally subanalytic  set. Note that a compact, locally subanalytic set is globally subanalytic.\qed
\end{remark}

If $S\subset \bR^n$ is an orientable, locally subanalytic, $C^1$
submanifold of $\bR^n$ of dimension $k$,  then any orientation
$\ori_S$ on $S$ determines a  $k$-dimensional current $[S,\ori_S]$
via the equality
\[
\lan \alpha, [S, \ori_S]\ran:=\int_S \alpha,\;\;\forall \alpha\in
\Omega^k_{cpt}(\bR^n).
\]
The integral in the right-hand side is well defined because any
bounded, $k$-dimensional  globally subanalytic set has finite
$k$-dimensional Hausdorff  measure.  For any open, locally
subanalytic   subset $U\subset \bR^n$ we  denote by $[S,\ori_S]\cap
U$  the   current $[S\cap U, \ori_S]$.

For any  locally subanalytic subset $X\subset \bR^n$ we denote by
$\eC_k(X)$ the    Abelian subgroup of $\Omega_k(\bR^n)$
generated  by currents of the form $[S,\ori_S]$,  as above, where
$\cl(S)\subset X$. The above operation $[S,\ori_S]\cap U$, $U$ open
subanalytic extends to a morphism of  Abelian groups
\[
\eC_k(X)\ni T\mapsto T\cap U\in\eC_k(X\cap U).
\]
We will refer to the elements of $\eC_k(X)$ as \emph{subanalytic
(integral) $k$-chains} in $X$.

Given compact subanalytic sets $A\subset X\subset \bR^n$ we set
\[
\eZ_k(X,A)=\bigl\{ T\in \eC_k(\bR^n);\;\;\supp T\subset X,\;\;\supp
\pa T\subset A\,\bigr\},
\]
and
\[
\eB_k(X, A)=\bigl\{ \pa T + S;\;\;T\in \eZ_{k+1}(X,A)), \;\;S\in
\eZ_k(A)\,\bigr\}.
\]
We set
\[
\eH_k(X,A):=\eZ_k(X,A)/\eB_k(X,A).
\]
R. Hardt has proved in \cite{Hardt1, Hardt2} that the assignment
\[
(X,A)\longmapsto \eH_\bullet(X,A)
\]
satisfies the Eilenberg-Steenrod   homology axioms with
$\bZ$-coefficients from which we conclude  that $\eH_\bullet(X,A)$
is naturally isomorphic  with the integral homology  of the pair.
In fact, we can be much more precise.

If $X$ is a compact subanalytic we can form the chain complex
\[
\cdots \stackrel{\pa}{\ra}\eC_k(X)\stackrel{\pa}{\ra}
\eC_{k-1}(X)\stackrel{\pa}{\ra}\cdots.
\]
whose homology is $\eH_\bullet(X)$.

If we choose  a subanalytic triangulation  $\Phi: |K|\ra\ra X$, and
we linearly orient the vertex set $V=V(K)$, then  for any
$k$-simplex $\si\subset K$ we get a subanalytic map   from the
standard  affine $k$-simplex $\Delta_k$  to $X$
\[
\Phi^\si:\Delta_k\ra X.
\]
This defines a current $[\si]=\Phi^\si_*([\Delta_k])\in \eC_k(X)$.
By linearity  we obtain a    morphism from the group of   simplicial
chains  $C_\bullet(K)$ to $\eC_\bullet(X)$ which commutes with the
respective boundary operators. In other words, we obtain  a morphism
of chain complexes
\[
C_\bullet(K)\ra \eC_\bullet(\Phi|K|).
\]
The arguments in \cite[Chap.III]{ES} imply that this induces  an
isomorphism in homology.

To describe the intersection theory of subanalytic chains  we need
to recall a  fundamental result of R. Hardt, \cite[Theorem
4.3]{Hardt}.  Suppose $E_0, E_1$ are two oriented  real Euclidean
spaces of dimensions $n_0$ and respectively $n_1$,  $f:E_0\ra E_1$
is a real analytic map, and $T\in \eC_{n_0-c}(E_0)$ a subanalytic
current of codimension $c$.  If $y$ is a regular value of $f$,  then
the fiber $f^{-1}(y)$  is  a submanifold  equipped with a natural
coorientation and thus defines a subanalytic  current $[f^{-1}(y)]$
in $E_0$  of codimension $n_1$, i.e., $[f^{-1}(y)]]\in
\eC_{d_0-d_1}(E_0)$. We would like to define the intersection of $T$
and $[f^{-1}(y)]$  as a subanalytic current  $T\bullet
[f^{-1}(y)]\in \eC_{n_0-c-n_1}(E_0)$.     It  turns out that  this
is possibly quite often, even in cases when $y$ is  not a regular
value.

\begin{theorem}[Slicing Theorem]  Let $E_0$, $E_1$, $T$ and $f$ be  as above, denote by $dV_1$ the Euclidean volume form on $E_1$,   by $\bom_{n_1}$ the volume of the unit ball in $E_1$, and set
\[
\eR_f(T):=\bigl\{ y\in E_1;\;\;\codim (\supp T )\cap f^{-1}(y) \geq
c+ n_1,\;\;\codim (\supp \pa T )\cap f^{-1}(y)\geq
c+n_1+1\,\bigr\}.
\]
For  every $\ve>0$ and $y\in E_1$ we define $T\bullet_\ve
f^{-1}(y)\in \Omega_{n_0-c-n_1}(E_0)$ by
\[
\bigl\lan\, \alpha, T\bullet_\ve f^{-1}(y)\,\bigr\ran
=\frac{1}{\bom_{n_1}\ve^{n_1}}\bigl\lan\, (f^*dV_1)\wedge\alpha ,
T\cap \bigl(\,f^{-1}(B_\ve(y)\,\bigr)\,\bigr\ran,\;\;\forall
\alpha\in \Omega^{n_0-c-n_1}_{cpt}(E_0).
\]
Then  for every $y\in \eR_f(T)$, the currents $T\bullet_\ve
f^{-1}(y)$ converge weakly as $\ve>0$   to a  subanalytic  current
$T\bullet f^{-1}(y)\in \eC_{n_0-c-n_1}(E_0)$ called  the
\emph{$f$-slice} of $T$ over $y$, i.e.,
\[
\bigl\lan\, \alpha, T\bullet f^{-1}(y)\,\bigr\ran =\lim_{\ve\searrow
0}\frac{1}{\bom_{n_1}\ve^{n_1}}\bigl\lan\, (f^*dV_1)\wedge\alpha ,
T\cap \bigl(\,f^{-1}(B_\ve(y)\,\bigr)\,\bigr\ran,\;\;\forall
\alpha\in \Omega^{n_0-c-n_1}_{cpt}(E_0).
\]
Moreover,   the map
\[
\eR_f\ni y\mapsto T\bullet f^{-1}(y)\in \eC_{d_0-c-d_1}(\bR^n)
\]
is continuous in the   locally   flat  topology.\qed

\end{theorem}

We will refer to the points $y\in \eR_f(T)$ as the
\emph{quasi-regular values of $f$ relative to $T$}.

Consider an oriented real analytic manifold $M$ of dimension $m$,
and $T_i\in\eC_{m-c_i}(M)$, $i=0,1$. We would like to define an
intersection  current $T_0\bullet T_1\in \eC_{m-c_0-c_1}(M)$.  This
will require  some very mild  transversality    conditions.

The slicing theorem  describes this intersection current when $T_1$
is     the integration current defined by the fiber of a real
analytic map. We want to reduce the general situation to this  case.
We will achieve this in two steps.

\begin{itemize}

\item  Reduction to the diagonal.

\item Localization.
\end{itemize}

To understand the reduction to the diagonal let us observe that if
$T_0, T_1$ were homology classes  then their intersection
$T_0\bullet T_1$   satisfies the  identity
\[
j_*(T_0\bullet T_1)=(-1)^{c_0(m-c_1)} (T_0\times T_1)\bullet
\Delta_M,
\]
where $\Delta_M$ denotes the diagonal class in $M\times M$, and $j:
M\ra M\times M$ denotes the diagonal embedding.

We use this fact to define the intersection current  in the special
case when $M$ is an open subset of $\bR^m$. In this case the
diagonal $\Delta_M$ is the  fiber over $0$ of the difference map
\[
\delta: M\times M\ra \bR^m,\;\; \delta(m_0,m_1)=m_0-m_1.
\]
If   the currents $T_0, T_1$   are \emph{quasi-transversal}, i.e.,
\begin{subequations}
\begin{equation}
\codim  (\supp T_0)\cap (\supp T_1)\geq c_0+c_1, \label{eq: qt0}
\end{equation}
\begin{equation}
\codim\bigl(\, (\supp T_0\cap \supp\pa T_1)\,\cup\,(\supp\pa T_0\cap
\supp T_1)\,\bigr)\geq c_0+c_1+1, \label{eq: qt1}
\end{equation}
\end{subequations}
then  $0\in \bR^m$ is a $T_0\times T_1$-quasiregular value of
$\delta$ so that  the intersection
\[
(T_0\times T_1)\bullet \delta^{-1}(0)= (T_0\times T_1)\bullet
\Delta_M
\]
is well defined.

The intersection current $T_0\bullet T_1$ is then the  unique
current in $M$ such that
\[
j_*(T_0\bullet T_1)=(-1)^{c_0(m-c_1)}(T_0\times T_1)\bullet
\delta^{-1}(0).
\]
If $M$ is an arbitrary  real analytic manifold  and the  subanalytic
currents are  quasi-transversal then we define $T_0\bullet T_1$ to
be the unique   subanalytic current such that for any  open  subset
$U$ of $M$  real analytically diffeomorphic to an open ball in
$\bR^m$ we have
\[
(T_0\bullet T_1)\cap U= (T_0\cap U)\bullet (T_1\cap U).
\]
One can prove that
\begin{equation}
\pa(T_0\bullet T_1)=(-1)^{c_0+c_1}(\pa T_0)\bullet \pa T_1 +
T_0\bullet(\pa T_1), \label{eq: pa}
\end{equation}
whenever    the  various pairs of  chains in the above formula are
quasi-transversal.

One of the key results in \cite{Hardt1, Hardt2}     states that this
intersection  of quasi-transversal chains induces a  well defined
intersection pairing
\[
\bullet: \eH_{m-c_0}(M)\times \eH_{m-c_1}(M)\ra \eH_{m-c_0-c_1}(M).
\]
These intersections pairings coincide with the intersection pairings
defined via Poincar\'e duality.  This follows by combining two
facts.

\begin{itemize}

\item The  subanalytic homology groups can be computed   via a triagulation, as explained above.

\item The classical proof\footnote{I cannot help but remark  that  it is  hard to find a  book written during the past three decades which presents  a \emph{complete} proof  of the Poincar\'e duality the ``old fashion way'', via triangulations and their dual cell decompositions. }  of the Poincar\'e duality  via triangulations (see \cite[Chap. 5]{Maun}).

\end{itemize}

For a submanifold $S\subset M$  of dimension $k$ we  define the
conormal bundle $T_S^*M$ to be the  kernel of the natural  bundle
morphism
\[
i^*:T^*M|_S\ra T^*S,
\]
where $i:S\hra M$ is the inclusion map. A co-orientation of $S$  is
then an orientation of  the conormal bundle. This induces an
orientation on  the cotangent bundle  of $S$ as follows.

\begin{itemize}

\item Fix $s_0\in S$, and   a positively  basis $\underline{b}_0=\{e^1,\dotsc, e^k\}$ of the   fiber of $T^*SM$ over $s_0$.

\item Extent the basis $\underline{b}_0$ to a positively oriented basis $\underline{b}=\{e^1,\dotsc, e^n\}$ of $T^*_{s_0}M$.

\item Orient $T_{s_0}^*S$ using the ordered basis $\{i^*(e^{k+1}),\dotsc, i^*(e^m)\}$.

\end{itemize}

We see  that  a pair $(S,\ori^\perp)$ consisting of a $C^1$, locally
subanalytic submanifold $S\hra M$, and a co-orientation $\ori^\perp$
defines a subanalytic chain $[S,\ori^\perp]\in\eC_k(M)$. Observe
that
\[
\supp \pa[S,\ori^\perp]\subset \cl(S)\setminus S.
\]
Thus, if $\dim\bigl(\,  \cl(S)\setminus S\,\bigr) < \dim S -1$ then
$\pa[S,\ori^\perp]=0$.

\begin{definition} An \emph{elementary cycle}  of $M$ is a co-oriented  locally subanalytic submanifold $(S,\ori^\perp)$ such that  $\pa[S,\ori^\perp]=0$.

We say that two elementary cycles $(S_i,\ori_i^\perp)$,  $i=0,1$
intersect \emph{quasi-transversally} if  the following hold.

\begin{itemize}

\item   The submanifolds $S_0$, $S_1$ intersect  transversally.

\item   $\cl(S_0)\cap \cl(S_1)=\cl(S_0\cap S_1)$.

\end{itemize}\qed
\label{def: elcyc}
\end{definition}

Observe that if two elementary cycles $(S_i,\ori_i^\perp)$, $i=0,1$,
intersect transversally , then the associated  subanalytic chains
$[S_i,\ori_i^\perp]$ are quasi-transversal.    The conormal bundle
of $S_0\cap S_1$ is the direct sum of the restrictions of the
conormal bundles of $S_0$ and $S_1$,
\[
T^*_{S_0\cap S_1}M=\bigl(T^*_{S_0}M\,\bigr)|_{S_0\cap S_1}\oplus
\bigl(T^*_{S_1}M\,\bigr)|_{S_0\cap S_1}
\]
There is  natural induced co-orientation $\ori_0^\perp\wedge
\ori_1^\perp$  on $S_0\cap S_1$ given by the above \emph{ordered}
direct sum.

\begin{proposition}  Suppose $(S_i,\ori_i^\perp)$, $i=0,1$ are elementary cycles intersecting quasi-transversally.   Then
\[
[S_0,\ori_0^\perp]\bullet[S_1,\ori_1^\perp]=[S_0\cap
S_1,\ori_0^\perp\wedge\ori_1^\perp].
\]
\label{prop: inter}
\end{proposition}

\proof From  (\ref{eq: pa}) we deduce that
\[
\pa [S_0,\ori_0^\perp]\bullet[S_1,\ori_1^\perp]=0.
\]
On the other hand,
\[
\supp [S_0,\ori_0^\perp]\bullet[S_1,\ori_1^\perp] \subset
\cl(S_0)\cap \cl(S_1)=\cl(S_0\cap S_1).
\]
Thus, to  find  the intersection current
$[S_0,\ori_0^\perp]\bullet[S_1,\ori_1^\perp]$ it suffices to test it
with  differential forms $\alpha\in\Omega^{c_0+c_1}(M)$ such that
\[
\supp \alpha \cap \cl(S_0)\cap \cl(S_1)\subset S_0\cap S_1.
\]
Via  local coordinates this reduces the problem  to the special case
when $S_0, S_1$ are  co-oriented subspaces of $\bR^n$ intersecting
transversally in which case the result follows by direct computation
from the definition.  We leave the details to the reader. \qed


\end{document}